\documentclass[12pt]{article}
\usepackage{jheppub}
\pdfoutput=1
\usepackage{amsmath,bbm,array,amsfonts,graphicx,wrapfig,lscape,float,mathtools,multirow,longtable,amsthm}
\usepackage[dvipsnames]{xcolor}
\usepackage{stackrel}
\usepackage[all]{xy}
\usepackage{caption}
\usepackage{subcaption}
\usepackage{multirow}
\usepackage[hang,flushmargin]{footmisc}
\usepackage{mathrsfs}
\usepackage{tikz}
\usepackage{bm}
\usepackage{color}
\usepackage{enumerate}
\usetikzlibrary{decorations.pathreplacing}
\usepackage{diagbox}
\numberwithin{equation}{section}
\numberwithin{figure}{section}
\numberwithin{table}{section}
\usepackage{pgfplots}
\pgfplotsset{compat=1.14}
\usepackage{makecell}
\usepackage{lscape}
\usepackage[utf8]{inputenc}
\usepackage{mathtools}

\usepackage[autostyle=true]{csquotes}
\newtheorem{definition}{Definition}[section]
\newtheorem{theorem}{Theorem}[section]

\usepackage[toc,page]{appendix}

\usepackage{tocloft}

\usetikzlibrary{external}
\tikzset{external/system call={pdflatex \tikzexternalcheckshellescape 
		-halt-on-error
		-interaction=batchmode 
		-jobname "\image" "\texsource"
		&& pdftops -eps "\image.pdf"}}
\newcommand{\comment}[1]{}

\newcommand\blfootnote[1]{
  \begingroup
  \renewcommand\thefootnote{}\footnote{#1}
  \addtocounter{footnote}{-1}
  \endgroup
}

	\title{Neurons on Amoebae}
	
	\author[a]{Jiakang Bao,}
	\author[b,a,c,d]{Yang-Hui He,}
	\author[a]{Edward Hirst}
	
	\affiliation[a]{
		Department of Mathematics, City, University of London, EC1V 0HB, UK}
	\affiliation[b]{
		London Institute for Mathematical Sciences, Royal Institution of GB, W1S 4BS, UK}
	\affiliation[c]{
		Merton College, University of Oxford, OX1 4JD, UK}
	\affiliation[d]{
		School of Physics, NanKai University, Tianjin, 300071, P.R. China}
	
	\emailAdd{jiakang.bao@city.ac.uk}
	\emailAdd{hey@maths.ox.ac.uk}
	\emailAdd{edward.hirst@city.ac.uk}
	
	\preprint{
		\begin{flushright}
			LIMS-2021-007
		\end{flushright}
	}
	
	\abstract{We apply methods of machine-learning, such as neural networks, manifold learning and image processing, in order to study 2-dimensional amoebae in algebraic geometry and string theory. With the help of embedding manifold projection, we recover complicated conditions obtained from so-called lopsidedness. For certain cases it could even reach $\sim99\%$ accuracy, in particular for the lopsided amoeba of $F_0$ with positive coefficients which we place primary focus. Using weights and biases, we also find good approximations to determine the genus for an amoeba at lower computational cost. In general, the models could easily predict the genus with over $90\%$ accuracies. With similar techniques, we also investigate the membership problem, and image processing of the amoebae directly.\protect\blfootnote{Contribution to special volume ``Algebraic Geometry and Machine Learning'' in J.~Symbolic Computation, Springer, Hauenstein et al. Ed.}
	}

\begin{document}
	\maketitle

\section{Introduction}\label{intro}
Amoebae of algebraic varieties along with closely related tropical geometry constitutes a lively branch of contemporary mathematics \cite{gelfand2008discriminants,itenberg2009tropical,zbMATH02223067}.
One notable direction is the relation to dimer models, which bring together the combinatorics of toric varieties, the geometry of complex manifolds, the theory of partitions and statistical mechanics \cite{zbMATH05051319,Kenyon:2003ui}.
In string theory and supersymmetric quantum field theory, amoebae and dimers naturally encode quiver gauge theories whose space of vacua are affine toric Calabi-Yau manifolds \cite{Feng:2005gw} as well as topological strings on these geometries \cite{Heckman:2006sk,Ooguri:2008yb,Cirafici:2009ga,Zahabi:2020hwu}.
The reader is referred to 
\cite{zbMATH02223067,Yamazaki:2008bt,Yamazaki:2011wy,bogaardintroduction,He:2016fnb,maclagan2015introduction,Bao:2020sqg} for various reviews, directed to mathematicians and physicists alike.

More recently, there has been much activity on using techniques of modern data science, especially machine learning (ML), to analyse pure mathematical data, as tools of conjecture formulation and aids in computation and proof.
Whilst this arose from the investigation of the string theory/algebraic geometry landscape \cite{He:2017aed,Krefl:2017yox,Ruehle:2017mzq,Carifio:2017bov,He:2020eva,Bao:2021auj}, one is led to see how such a methodology can be extended to different disciplines of mathematics (q.v.~reviews in \cite{He:2018jtw,He:2021oav}).

Given these two skeins of enquiry, 
it is only natural to ask how one might machine-learn pertinent features of amoeba geometry.
As outlined in \cite{He:2021oav}, whatever sub-field of mathematics, there often emerges theoretical data representable as tensors, labelled or unlabelled. This is particularly pronounced due to the recent vast growth in computing power and experimental mathematics.
One is thence encouraged to immediately turn to a neural network classifier or regressor to treat such data.
Algebraic geometry is certainly a discipline rich in data of this type,
for example varieties can be computationally recorded by lists of its polynomial coefficients; and its morphism maps in cohomology are integer matrices.

As we will see, amoebae and tropical geometry are endowed with the perfect data structure for machine learning: they are, by nature, graphical, and image-processing is ML's forte.
The purpose of this paper is to initiate the study of ML of tropical geometry, much like how \cite{He:2017aed} attempted to learn complex algebraic geometry or how \cite{He:2020fdg} discrete geometry.
We will see how otherwise difficult problems, such as deciding membership and computing genus, can be solved to high accuracy with neural networks (NNs).
In fact, the problem can be effectively redescribed as the identification of bounding hypersurfaces in the amoebae coefficient abstract space --- a feat NNs famously perform exceptionally well at.
We should emphasize what we are doing here is paradigmatically different from \cite{zhang2018tropical,charisopoulos2019tropical}, which finds tropical-geometric formulations of deep ReLU networks.
Our emphasis is the use on ML to study tropical geometry, as opposed to the use of tropical geometry to study ML; and here we focus in on learning the genus of 2-dimensional amoebae within tropical geometry.
It is worth noting that this work puts particular focus on the prototypical $F_0$ example, with other examples considered briefly. Training primarily uses lopsided amoebae, which act as approximations to the true amoebae, where the limit of convergence to the true amoebae considered in specific cases.

\paragraph{Outline} The paper is organized as follows. In \S\ref{rudiments}, we briefly introduce amoebae, lopsidedness, and relevant concepts in mathematics and physics. In \S\ref{ml}, we move onto the ML investigations. Starting with the simplest lopsided amoebae in \S\ref{n1} where NNs are applied to learning the genus for various examples. From manifold learning and the neural network structure, we will try to recover or approximate the conditions that determine the genus. We then conduct similar study for more subtle lopsidedness in \S\ref{largern}. In particular, we will focus on the zeroth Hirzebruch surface $F_0$ as an illustration. In \S\ref{membership}, we will briefly comment on the membership problem. Finally, in \S\ref{ImageCNN} convolutional neural networks (CNNs) learn to identify the genus for $F_0$ from Monte Carlo generated amoebae at varying image resolutions. More background and related calculations for amoebae can be found in the appendices.
One thing to note is that in Monte Carlo, sampling of points is an important issue, which we discuss in Appendix \ref{ap:transformation}.

This paper is a contribution to the special volume ``Algebraic Geometry and Machine Learning'' in the Journal of Symbolic Computation \cite{JSC}, Springer 2021.

\section{Rudiments of Amoebae}\label{rudiments}
In this section, we review some basics of amoebae, with a view toward how they arise in physics.
Suppose we are given a multi-complex-variable polynomial
\begin{equation}\label{poly}
    P(\bm{z})=\sum_{\bm{p}}c_{\bm{p}}\bm{z}^{\bm{p}}=\sum_{i_1,\dots,i_r}c_{i_1,\dots,i_r}z_1^{i_1}\dots z_r^{i_r} \ , \quad
    \bm{z}\in(\mathbb{C}^*)^r \ .
\end{equation}
As is customary, we use bold-face to denote the multi-index notation for monomials.
This defines an affine algebraic variety as a hyper-surface $P(\bm{z}) = 0$ in $(\mathbb{C}^*)^r$.
Then, we have the following. 
\begin{definition}
    Writing the coordinates as
    \begin{equation}
        (z_1,\dots,z_r):=
        (\exp(s_1+i\theta_1),\dots,\exp(s_r+i\theta_r)),
    \end{equation}
    where $s_i\in\mathbb{R}$ and $\theta_i\in[0,2\pi)$,
    the {\bf amoeba} $\mathcal{A}_P$ is the image under the projection
    \begin{equation}
        \textup{proj} : (z_1,\dots,z_r)
        \mapsto
        (s_1,\dots,s_r) \ .
    \end{equation}
    Equivalently, \textup{proj} is the map
    \begin{equation}
        (z_1,\dots,z_r)
        \mapsto
        (\log|z_1|,\dots,\log|z_r|).
    \end{equation}\label{ambproj}
\end{definition}\label{ambdef}
We would often abbreviate this as \textup{Log}$|\bm{z}|$ in our multi-index notation $\bm{z}$.

While the above definition is general, we will exclusively focus on $P(\bm{z})$ being a Newton polynomial, which we recall as follows:
\begin{definition}
    Let $\Delta$ be a convex lattice polytope, which we will also refer to as a {\bf toric diagram} \footnote{We refer the reader to the modern classic \cite{cox2011toric} for toric varieties and convex polytopes.}, with lattice points $\{ \bm{p}  \in \Delta \}$.
    We can always write a polynomial, called the {\bf Newton polynomial}, each of whose monomial terms is $\bm{z}^{\bm{p}}$.
    The coefficient is an arbitrary non-zero complex number and constitutes one of the {\it moduli}.
    Conversely, for any polynomial in \eqref{poly}, one can construct a Newton polytope by forming the convex hull of all the powers which appear in the monomials.
\end{definition}

\begin{figure}[h]
	\centering
		(a)
		\includegraphics[width=4cm]{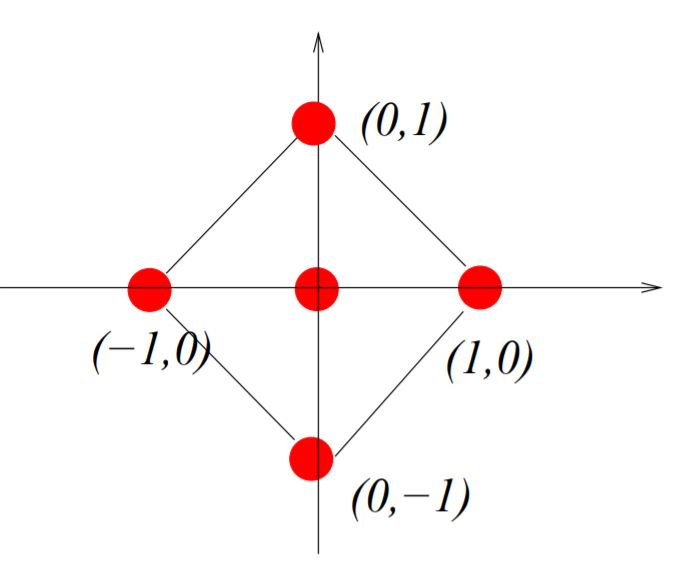}
		(b)
		\includegraphics[width=10cm]{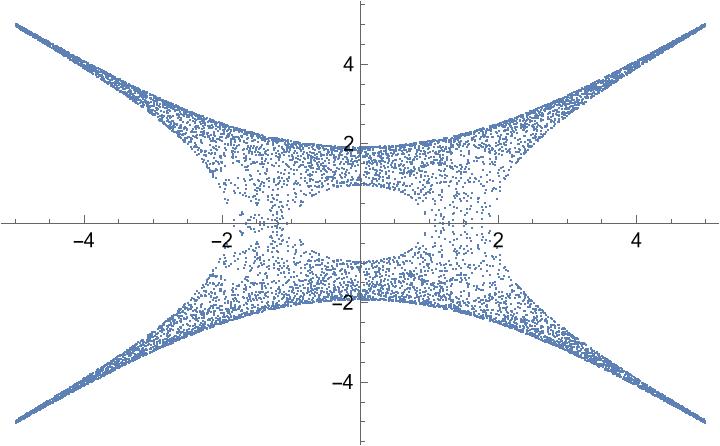}
	\caption{
	(a) A lattice polygon (polytope in 2-D), giving the Newton polynomial  $P(z_1,z_2) = c_{(0,0)} + c_{(1,0)} z_1 + c_{(-1,0)} z_1^{-1} + c_{(0,1)} z_2 + c_{(0,-1)} z_2^{-1}$.
	(b) Taking the coefficients $c_{(0,0)} = 5$ and all others to be 1, the amoeba is the region in blue (the figure is obtained by Monte Carlo).
	}\label{egAmoeba}
\end{figure}
Let us consider an example.
Take the 2-dimensional polytope in Fig.~\ref{egAmoeba} (a).
The lattice points contained therein are $(0,0),\ (1,0),\ (0,1),\ (-1,0),\ (0,-1)$.
Thus, the Newton polynomial is
$P(z_1,z_2) = c_{(0,0)} + c_{(1,0)} z_1 + c_{(-1,0)} z_1^{-1} + c_{(0,1)} z_2 + c_{(0,-1)} z_2^{-1}$, where the five $c$'s are complex coefficients.
Note that strictly speaking, the Newton polynomial is a Laurent polynomials since negative powers appear; these can be removed by shifting the origin of the polytope.
The amoeba for the choice of coefficients $5+ z_1 + z_1^{-1} + z_2 + z_2^{-1}$ is shown in Fig.~\ref{egAmoeba} (b).
This star-shaped region with a hole looks like an amoeba from biology, hence the name. 
The hole in the middle, along with the four unbounded complementary regions are real points $(s_1,s_2)$, where no complex numbers $(z_1,z_2)$ satisfying $P(z_1,z_2)=0$ map to them under the projection \ref{ambproj}. One may also notice there are four \emph{tentacles} of this amoeba that go to infinity exponentially; these four asymptotic lines of the tentacles are known as the \emph{spines}.

Since the amoeba is obtained from the solutions to $P(\bm{z})=0$, it is a subset in $\mathbb{R}^r$.
Thus, for one thing, this half-reduction from complex dimension $r$ to real dimension $r$ is good for {\it visualization}.
In particular, for $r=2$ (where the Newton polytope is a convex lattice polygon), the amoeba is a great way to picture Riemann surfaces. We will exclusively work with $r=2$.
For these, we have some important definitions.
\begin{definition}
Consider the amoeba $\mathcal{A}$ of the complex curve/Riemann surface \newline$P(z_1,z_2)=0$.
\begin{itemize}
    \item There are (potentially separated) regions complement to $\mathcal{A}$, which may or may not be bounded. We shall call the number of such bounded complementary regions the {\bf genus}\footnote{It is curiously appealing that the phrase ``the genus of an amoeba'' makes sense both in biology and in mathematics!
    Perhaps the most common species, Amoeba proteus, should be associated to a particular class of Riemann surfaces.
    } of $\mathcal{A}$.
    For instance, the amoeba in Fig.~\ref{egAmoeba} (b) is genus 1.
    
    \item The genus counts the number of ``holes'' in $\mathcal{A}$, and it explicitly depends on the choice of the coefficients in the Newton polynomial.
    As we will focus on spectral curves $P(z_1,z_2)=0$ with real coefficients, these curves are real plane curves in $\mathbb{RP}^2$. In particular, the curves for which the numbers of holes (as well as components) are maximal amongst the choices of coefficients are called the {\bf Harnack curves}.
\end{itemize}
\end{definition}

\paragraph{Membership Problem: }
Given an amoeba, it is natural to ask how one can determine whether a point in the Log space belongs to the amoeba. 
This is known as the \emph{membership problem}. From \cite{collins1975quantifier,ben1984complexity,theobald2002computing}, we learn that the membership problem can be solved in polynomial time for a fixed dimension. However, to determine the full boundaries of a general amoeba, the typical approach would be approximations using lopsidedness as we shall explain now. The algorithm is more time-consuming and the fastest algorithm so far to our best knowledge was proposed in \cite{forsgard2017lopsided} (see Footnote \ref{cycresfootnote}).

\subsection{Amoebae and Lopsidedness} To learn the complementary regions of amoebae --- and hence its boundaries, genus, and membership decisions --- it is useful to introduce the concept of lopsidedness \cite{purbhoo2006nullstellensatz}.
\begin{definition}
	Let $f\in\mathbb{C}\left[z_1,z_1^{-1},\dots,z_r,z_r^{-1}\right]$ be a sum of (Laurent) monomials $m_i$ as  $f(\bm{z})=m_1(\bm{z})+\dots+m_k(\bm{z})$. 
	For $\bm{x}\in\mathbb{R}^r$, define the list of positive numbers
	\begin{equation}\label{lopPoly}
		f\{\bm{x}\}:=\left\{\left|m_1(\textup{Log}^{-1}(\bm{x}))\right|,\dots,\left|m_k(\textup{Log}^{-1}(\bm{x}))\right|\right\}.
	\end{equation}
	We say a list of positive numbers is {\bf lopsided} if one of the numbers is greater than the sum of all the others.
	This definition can then be applied to $f\{\bm{x}\}$.
\end{definition}
Importantly, a list of positive numbers $\{b_i\}$ is {\it not lopsided} if one could find a list of phases $\{\phi_i\} \subset \mathbb{C}$ with $|\phi_i|=1 \;\forall i$ such that $\sum_i \phi_i b_i = 0$; this follows from the triangular inequality. These exact phases allow a real point on the plane projected to, to have a $\text{Log}^{-1}$ lift into complex space that satisfies the amoebae equation $P(\bm{z})=0$, thus making this real point a member of the amoeba. It is here how not-lopsidedness connects to amoebae membership.

It is then natural to define the following \cite{forsgard2017lopsided}:
\begin{definition}
    Given a Newton polynomial $P$, the {\bf lopsided amoeba} \footnote{The appellation might seem a bit confusing since the lopsided amoeba is the set where the $P\{\bm{x}\}$ is {\it not} lopsided.} is
\begin{equation}
    \mathcal{LA}_P :=
    \{\bm{x}\in\mathbb{R}^r
    |~
    P\{\bm{x}\}~\textup{is not lopsided}\}.
\end{equation}
\end{definition}
For some cases, such as $P=z_1+z_2+1$, $\mathcal{LA}_P=\mathcal{A}_P$. However, in general, they do not need to coincide. Nevertheless,
\begin{equation}
    \mathcal{A}_P \subseteq \mathcal{LA}_P \ ,
\end{equation}
so that $\mathcal{LA}_P$ can be constructed as a crude approximation to $\mathcal{A}_P$.  This can be made precise as follows.

Let $n$ be a positive integer,
$\bm{x}\in\mathbb{R}^r$,
and $P(\bm{x})$ a (Newton) polynomial, define $\Tilde{P}_n$ to be \footnote{In \cite{forsgard2017lopsided}, a faster algorithm was proposed to compute $\mathcal{LA}_{\Tilde{P}_n}$ at level $k$ where $n=2^k$ using the properties of cyclic resultants. The time complexity is $\mathcal{O}(kd^2)$ with $d$ being the degree of $P(z_1,z_2)$.\label{cycresfootnote}}
\begin{equation}
    \Tilde{P}_n(\bm{x}):=
    \prod_{k_1=0}^{n-1}\cdots\prod_{k_r=0}^{n-1}P\left(\text{e}^{2\pi ik_1/n}x_1,\dots,\text{e}^{2\pi ik_r/n}x_r\right).
\end{equation}
Clearly, $\Tilde{P}_1=P$.
Such $\Tilde{P}_n$ is in fact a \emph{cyclic resultant}
\begin{equation}
    \Tilde{P}_n=\text{res}_{u_r}\left(\text{res}_{u_{r-1}}\left(\dots\text{res}_{u_1}\left(P(u_1x_1,\dots,u_rx_r),u_1^n-1\right)\dots,u_{r-1}^n-1\right),u_r^n-1\right)
\end{equation}
where $\text{res}_{u}(f,g)$ is the resultant of $f,g$ with respect to the variable $u$.

The lopsided amoeba
$\mathcal{LA}_{\Tilde{P}_n}$ 
for $\Tilde{P}_n$ approximates
$\mathcal{A}_P$ itself \cite{purbhoo2006nullstellensatz}:
\begin{theorem}
    For an $r$-dimensional Newton polytope $\Delta(P)$, with polytope coordinates $p_i$ for each $i^\text{th}$ direction in the $\mathbb{Z}^r$ lattice which the polytope is defined in, one defines $c_i:=\max(p_i)-\min(p_i)$ over the polytope vertices; then $c=\max(c_i)$. Suppose $\bm{x}\in\mathbb{R}^r\setminus\mathcal{A}_P$ is a point in the amoeba complement whose distance from $\mathcal{A}_P$ is at least $\epsilon>0$.
    If $n$ is large enough so that
    \begin{equation}
        n\epsilon\geq(r-1)\log n+\log((r+3)2^{r+1}c),
    \end{equation}
    then $\Tilde{P}_n\{\bm{x}\}$ is lopsided\footnote{In this paper, as $r$ is always 2, we have $n\epsilon\geq\log n+\log(8c)$.}
    and
    $\mathcal{LA}_{\Tilde{P}_n}$ converges uniformly to $\mathcal{A}_P$ as $n\rightarrow\infty$.
\end{theorem}
A consequence of this is one way to solve the membership problem:
\begin{theorem}
	Let $I\subset\mathbb{C}\left[z_1,z_1^{-1},\dots,z_r,z_r^{-1}\right]$ be an ideal. The point $\bm{x}\in\mathbb{R}^r$ is in the amoeba $\mathcal{A}_I$ if and only if $f\{\bm{x}\}$ is not lopsided for every $f\in I$.
\end{theorem}
Therefore, to fully determine the boundary of an amoeba, we need to consider all the Laurent polynomials in the ideal generated by our Newton polynomial. Equivalently, we need to take $n\rightarrow\infty$ for the cyclic resultant. As a result, we often approximate the boundary with some finite large $n$ in practice. This would also be our basic strategy to study the genus using neural networks later, although unlike finding the boundary, sometimes we can count the genus in other ways, as we will see.

\begin{figure}[h]
	\centering
	\includegraphics[width=10cm]{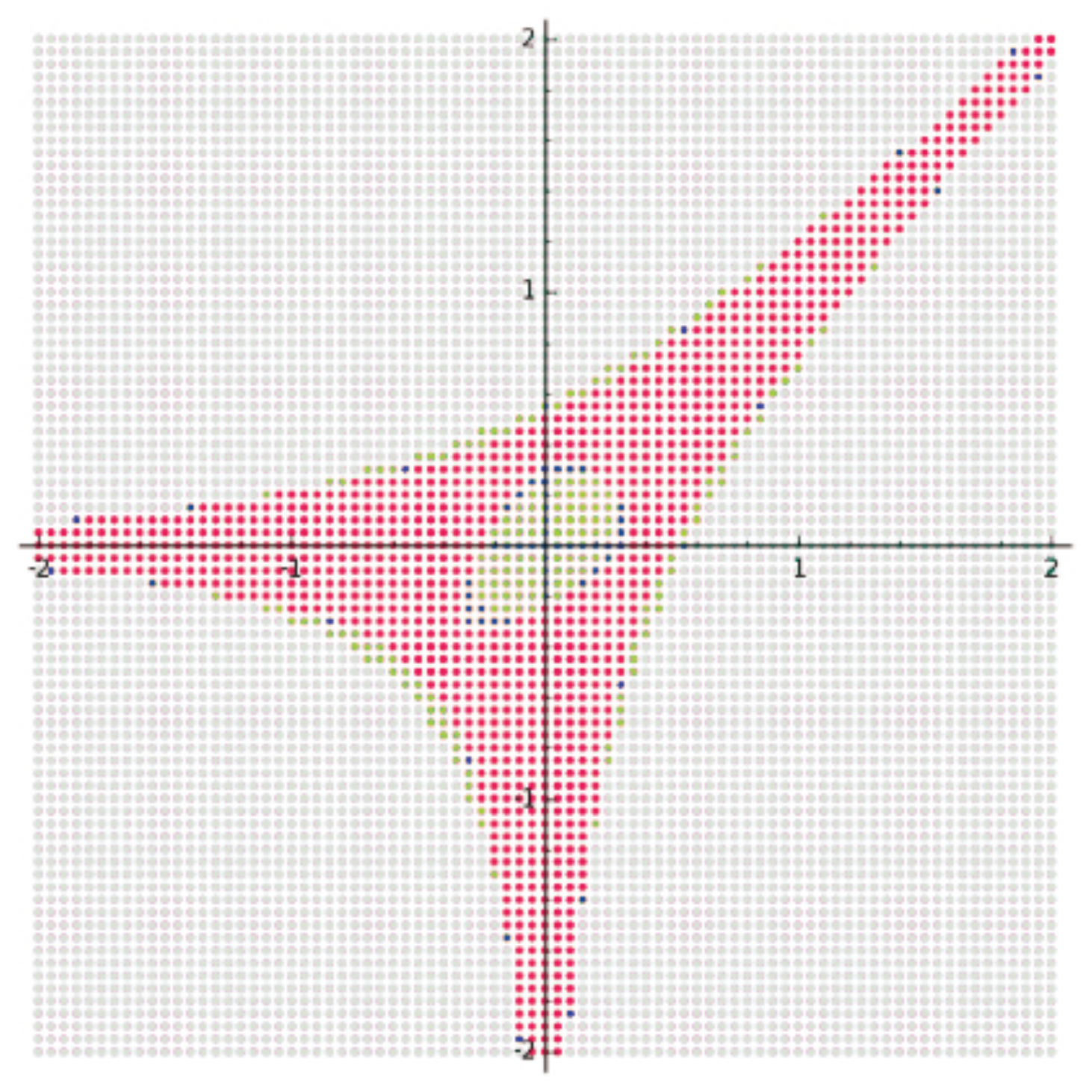}
	\caption{It is clear from this example that the lopsided amoeba contains the amoeba as a subset. Figure taken from \cite[Figure 1]{forsgard2017lopsided}.}\label{lopsidedamoeba}
\end{figure}

\paragraph{Example} For instance, consider the Newton polynomial $P=z^3+w^3+2zw+1$ whose amoeba is plotted red in Figure \ref{lopsidedamoeba}. The dark blue points (plus the red ones) form the lopsided amoeba $\mathcal{LA}_{\widetilde{P}_{16}}$ while the chartreuse points (plus the red and dark blue ones) give the region of $\mathcal{LA}_{\widetilde{P}_{8}}$.

\subsection{Amoebae from D-branes}\label{s:branes}
As mentioned in the introduction, amoebae are also interesting objects from the viewpoints of physics.
Historically, amoebae have a strong connection to dimer models where $P=0$ acts as the spectral curve of the bipartite graph \cite{zbMATH05051319,Kenyon:2003ui}. In particular, for a dimer with positive weights on its edges, $P=0$ is a \emph{Harnack curve}. 
D-branes probing toric Calabi-Yau (CY) manifolds has been very well studied in the past few decades especially in the context of quiver \footnote{q.v.~\cite{Bao:2020nbi} for a recent ML treatment of quiver gauge theories.} gauge theories \cite{Feng:2000mi,He:2001ey,Feng:2004uq}. 
Brane tilings are precisely dimer models that successfully show how one can connect, for example, the world-volume theories of a stack of D3-branes and the toric geometry of CY threefolds \cite{Feng:2001bn,Hanany:2005ve,Franco:2005rj}. It is also well-known that under mirror symmetry, this gives rise to a description in terms of D6-branes wrapping the $\mathbb{T}^3$ fibres \cite{Hori:2000kt,Hori:2000ck}.

Briefly, as we are considering 2d lattice polygons, we will denote the Newton polytope as $P(z,w)$ with $z,w\in\mathbb{C}^*$, note this is the convention where the surface is in $(\mathbb{C}^*)^2$ such that subsequent use of $(z,w)$ over $(z_1,z_2)$ implies this. As shown in \cite{Hori:2000kt}, the mirror geometry of a CY$_3$ whose toric diagram is $\Delta(P)$ is the local threefold $P(z,w)=uv$ where $u,v$ are complex variables. As discussed in \cite{Feng:2005gw}, the curve $P(z,w)=0$ plays a crucial role in geometry and quiver theories. This is a punctured Riemann surface $\Sigma$ of genus $g$ where $g$ equals the number of interior points of the Newton polygon.

Furthermore, the ($p,q$)-webs are also related to the above CY$_3$ and the Riemann surface $\Sigma$. If we wrap the 5-branes on $\Sigma$ and compactify the theory on a torus, we can get the theory for D3-branes probing CY$_3$ after performing T-dualities on the two directions of the torus. In fact, the toric diagram reveals a simple connection: the dual graph of this toric polygon is exactly the ($p,q$)-web diagram. The number of the boundary points of the toric diagram is equal to the number of NS5 cycles in the brane system. 
As shown in \cite{Feng:2005gw}, the Riemann surface $\Sigma$ can also be thought of as the thickening of the web diagram; conversely, the web is a deformation retract of the amoeba.

\begin{figure}[t]
    $\begin{array}{c}
		\includegraphics[width=12cm]{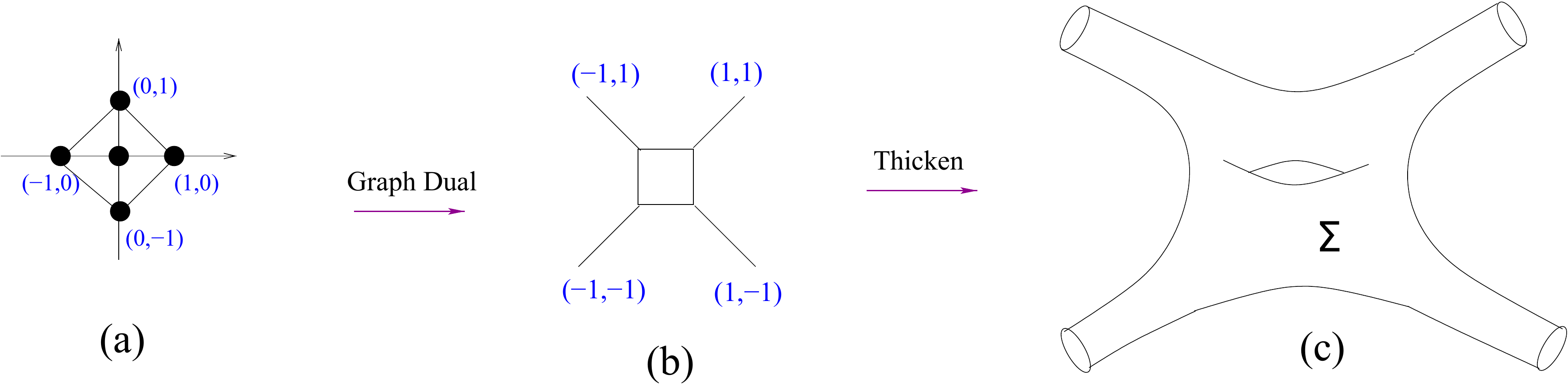}
	\end{array}
		\Rightarrow
		\begin{array}{c}
		\includegraphics[width=2.5cm]{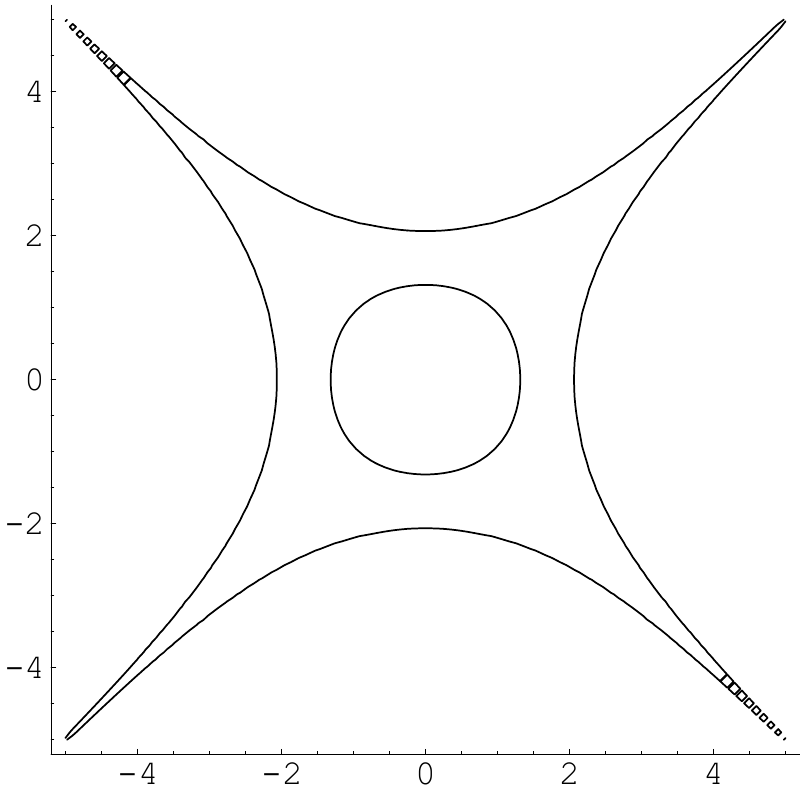}		\end{array}$
	\caption{Figure taken from \cite{Feng:2005gw} (figure 3): (a) The toric diagram of $F_0$. (b) The $(p,q)$-web is the dual graph of the toric diagram. (c) The holomorphic Riemann surface $\Sigma$ is a thickening of the web diagram.
	To the right we also include the amoeba for reference, this is a deformation retract of $\Sigma$.}\label{F0pq}
\end{figure}

As an example, the lattice polygon in Figure \ref{egAmoeba} (a) is actually the toric diagram of the affine Calabi-Yau cone over the zeroth Hirzebruch surface $F_0 \cong \mathbb{P}^1 \times \mathbb{P}^1$.
We summarize the above concepts in 
Figure \ref{F0pq}, using this example.
How one can obtain the brane tiling on $\mathbb{T}^2$ and the intersection locus of D6s on $\Sigma$ was studied in \cite{Feng:2005gw}. Following Definition \ref{ambdef}, we can take the projection from $\left(\mathbb{C}^*\right)^2$ to $\mathbb{R}^2$ to obtain the corresponding amoeba and plot it on the Log plane. 

\comment{
The example of $F_0$ is shown in Figure \ref{egF0amoebabound}.
\begin{figure}[h]
	\centering
		\centering
		\includegraphics[width=5cm]{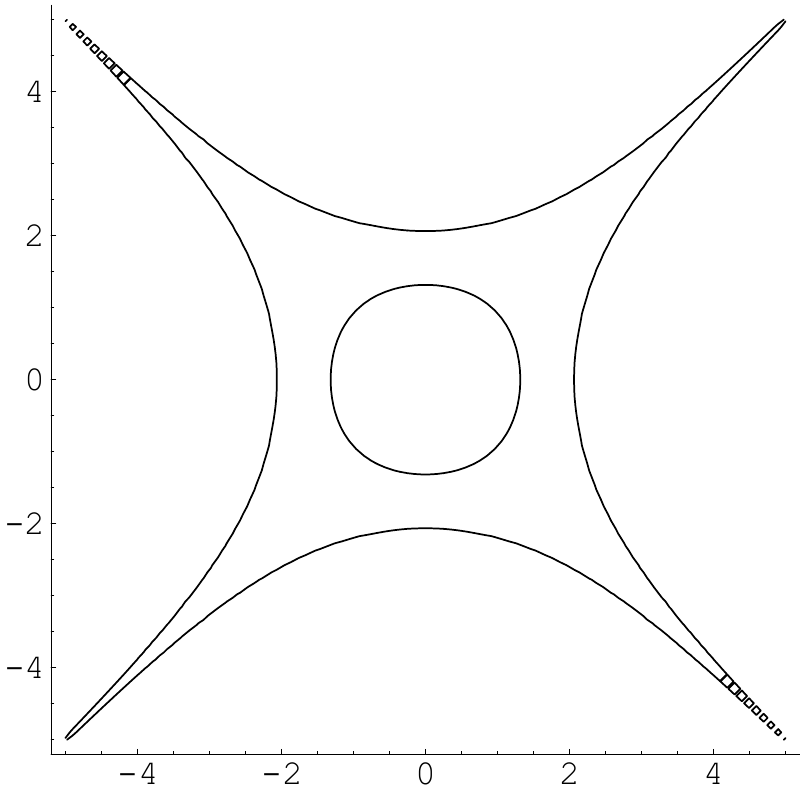}
	\caption{The amoeba for $F_0$. In this example, the Newton polynomial is $P(z,w)=z+w+1/z+1/w+6$. Figure taken from \cite{Feng:2005gw} (figure 15(a)).}\label{egF0amoebabound}
\end{figure}
}

We can see our familiar concepts such as genus, tentacles, spine, and thickening from the figure.
If this thickening is thick enough, the hole in the above example might disappear. 
For any general amoeba, the number of its genus is controlled by the coefficients of the Newton polynomial. For instance, the amoeba in Figures \ref{egAmoeba} and \ref{F0pq} would become genus 0 if $c\leq4$ for $P(z,w)=z+w+1/z+1/w+c$. 

Therefore, the coefficients in front of the monomials in the Newton polynomial are of particular interest in the study of amoebae as they determine the amoebae boundaries. In general, the coefficients can be any complex numbers and they are mirror to the K\"ahler moduli of the toric geometry. In dimer models, we are mainly interested in positive integer coefficients since they count perfect matchings \cite{Franco:2005rj}. For more details, one is referred to \cite{Feng:2005gw} as well as the review \cite{Yamazaki:2008bt}. In this paper, we will consider real coefficients as our moduli.

\section{ML Amoebae from Coefficients}\label{ml}
Having introduced all the basics and hopefully having motivated the readers, we shall study amoebae from the novel perspectives of machine learning.
We will consider the problems of computing genus and finding membership using feed-forward and convolutional neural networks, as well as using manifold learning. Unless specified,
the architectures are as follows:
\begin{description}
    \item[MLP] The multilayer perceptrons (MLPs) or more commonly known today as feed-forward NN, all have the structure with one hidden layer with 100 perceptrons/single-neurons. We will always use ReLU as the activation function\footnote{The rectified linear unit (ReLU) function is defined as ReLU$(q):=\text{max}(q,0)$. Similarly, LeakyReLU$(q)$ is defined to be $\alpha q$ when $q<0$ for some constant $\alpha$ and $q$ when $q\geq0$. 
    In our CNN models, we will take $\alpha=0.1$. Inspired by \cite{relu} single-layer ReLU NNs were selected as universal function approximators.} for MLPs.
    \item[CNN] The CNNs consist of four 1d convolutional layers, each followed by a LeakyReLU layer then a 1d MaxPooling layer.
\end{description}
 For all our NNs, the {\it learning rates} are all 0.001 and we always use the Adam optimizer \cite{kingma2017adam}.
 All computation time is matter of a few minutes even though the structure is not minimized. 
 Notice that it is possible to further reduce the structures for these MLPs and CNNs. In general, a not-carefully-tuned network can still give similar good performance. Therefore, the structures mentioned here are not so important. In the following, we will reduce the network to find the minimal structure whenever necessary.
 This independence of detailed architecture is again consistent with \cite{He:2021oav} which observed that, generically across disciplines, the robustness and absence of noise in pure mathematical data tend to give similar accuracies for rather different methods of NNs and classifiers.

\subsection{Lopsided Amoebae: $n=1$}\label{n1}
 As discussed above, we will use lopsided amoebae since they are more amenable to computation.
 We begin with the simplest case of $n=1$ where $\Tilde{P}_1 = P$ by definition, so that
 $\mathcal{LA}_{\Tilde{P}_1}=\mathcal{LA}_P$. 
 Can one predict the genus without explicit computation?
 
 For a fixed Newton polynomial $P(z,w)=\sum c_kz^iw^j$, the input is the vector composed of the coefficients, i.e., $\{c_1,c_2,\dots,c_n\}$, and the output is the genus.
 In other words, we have labelled data of the form
\begin{equation}
    \{c_1,c_2,\dots,c_n\} \longrightarrow g \ .
\end{equation}
Let us use some standard classifiers on this problem, in the spirit of \cite{He:2021oav}.

\subsubsection{Example 1: $F_0$}\label{F0}
We start with our simple running example $F_0$ whose toric diagram
was given in Figure \ref{egAmoeba} (a).
\comment{
\begin{equation}
	\includegraphics[width=1.5cm]{F0.pdf}.
\end{equation}
}
The Newton polynomial is
\begin{equation}
	P(z,w)=c_1z+c_2w+c_3z^{-1}+c_4w^{-1}+c_5.
\end{equation}
Hence, our input is $\{c_1,c_2,c_3,c_4,c_5\}$. Since the resulting lopsided amoeba could have at most one genus (corresponding to its single interior point), this is a binary classification where the output is either $g=0$ or $g=1$.
For this simple example, one can analytically derive (see Appendix \ref{lopsided}) the genus $g$ as a function of the coefficients:
\begin{equation}
	g=\begin{cases}
		0,&|c_5|\leq2|c_1c_3|^{1/2}+2|c_2c_4|^{1/2}\\
		1,&|c_5|>2|c_1c_3|^{1/2}+2|c_2c_4|^{1/2}
	\end{cases}.\label{gF0}
\end{equation}
Nevertheless, we can see from the RHS that the boundaries of decision are already non-trivial even for this simplest of examples.

\paragraph{Real Coefficients in Newton Polynomial: } We generate a balanced dataset with $\sim2000$ random samples, with $c_{1,2,3,4}\in[-5,5]$ and $c_5\in[-20,20]$.
Our forward-feed NN (MLP) can easily reach $0.957(\pm0.005)$ accuracy (with $95\%$ confidence interval) for a 5-fold cross validation\footnote{We have also tried some other classification models such as support vector machine, random forests etc. It turns out that these models would have nice performance with accuracies around $0.900(\pm0.010)$ for 5-fold cross validation although the MLP behaves slightly better. Henceforth, we shall always stick to MLP.}.

\begin{figure}[t]
	\centering
	\begin{subfigure}{10cm}
		\centering
		\includegraphics[width=10cm]{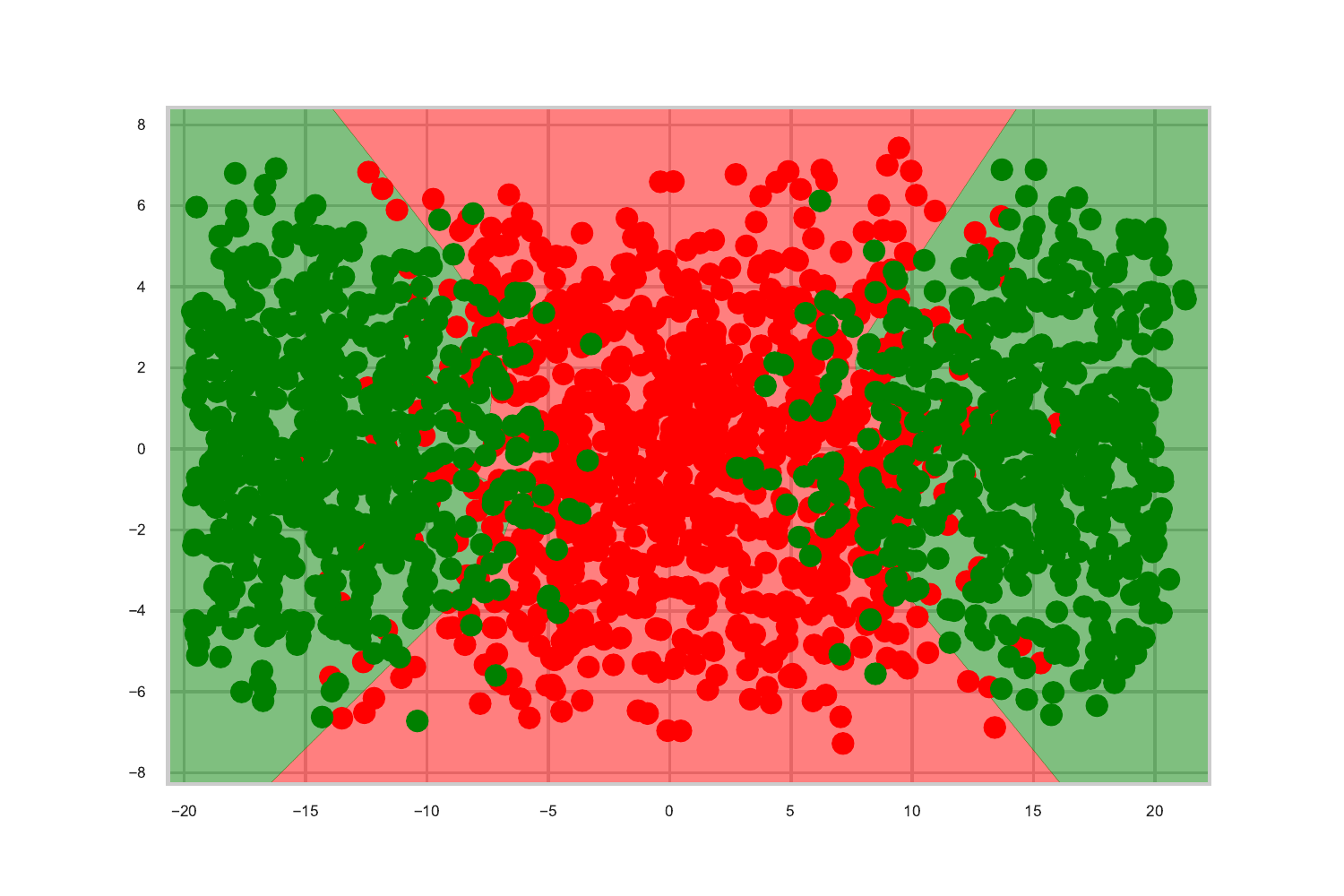}
		\caption{}
	\end{subfigure}
	\begin{subfigure}{10cm}
		\centering
		\includegraphics[width=10cm]{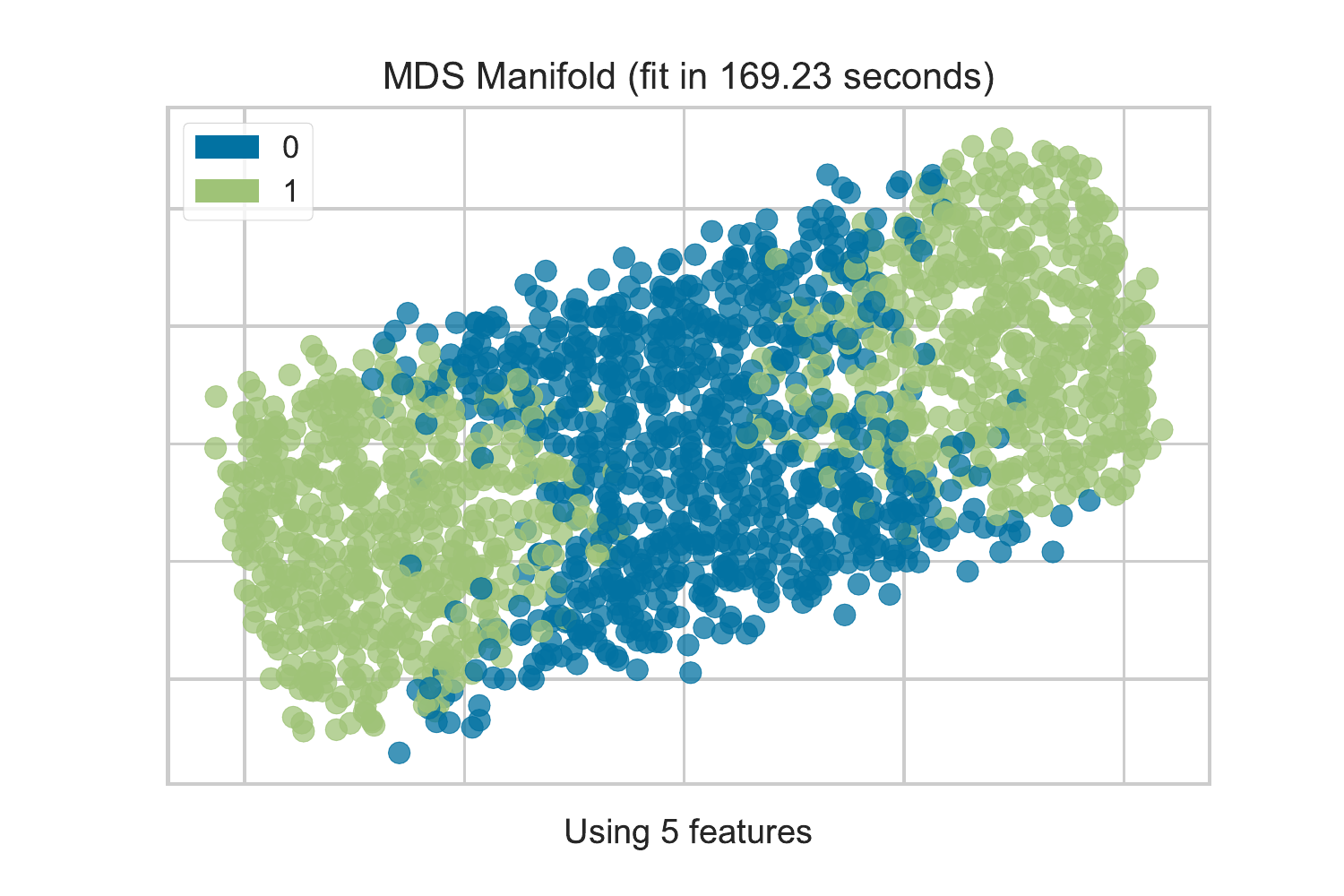}
		\caption{}
	\end{subfigure}
    \begin{subfigure}{5cm}
    	\centering
    	\includegraphics[width=5cm]{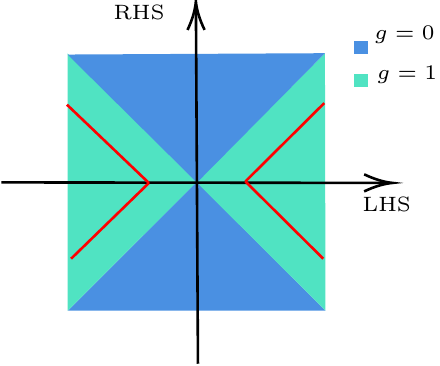}
    	\caption{}
    \end{subfigure}
	\caption{(a) The kernel PCA projection for an NN. (b) The MDS manifold projection which gives a better separation of the two classes of points. (c) Ideally, the blue and green regions would be separated by $y=\pm x$. In practice, due to the complication of square roots, the NN would get shifted. This shift, i.e., the actual separation of the blue and green points, is indicated by the red lines here.}\label{F0PCA}
\end{figure}

To see how well the NN is learning the above equation for $g$, we can a perform principal component analysis (PCA) projections, as visualized in Figure \ref{F0PCA}(a). With the help of the \texttt{Yellowbrick} package \cite{bengfort_yellowbrick_2019}, we can use multi-dimensional scaling (MDS) manifold projection to make the two types of data points further separated as in Figure \ref{F0PCA}(b).
For a brief introduction to different methods of manifold learning, see Appendix \ref{manifoldlearning}.

To get an idea about how these data points are distributed in the projection, let us define\footnote{As discussed in Appendix \ref{lopsided}, \eqref{gF0} works for all complex coefficients. However, we are only using real input vectors here, and more importantly, only the absolute values would matter in the condition. Hence, we can set $c_5=x\in\mathbb{R}$ in the projection.} $x\equiv c_5$ and $|y|\equiv2|c_1c_3|^{1/2}+2|c_2c_4|^{1/2}$. Then the inequalities in \eqref{gF0} has the boundary lines $y=\pm x$. As depicted in Figure \ref{F0PCA}(c), the two lines divide the projection plane into different regions, where in the blue region we have genus zero while in the green region we have genus one. The equivalent separation coming from our NN in Figure \ref{F0PCA}(a,b) is represented by the red lines in (c). We can see that they do not intersect at the origin and the $g=0$ region would occupy some of the $g=1$ parts. This could account for the $0.95$ accuracy (rather than $100\%$).

Why is there a shift of the boundary lines? We believe that this is due to the non-trivial expression of $y$ (especially the square roots therein). Although it is possible for an NN to learn such expression in principle, it could still be too complicated for a simple NN to fully simulate this.

\begin{figure}[h]
	\centering
	\begin{subfigure}{7cm}
		\centering
		\includegraphics[width=7cm]{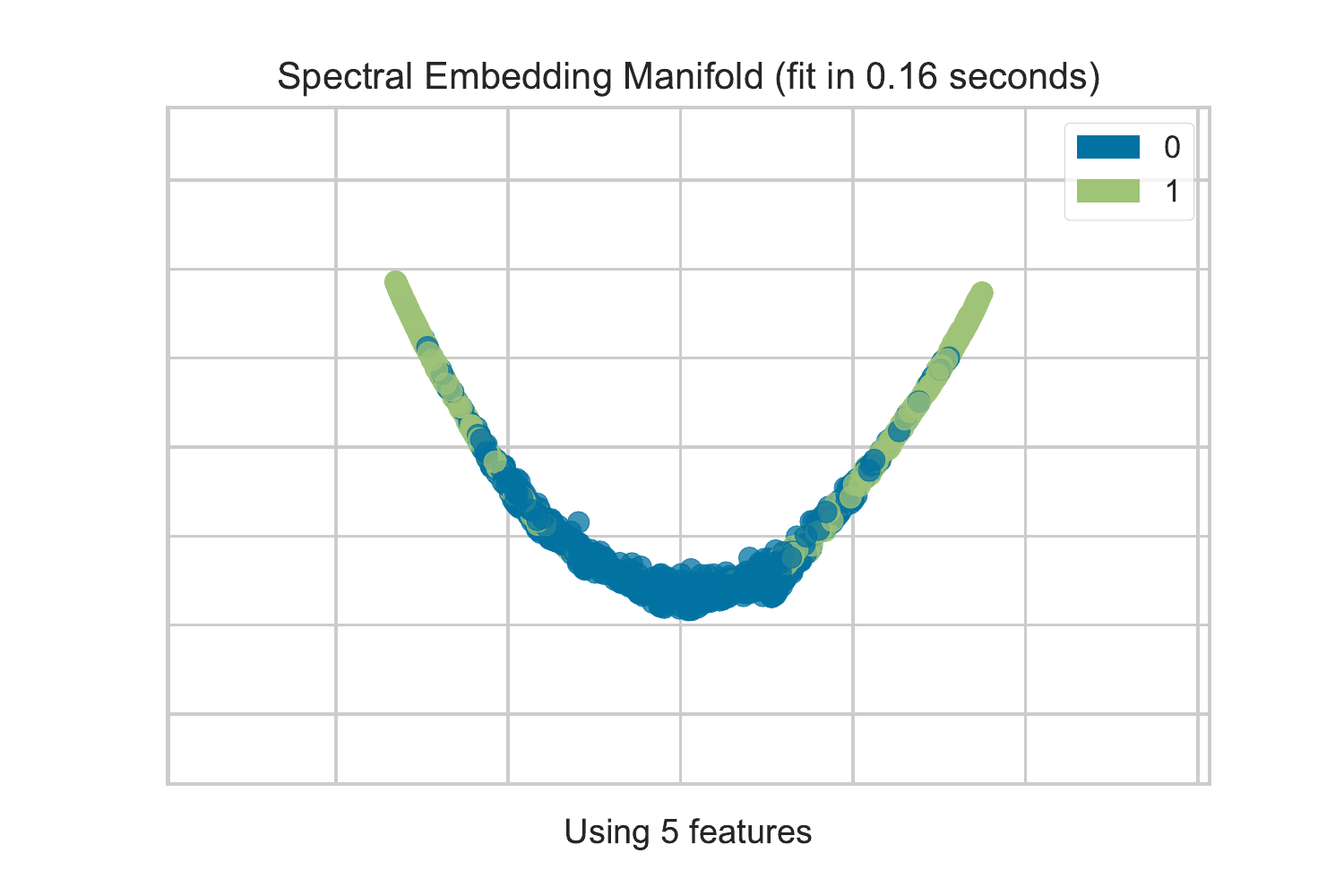}
		\caption{}
	\end{subfigure}
	\begin{subfigure}{7cm}
		\centering
		\includegraphics[width=7cm]{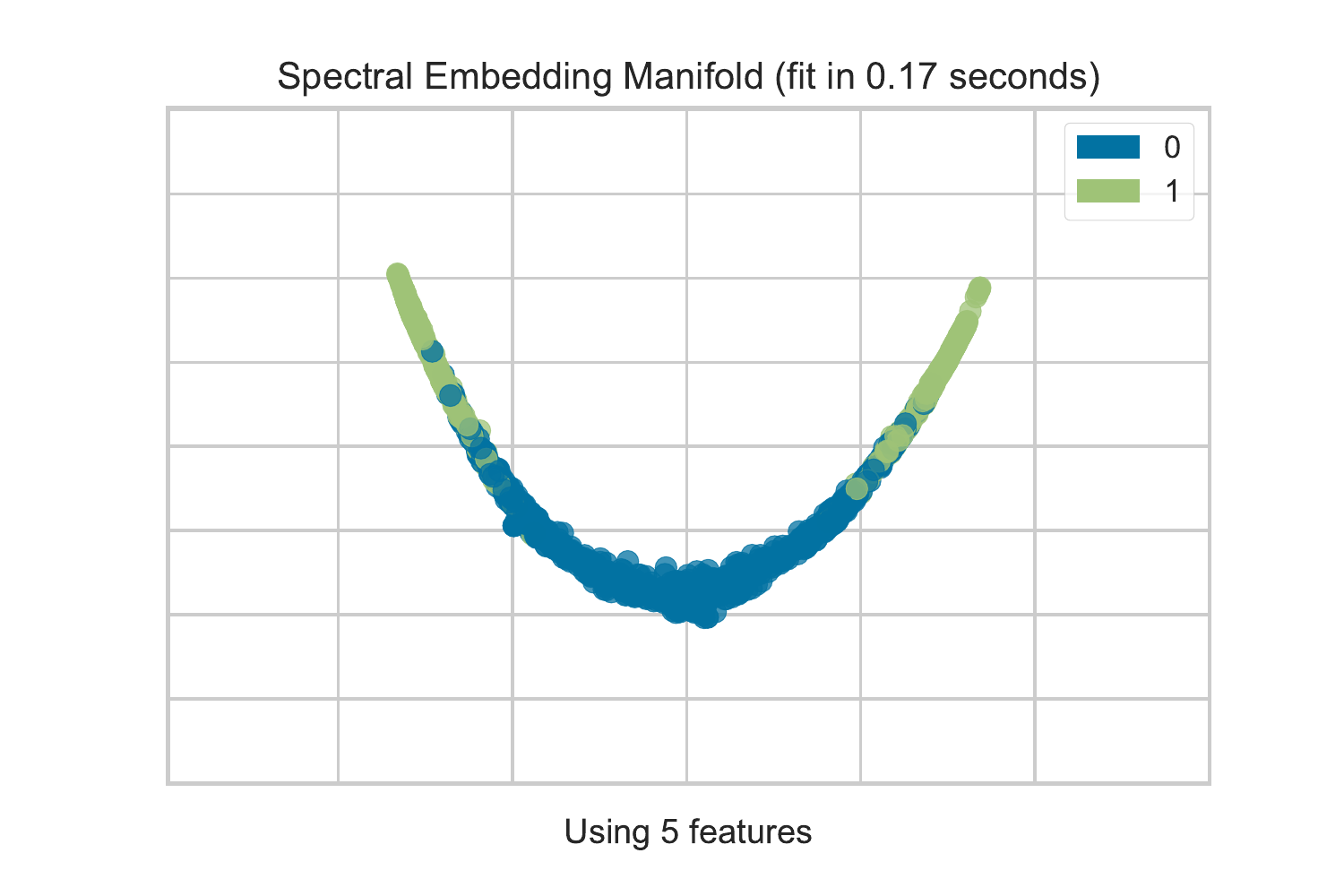}
		\caption{}
	\end{subfigure}
	\caption{(a) The spectral embedding manifold projection of the dataset for $F_0$. For the input vectors, $c_{1,2,3,4}$ range from $-5$ to 5 while $c_5$ ranges from $-20$ to 20. (b) To verify our explanation for (a), we generate a dataset whose two classes are separated by $c_5'^2=4|c_1'c_3'|+4|c_2'c_4'|+c_0$ (with $c_i'$ the same range as $c_i$), where $c_0=4\times2\times2.5^2$ as $2.5$ is the average of the possible values for $c_{1,2,3,4}$.}\label{parabola}
\end{figure}

There is another useful projection in this example, that is, the spectral embedding visualization \cite{bengfort_yellowbrick_2019} as shown in Figure \ref{parabola}(a).
As we can see, this gives a distribution in a parabola shape. We argue that this comes from squaring $y=\pm x$, i.e., $Y(\equiv y^2)=x^2$. However, again due to the complicated expression, it would be very hard for a simple NN to fully recover $Y$. Hence, there could still be some mixing at the boundary parts. As a validation, we generate a dataset whose binary classes are separated by the bound $c_5'^2=4|c_1'c_3'|+4|c_2'c_4'|+50$ and plot its spectral embedding as in Figure \ref{parabola}(b). Indeed, we obtain a similar parabola-shaped projection. Moreover, since the square root part is replaced by the constant 50, NN could also give higher accuracy ($\sim0.97$ at same training percentage).

From \eqref{gF0}, we also learn that only the absolute values of the coefficients would matter. In fact, this could be reflected in machine learning as well. We can use the same dataset, but with $|c_i|$'s as input. Indeed, this helps the model to improve its performance and a 5-fold cross validation now gives $0.987(\pm0.007)$ accuracy.

\begin{figure}[h]
	\centering
	\begin{subfigure}{7cm}
		\centering
		\includegraphics[width=7cm]{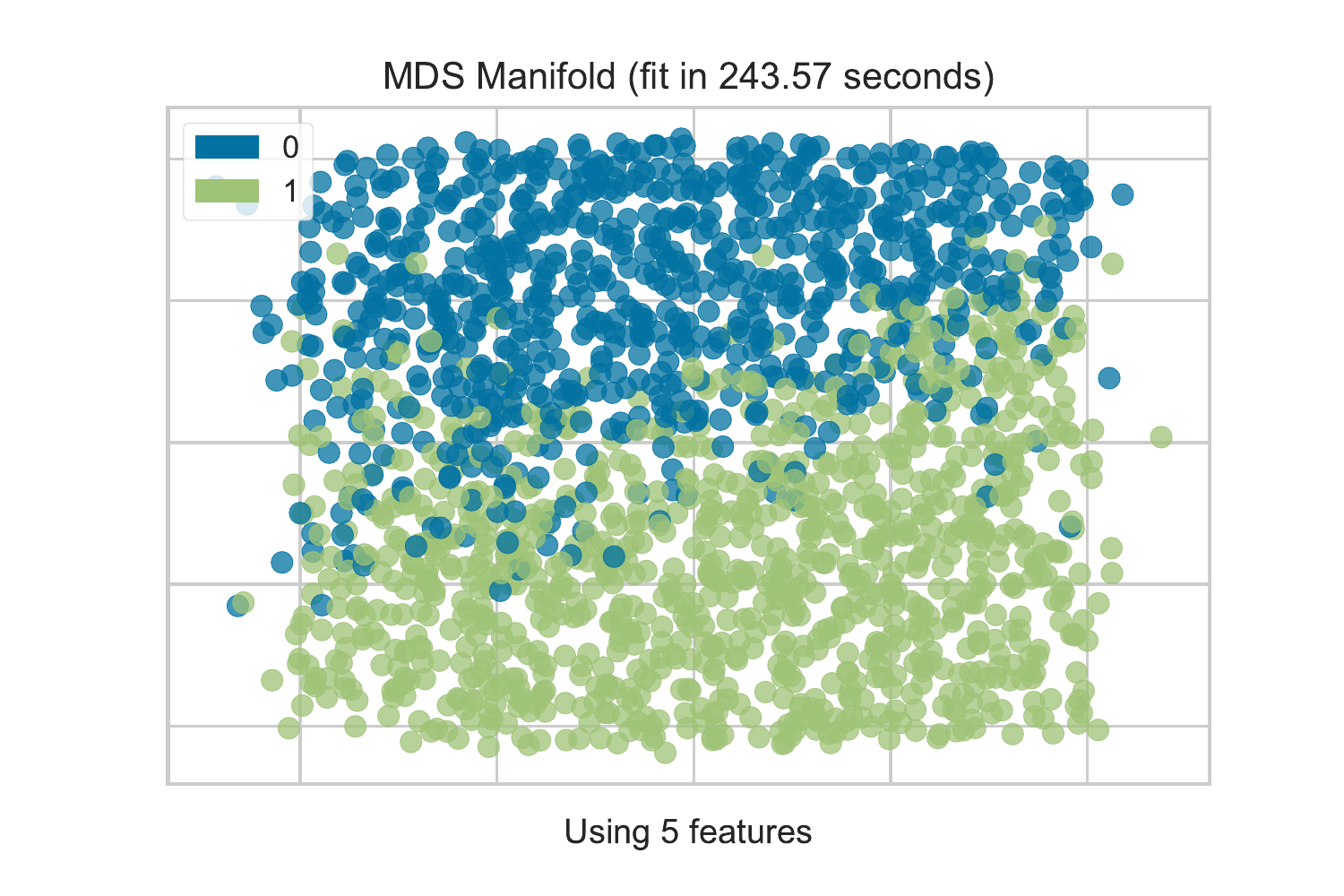}
		\caption{}
	\end{subfigure}
	\begin{subfigure}{7cm}
		\centering
		\includegraphics[width=7cm]{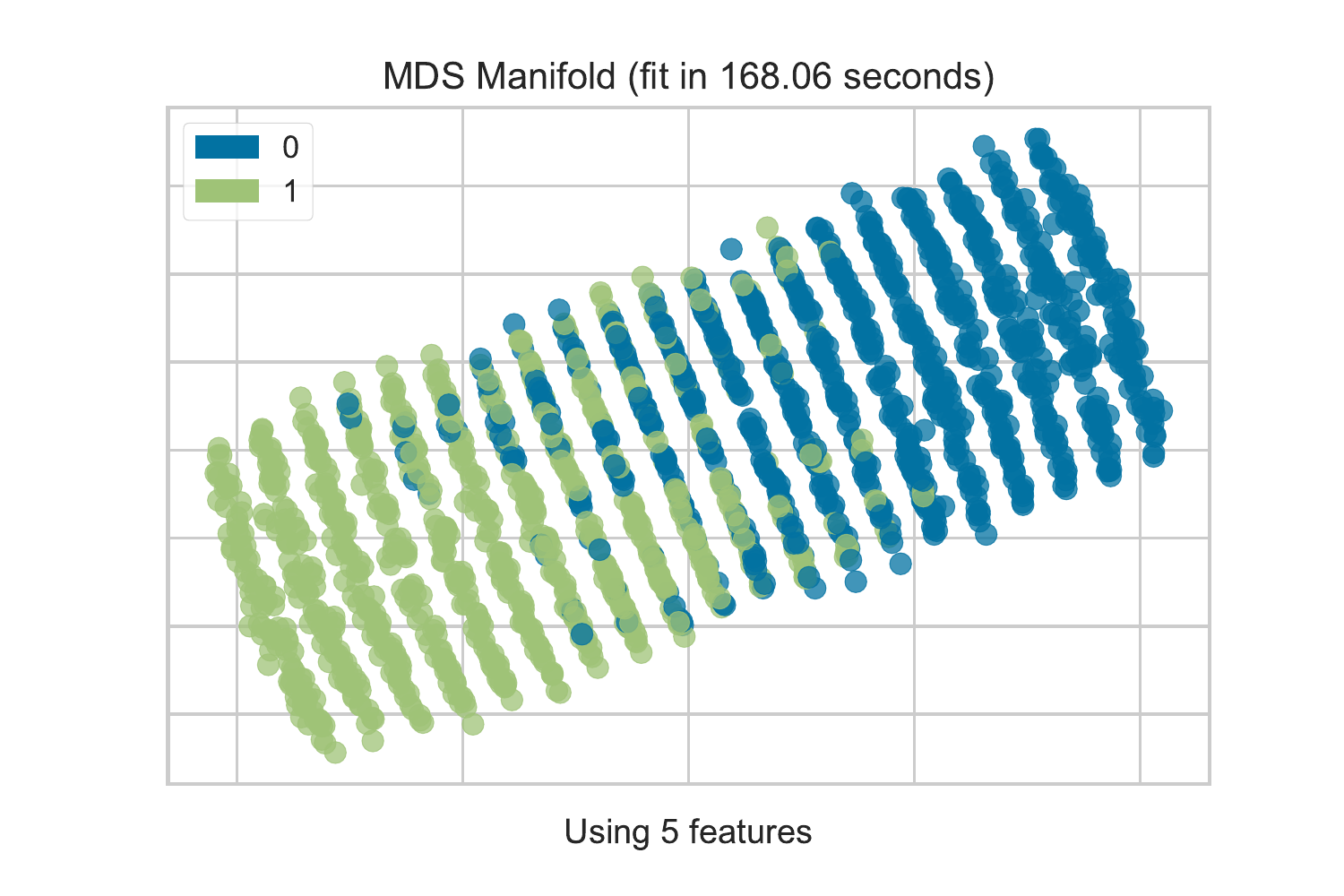}
		\caption{}
	\end{subfigure}
	\caption{(a) The MDS manifold projection for positive real input. (b) The MDS manifold projection for positive integer input. Since we are taking positive values, the plots would correspond to the first quadrant in Figure \ref{F0PCA}.}\label{f0PCApositive}
\end{figure}

\paragraph{Integer Coefficients in Newton Polynomial: }
As discussed in \S\ref{s:branes}, the coefficients $c_i$ have significance: they count perfect matching in dimer models and so-called GLSM quantum fields in brane tilings in string theory.
Thus, positive integer coefficients are of particular interest. We therefore generate $\sim2000$ data points with only positive integer inputs ($c_{1,2,3,4}\in[1,5],~c_5\in[1,20]$) for machine learning as well. It turns out that the accuracy is further improved to $0.992(\pm0.004)$ for 5-fold cross validation. We summarize the results, contrasting the input type, in Table \ref{f0table}.
\begin{table}[h]
\centering
\begin{tabular}{c||c|c|c}
 Input type & $\mathbb{R}^5$ & $(\mathbb{R}^+)^5$ & $(\mathbb{Z}^+)^5$ \\ \hline
 Accuracy & $0.957(\pm0.005)$ & $0.987(\pm0.007)$ & $0.992(\pm0.004)$
\end{tabular}
\caption{The performance for MLP using 5-fold cross validation with $95\%$ confidence interval.
The input is a five vector of the coefficients $c_i$ in the Newton polynomial for $F_0$, and we contrast the 3 different types of choices for $c_i$.}\label{f0table}
\end{table}
We also plot the MDS projections for the positive inputs in Figure \ref{f0PCApositive}.

\subsubsection{Example 2: $L^{3,3,2}$}\label{L332}
Let us consider a non-reflexive example with two interior points \footnote{The 16 reflexive polygons were considered in the context of brane-tiling in \cite{Hanany:2012vc,He:2017gam}, and the 45 non-reflexives with 2 interior points, in \cite{Bao:2020kji}.}, viz., $L^{3,3,2}$ whose toric diagram is
\begin{equation}
	\includegraphics[width=2cm]{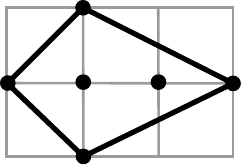}.
\end{equation}
The Newton Polynomial is
\begin{equation}
	P(z,w)=c_1z+c_2w+c_3z^{-1}+c_4w^{-1}+c_5z^2+c_6.
\end{equation}
Hence, our input is $\{c_1,c_2,c_3,c_4,c_5,c_6\}$. Since the resulting lopsided amoeba could have at most two holes (corresponding to its two interior points), this is a ternary classification where the output can be $g=0,1,2$.

As with $F_0$, we can analytically find the boundary decisions, though the expressions are much more complicated.
The details are presented in Appendix \ref{lopsided}, and we summarize them here:
\begin{equation}
	g=\begin{cases}
		0,&|c_1|\leq a_1~\text{and}~|c_6|\leq a_2\\
		1,&(|c_1|>a_1~\text{and}~|c_6|\leq a_2)~\text{or}~(|c_1|\leq a_1~\text{and}~|c_6|>a_2)\\
		2,&|c_1|>a_1~\text{and}~|c_6|>a_2
	\end{cases},
\end{equation}
where $a_1 := |c_2|w_0/z_0+|c_3|/z_0^2+|c_4|/(z_0w_0)+|c_5|z_0+|c_6|/z_0$ and $a_2 := |c_1|z_0'+|c_2|w_0+|c_3|/z_0'+|c_4|/w_0+|c_5|z_0'^2$ such that
$z_0 := \sqrt[3]{-\frac{q}{2}+\sqrt{\Delta}}+\sqrt[3]{-\frac{q}{2}-\sqrt{\Delta}}$, with
$\Delta :=\left|\frac{c_3}{c_5}\right|^2-\frac{\left(2|c_2c_3|^{1/2}+|c_6|\right)^3}{27|c_5|^3}$, and
$z_0'$ is the positive root of $2|c_5|z_0'^3+|c_1|z_0'^2-|c_3|=0$.

We generate a balanced set with only $\sim9000$ random samples. A 5-fold cross validation for MLP gives $0.912(\pm0.006)$ accuracy (with $95\%$ confidence interval).

\begin{figure}[h]
	\centering
	\begin{subfigure}{7cm}
		\centering
		\includegraphics[width=7cm]{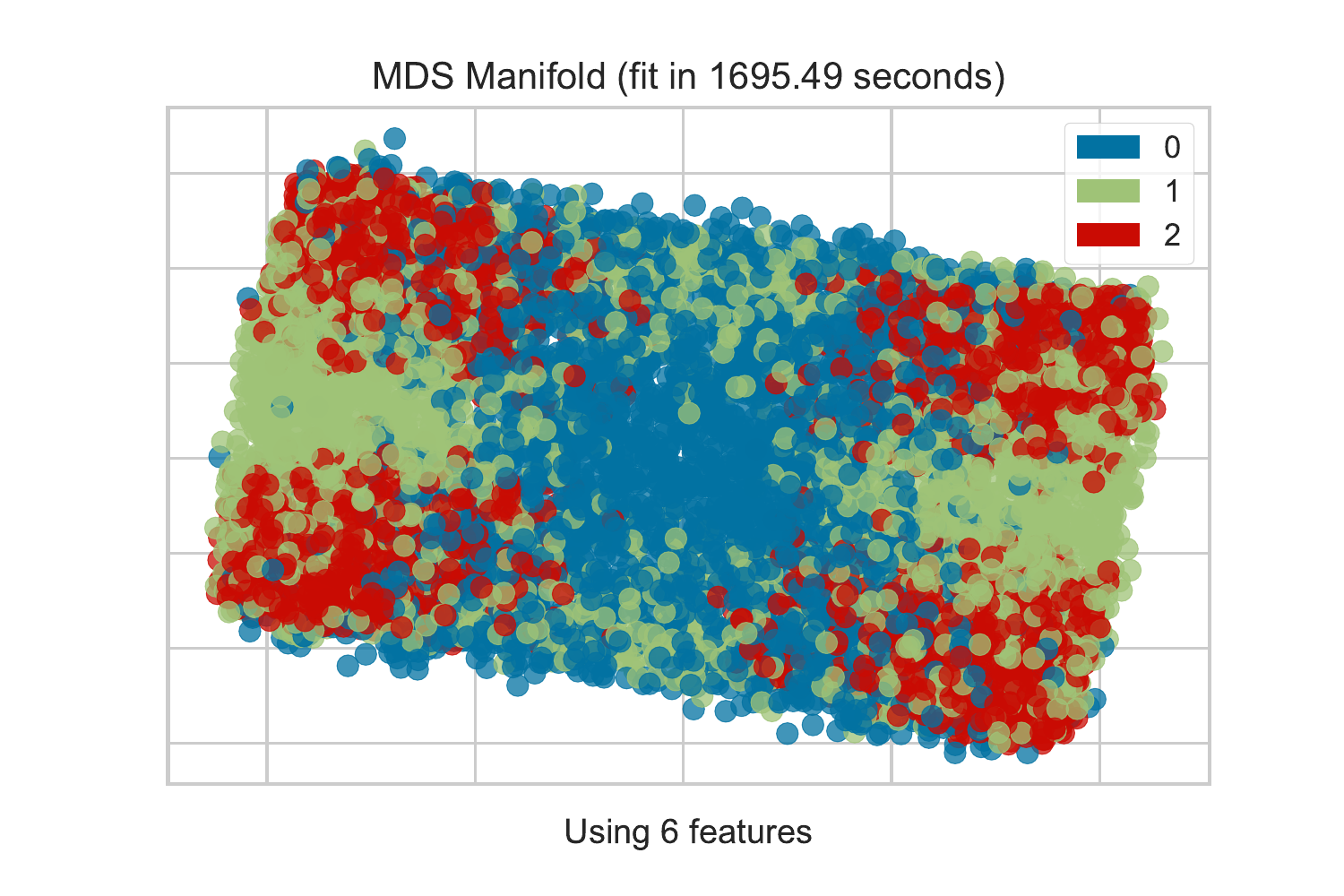}
		\caption{}
	\end{subfigure}
	\begin{subfigure}{7cm}
		\centering
		\includegraphics[width=7cm]{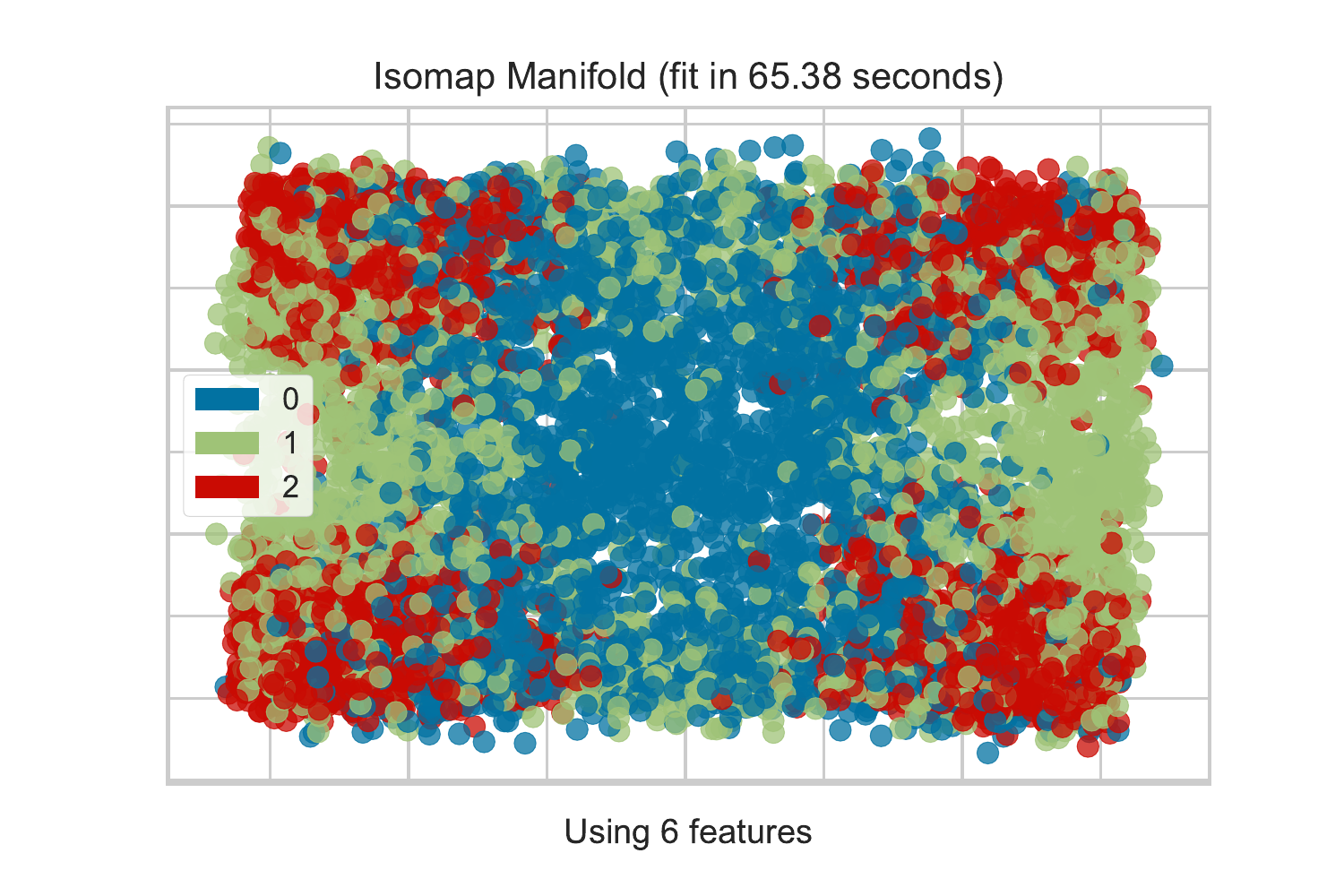}
		\caption{}
	\end{subfigure}
	\begin{subfigure}{8cm}
		\centering
		\includegraphics[width=7cm]{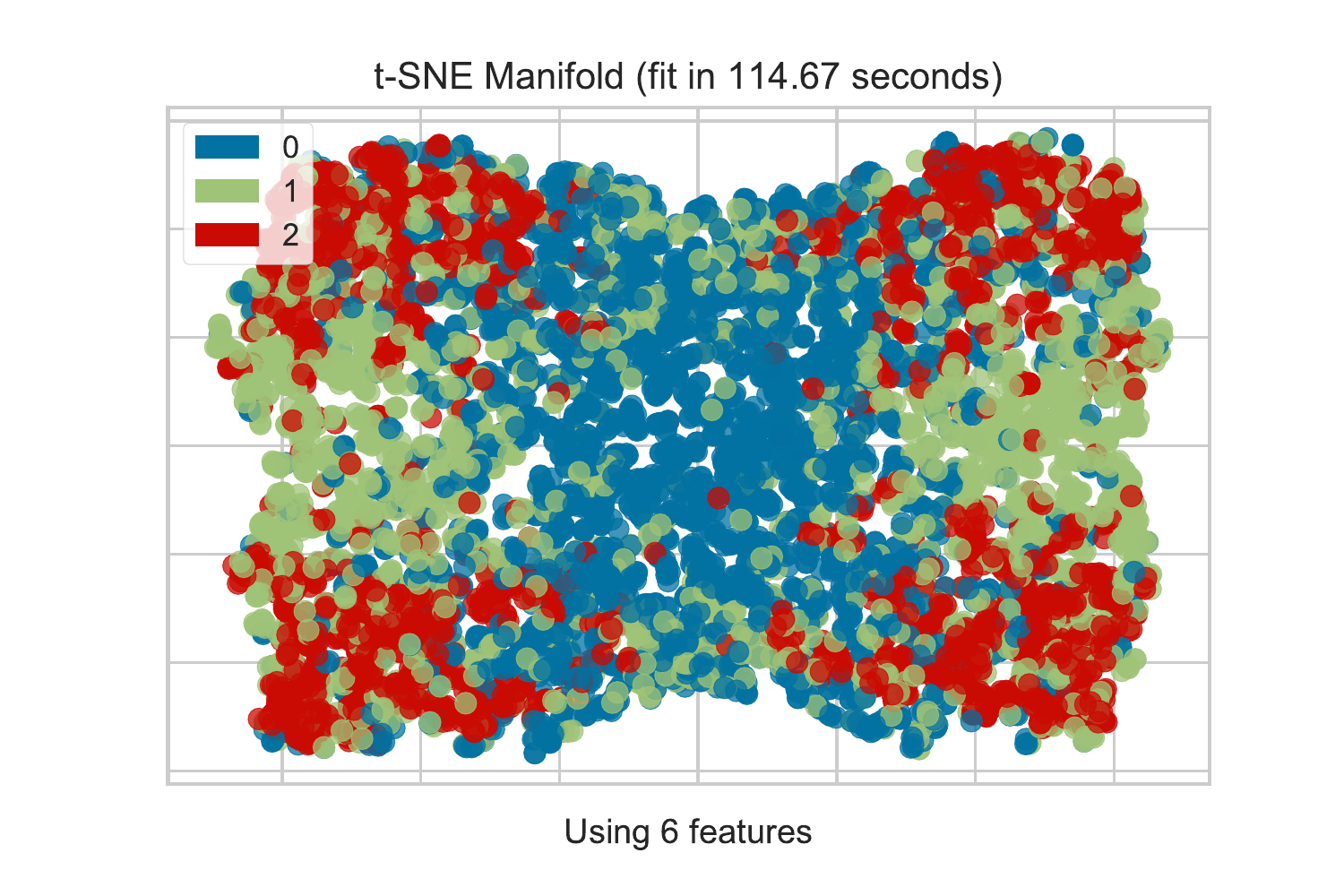}
		\caption{}
	\end{subfigure}
\begin{subfigure}{5cm}
	\centering
	\includegraphics[width=5cm]{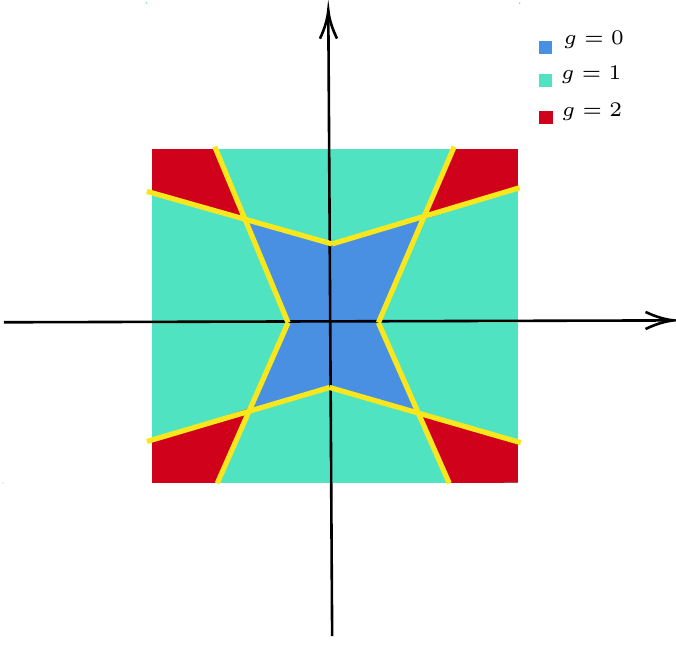}
	\caption{}
\end{subfigure}
	\caption{(a) The MDS embedding manifold projection of the dataset for $L^{3,3,2}$. (b) The Isomap embedding manifold projection. (c) The t-SNE embedding manifold projection. For instructions on these manifold projections, one is referred to \cite{bengfort_yellowbrick_2019}. (d) The sketch of an ideal separation of the data points.}\label{L332PCA}
\end{figure}

Again, we can plot different manifold projections to visualize how the NN is simulating the above condition as shown in Figure \ref{L332PCA}(a,b,c).
We can see that different projections give similar data separations. To understand such decision regions, let us first write the two bounds as
\begin{equation}
	\begin{split}
		|c_1|&=a_1\equiv|c_6|/z_0+b_1,\\
		|c_6|&=a_2\equiv|c_1|z_0'+b_2,
	\end{split}
\end{equation}
where $b_i>0$. In other words, when $g=2$, the region is bounded by $|c_6|/z_0+b_1<|c_1|<|c_6|/z_0'-b_2$. This corresponds to the red region in Figure \ref{L332PCA}(d) whose horizontal and vertical axes are $x\equiv\pm|c_6|$ and $y\equiv\pm|c_1|$ respectively. Notice that $z_0,z_0',b_i$ are not constants, but we can always draw some boundary lines (in yellow) as a sketch (assuming that the changes of the boundary lines are small compared to different coloured regions\footnote{For instance, here $c_{1,6}$ are generated in the range $[-20,20]$, and the means are $\mu(1/z_0)=0.31,\mu(1/z_0')=2.01,\mu(b_1)=6.48,\mu(b_2)=9.18$ with standard deviations $\sigma(1/z_0)=0.18,\sigma(1/z_0')=0.77,\sigma(b_1)=3.72,\sigma(b_2)=3.07$. We will not restate this explicitly for the examples discussed below. The complicated expressions for $b_i$ and the standard deviations could account for the mis-classifications in machine learning.}). 

Likewise, it is straightforward to find the regions for $g=0$ (in blue) and $g=1$ (in green). Indeed, this is the decision regions we get from those different projections. As $z_0$, $z_0'$ and $b_i$ are not really constants and have some complicated expressions, it is natural to see some mixings in the projections.
With the same dataset, we can also take the absolute values of the coefficients as input since only $|c_i|$ is relevant for counting the genus. For 5-fold cross validation, the accuracy is significantly improved to $0.968(\pm0.004)$.

\begin{figure}[h]
	\centering
	\begin{subfigure}{7cm}
		\centering
		\includegraphics[width=7cm]{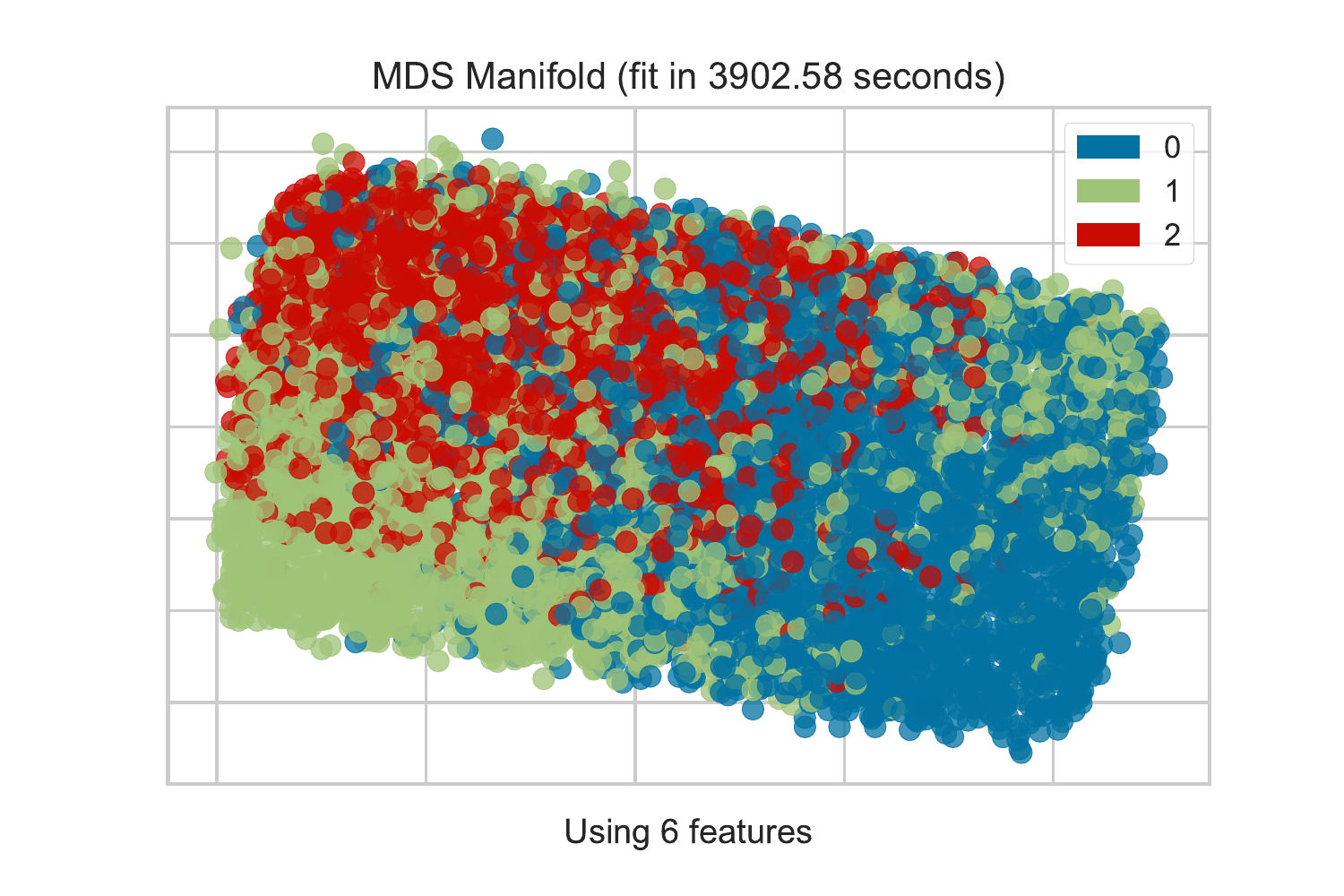}
		\caption{}
	\end{subfigure}
	\begin{subfigure}{7cm}
		\centering
		\includegraphics[width=7cm]{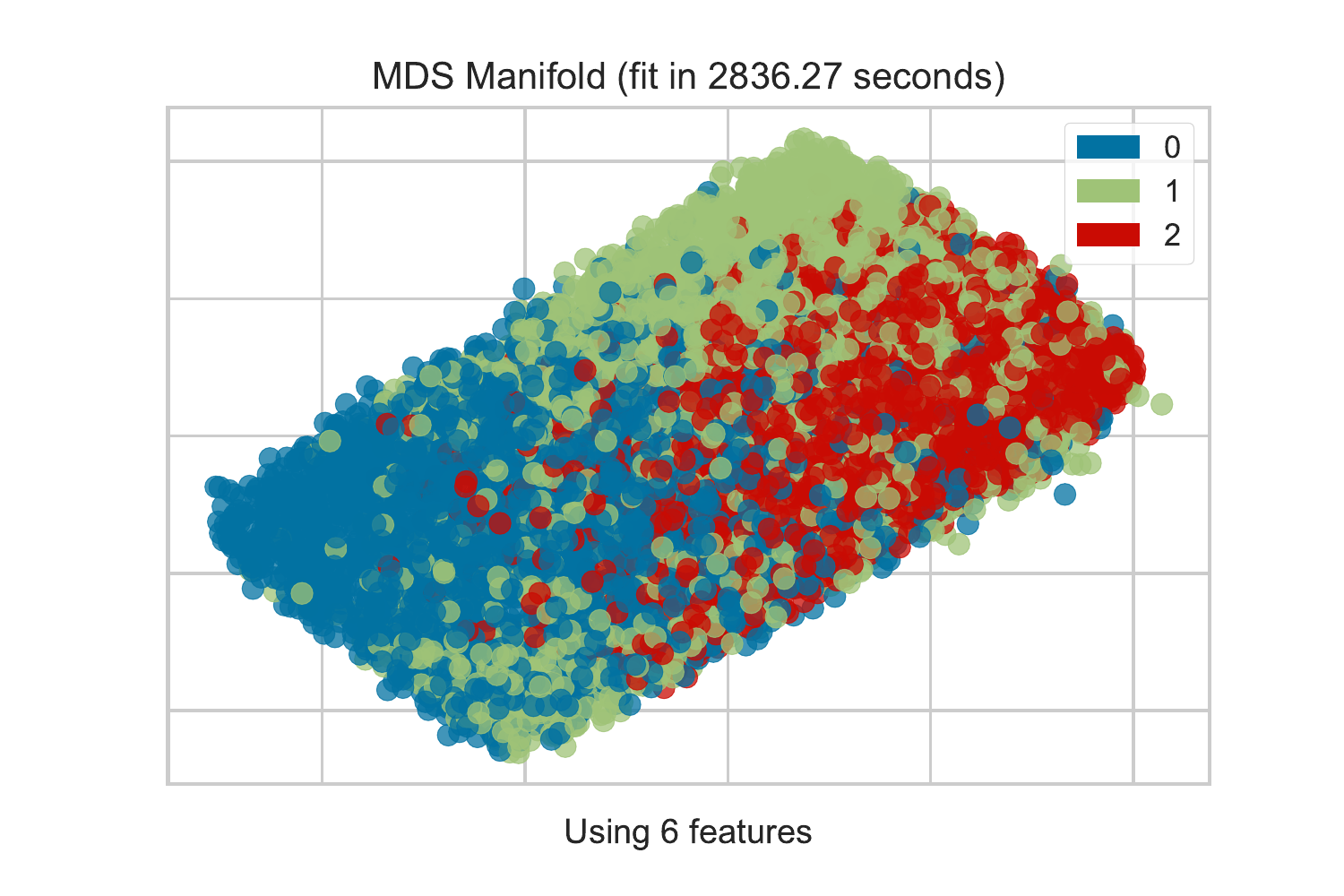}
		\caption{}
	\end{subfigure}
	\caption{(a) The MDS manifold projection for positive real input. (b) The MDS manifold projection for positive integer input. Since we are taking positive values, the plots would correspond to the first quadrant in Figure \ref{L332PCA}.}\label{l332PCApositive}
\end{figure}

As before, we may also generate a dataset of similar size with only positive integer coefficients. For 5-fold cross validation, the accuracy is improved to $0.990(\pm0.003)$.  We summarize the results in Table \ref{l332table}.
\begin{table}[h]
\centering
\begin{tabular}{c||c|c|c}
 Input type & $\mathbb{R}^6$ & $(\mathbb{R}^+)^6$ & $(\mathbb{Z}^+)^6$ \\ \hline
 Accuracy & $0.912(\pm0.006)$ & $0.968(\pm0.004)$ & $0.990(\pm0.003)$
\end{tabular}
\caption{The performance for MLP using 5-fold cross validation with $95\%$ confidence interval.}\label{l332table}
\end{table}
We also plot the MDS projections as an example for the positive inputs in Figure \ref{l332PCApositive}.

\subsubsection{Example 3: $\mathcal{C}/(\mathbb{Z}_2\times\mathbb{Z}_4)$}\label{CZ2Z4}
We now contemplate an example with three interior points, that is, $\mathcal{C}/(\mathbb{Z}_2\times\mathbb{Z}_4)$ with actions $(1,0,0,1); (0,1,3,0)$, an Abelian quotient of the conifold $\mathcal{C}$. Its toric diagram is
\begin{equation}
	\includegraphics[width=3cm]{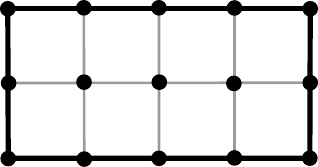}.
\end{equation}
and the Newton polynomial is
\begin{equation}
    \begin{split}
         P(z,w)&=c_0+c_1z+c_2z^2+c_3z^3+c_4z^4+c_5w+c_6zw+c_7z^2w+c_8z^3w+c_9z^4w\\
         &c_{10}w^2+c_{11}zw^2+c_{12}z^2w^2+c_{13}z^3w^2+c_{14}z^4w^2.
    \end{split}
\end{equation}
Hence, our input is $\{c_0,c_1,c_2,\dots,c_{14}\}$. Since the resulting lopsided amoeba could have at most three genera (corresponding to its three interior points), this is a classification where the output can be $g=0,1,2,3$.

We generate $\sim12000$ random samples with $c_i\in[-30,30]$, using the coefficient vectors and their absolute values respectively as inputs for two experiments. We also generate $\sim12000$ random coefficient vectors with only positive integer entries. It turns out that a CNN has a better performance than MLP. We list the results for 5-fold cross validation in Table \ref{cz2z4table}.
\begin{table}[h]
\centering
\begin{tabular}{cc||c|c|c}
\multicolumn{2}{c||}{Input type} & $\mathbb{R}^{15}$ & $(\mathbb{R}^+)^{15}$ & $(\mathbb{Z}^+)^{15}$ \\ \hline
\multicolumn{1}{c|}{\multirow{2}{*}{Accuracy}} & MLP & $0.792(\pm0.006)$ & $0.856(\pm0.005)$ & $0.890(\pm0.005)$ \\ \cline{2-5} 
\multicolumn{1}{c|}{}  & CNN & $0.803(\pm0.009)$ & $0.909(\pm0.007)$ & $0.907(\pm0.002)$
\end{tabular}
\caption{The performances for MLP and CNN using 5-fold cross validation with $95\%$ confidence interval.}\label{cz2z4table}
\end{table}

\begin{figure}[h]
	\centering
	\begin{subfigure}{7cm}
		\centering
		\includegraphics[width=7cm]{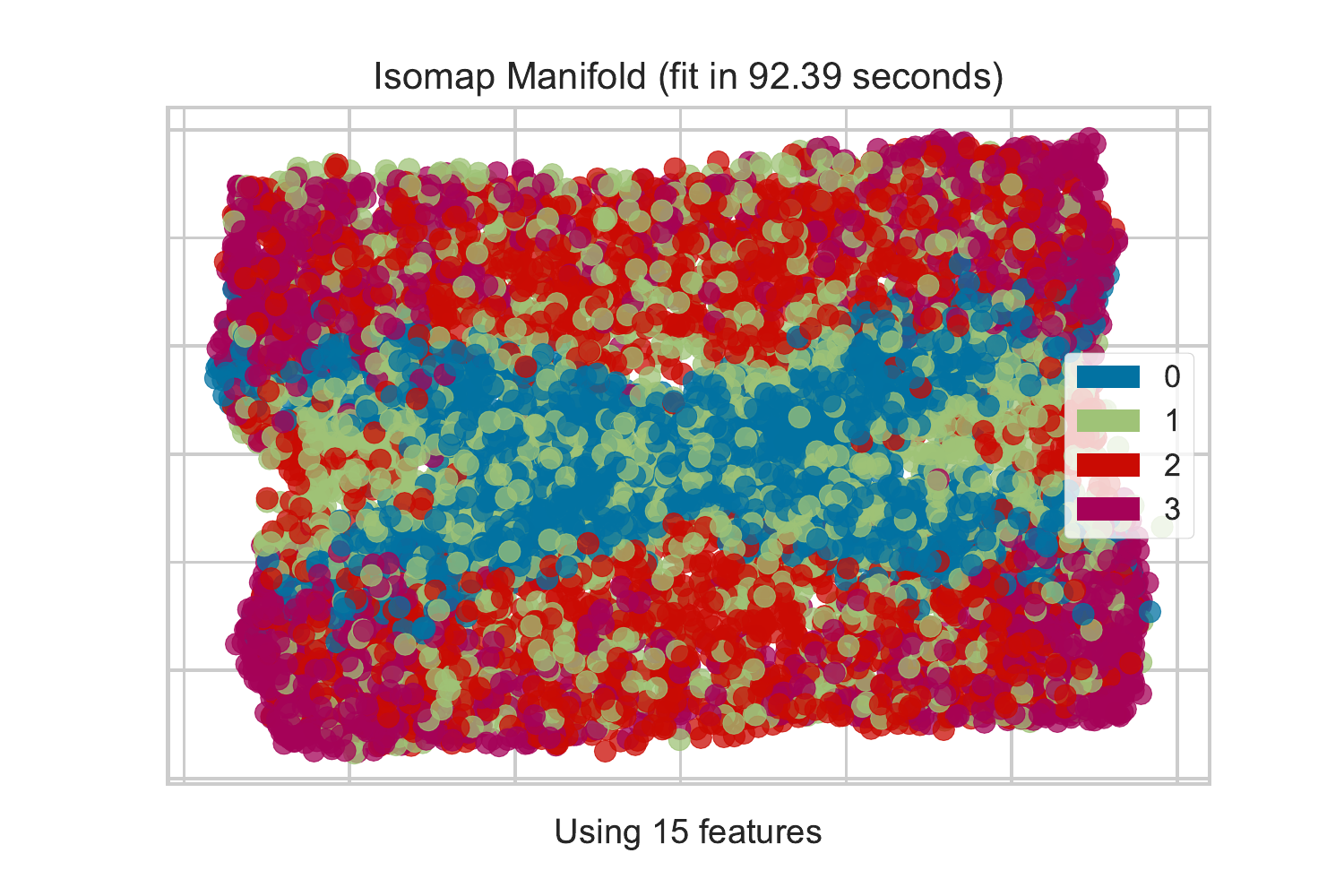}
		\caption{}
	\end{subfigure}
	\begin{subfigure}{7cm}
		\centering
		\includegraphics[width=7cm]{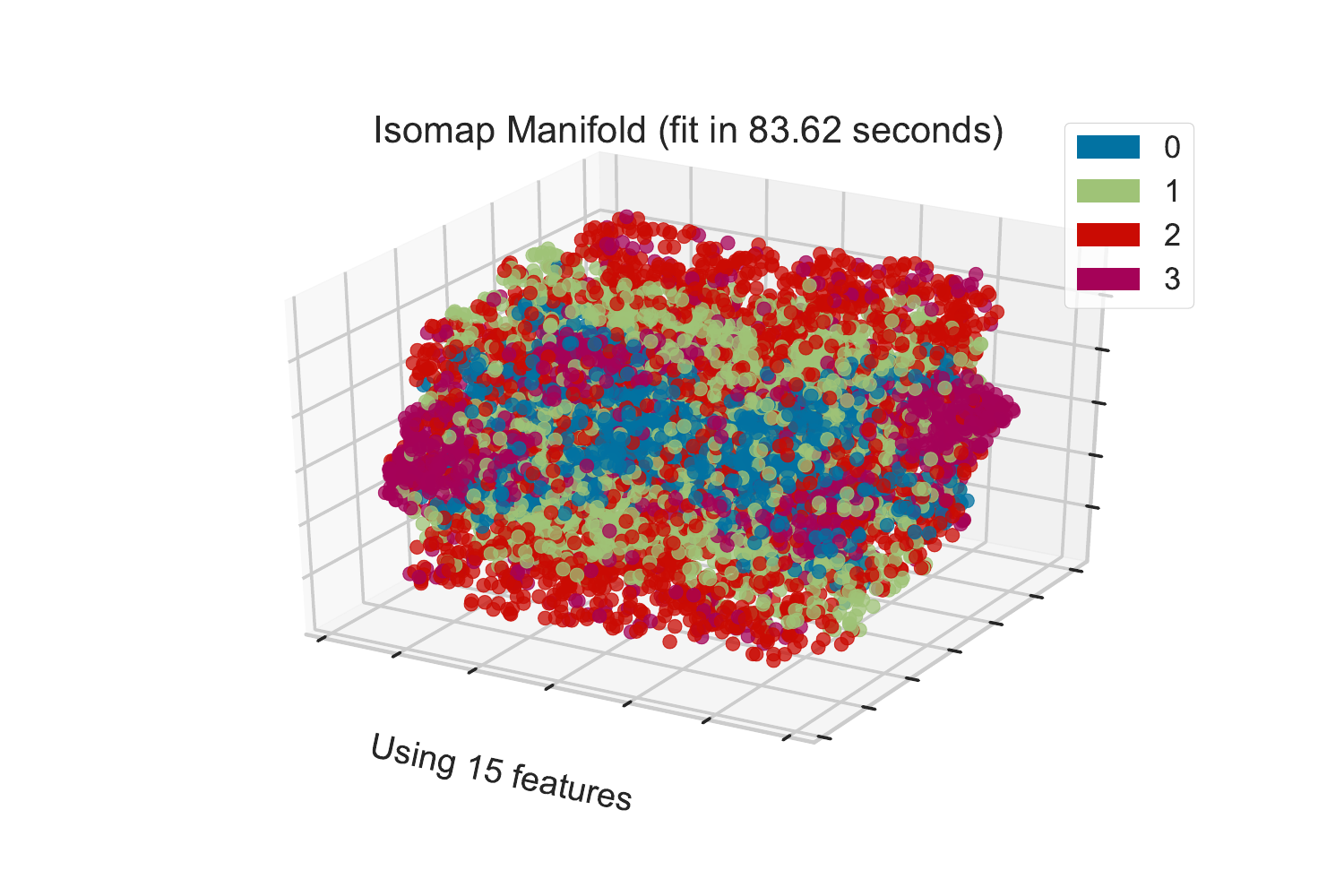}
		\caption{}
	\end{subfigure}
	\caption{(a) The Isomap manifold projection for real input. (b) The Isomap manifold 3d projection for the same input.}\label{cz2z4PCA}
\end{figure}

To understand how the model does the classification, we plot the Isomap embedding manifold projection in Figure \ref{cz2z4PCA}.
Likewise, we can write down the set of the monomials of $P(z,w)$ and check different lopsided cases to find the regions for different $g$ though there is no general formulae to solve the equations and find analytic expressions for different regions in terms of the coefficients as variables. As a sketch, the lopsidedness property leads to the following inequalities:
\begin{equation}
    \begin{split}
         |c_8|&<a_1|c_6|-a_2|c_7|-a_3,\\
         |c_8|&<-b_1|c_6|+b_2|c_7|-b_3,\\
         |c_8|&>d_1|c_6|+d_2|c_7|+d_3,
    \end{split}
\end{equation}
where $a_i,b_i,d_i$ are complicated \emph{positive} expressions in terms of $|c_i|$.

\begin{figure}[h]
	\centering
	\begin{subfigure}{7cm}
		\centering
		\includegraphics[width=7cm]{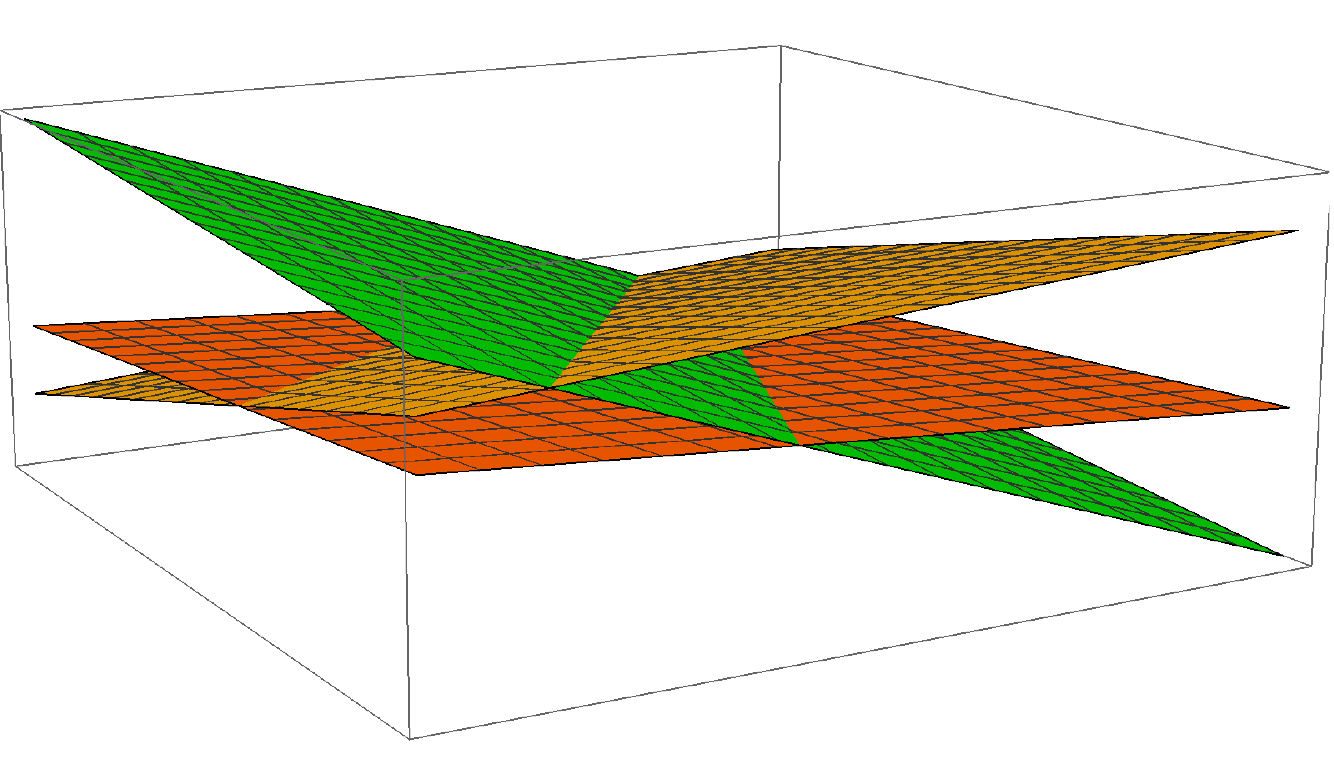}
		\caption{}
	\end{subfigure}
	\begin{subfigure}{7cm}
		\centering
		\includegraphics[width=7cm]{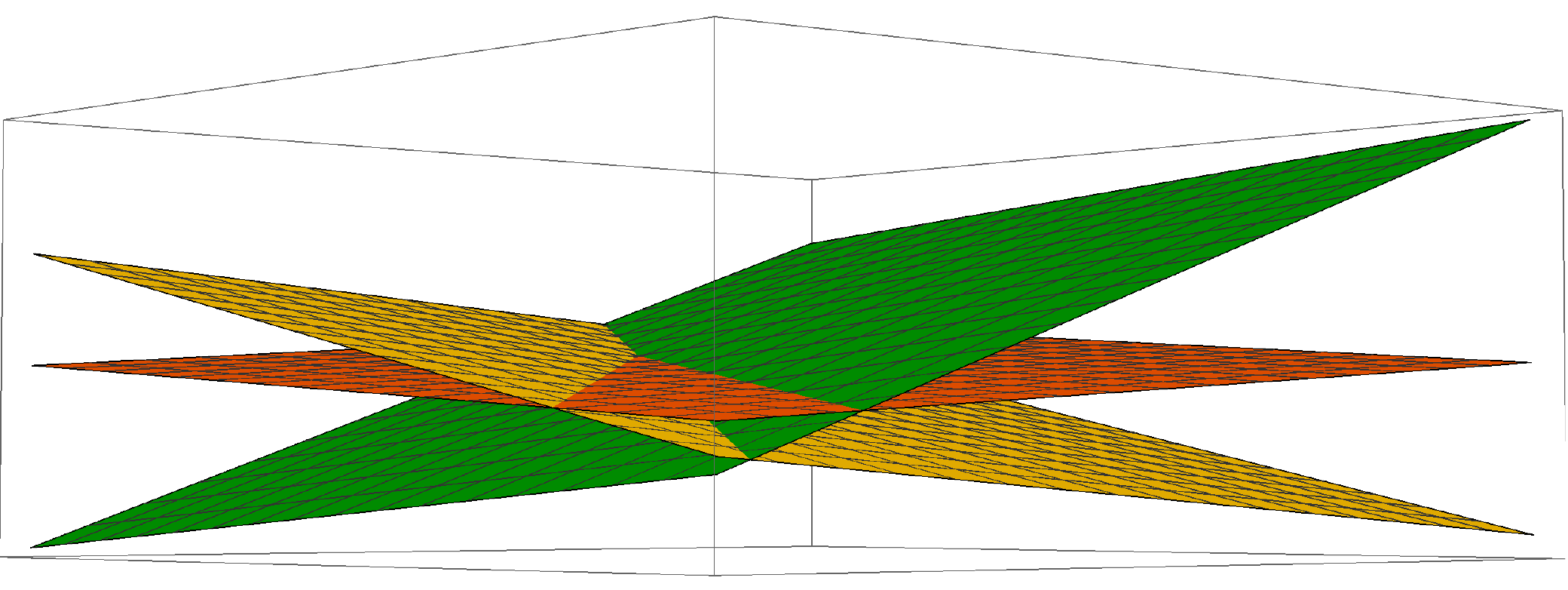}
		\caption{}
	\end{subfigure}
	\begin{subfigure}{7cm}
		\centering
		\includegraphics[width=4cm]{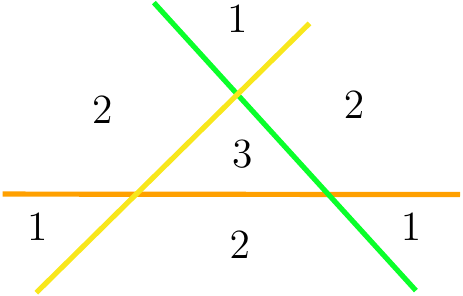}
		\caption{}
	\end{subfigure}
	\begin{subfigure}{7cm}
		\centering
		\includegraphics[width=4cm]{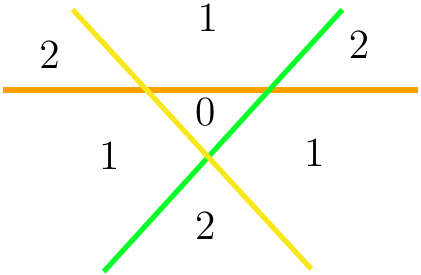}
		\caption{}
	\end{subfigure}
	\caption{The bounds for different $g$'s in $|c_{6,7,8}|$ coordinates. (a) The bounds for different regions viewed from the ``front''. (b) The bounds for different regions viewed from the ``back''. The orange plane corresponds to $|c_8|=d_1|c_6|+d_2|c_7|+d_3$. The green plane is $|c_8|=a_1|c_6|-a_2|c_7|z-a_3$. The yellow one is $|c_8|=-b_1|c_6|+b_2|c_7|-b_3$. (c) The 2d cross section of (a) viewed from the ``front''. (d) The 2d cross section of (b) viewed from the ``back''. The number in each region is the genus.}\label{cz2z4bound}
\end{figure}

Then for fixed coefficients, $g$ is equal to the number of these inequalities being satisfied. As our data is generated with coefficients in certain finite ranges, we can still assume that the changes of $a_i,b_i,d_i$ are small compared to the coloured regions. The schematic is depicted in Figure \ref{cz2z4bound}, where for convenience only the first quadrant is shown.
As we can see, the 2d projection in Figure \ref{cz2z4PCA}(a) is actually viewing the region from the ``side'' (rather than ``front'' or ``back'')\footnote{One consequence of viewing from the ``side'' is that $g=0$ and $g=3$ are very well separated. Indeed, if we only use data with $g=0,3$ for binary classification. The accuracy could easily reach over $0.99$.} while the 3d projection in Figure \ref{cz2z4PCA}(b) is also consistent with Figure \ref{cz2z4bound}.

\begin{figure}[h]
	\centering
	\begin{subfigure}{7cm}
		\centering
		\includegraphics[width=7cm]{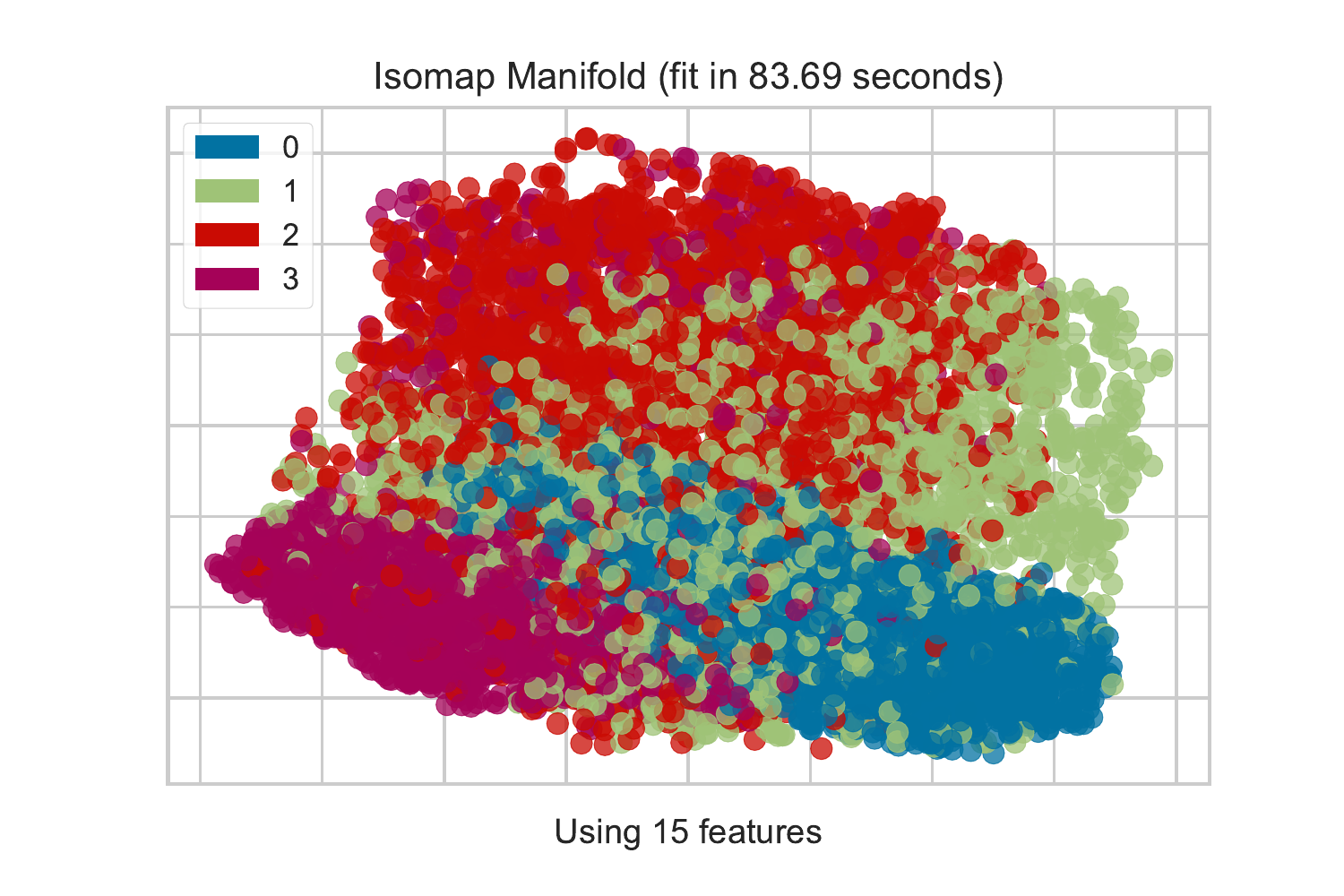}
		\caption{}
	\end{subfigure}
	\begin{subfigure}{7cm}
		\centering
		\includegraphics[width=7cm]{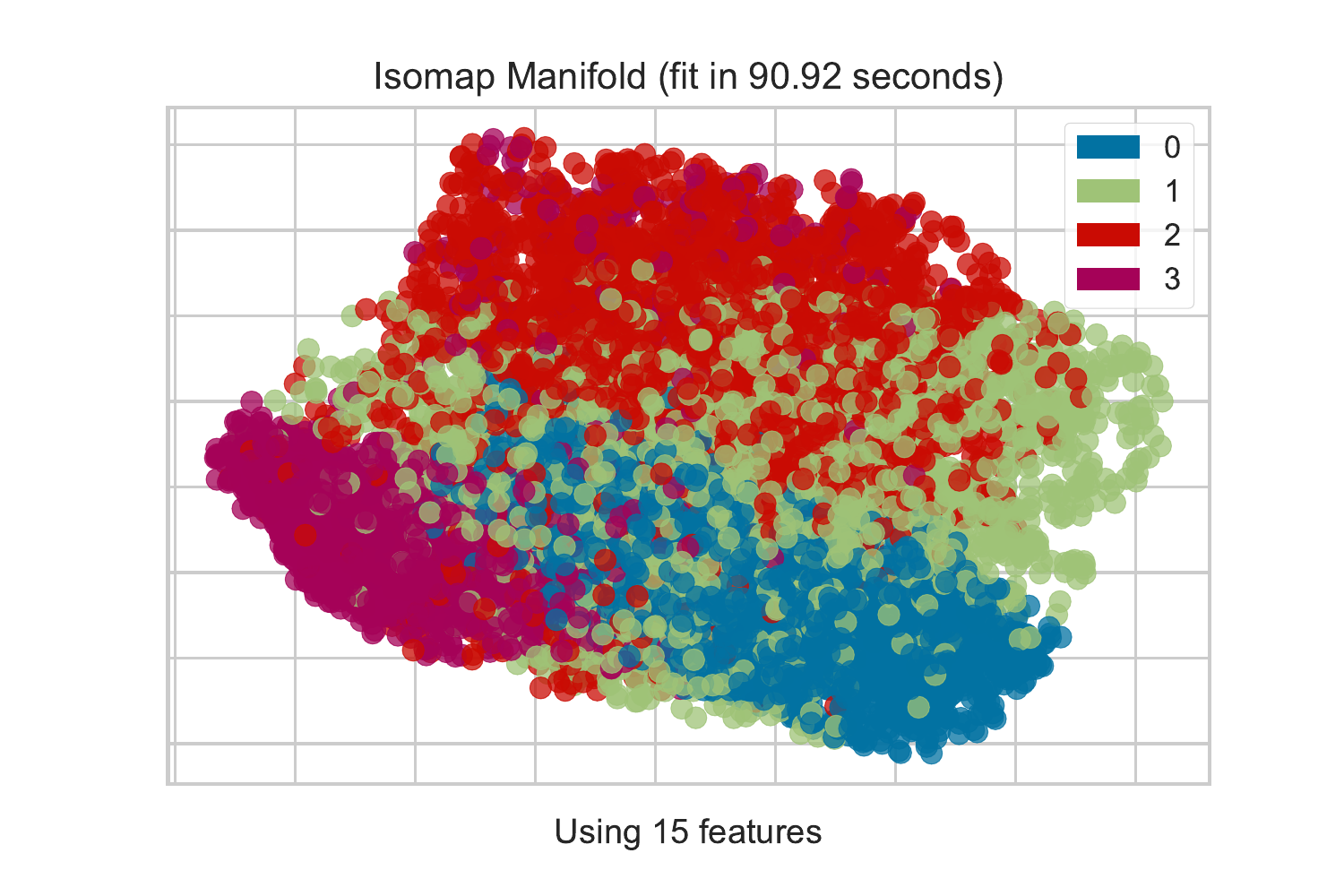}
		\caption{}
	\end{subfigure}
	\caption{(a) The Isomap manifold projection for positive real input. (b) The Isomap manifold projection for positive integer input.}\label{cz2z4PCApositive}
\end{figure}

\begin{figure}[h]
	\centering
		\centering
		\includegraphics[width=7cm]{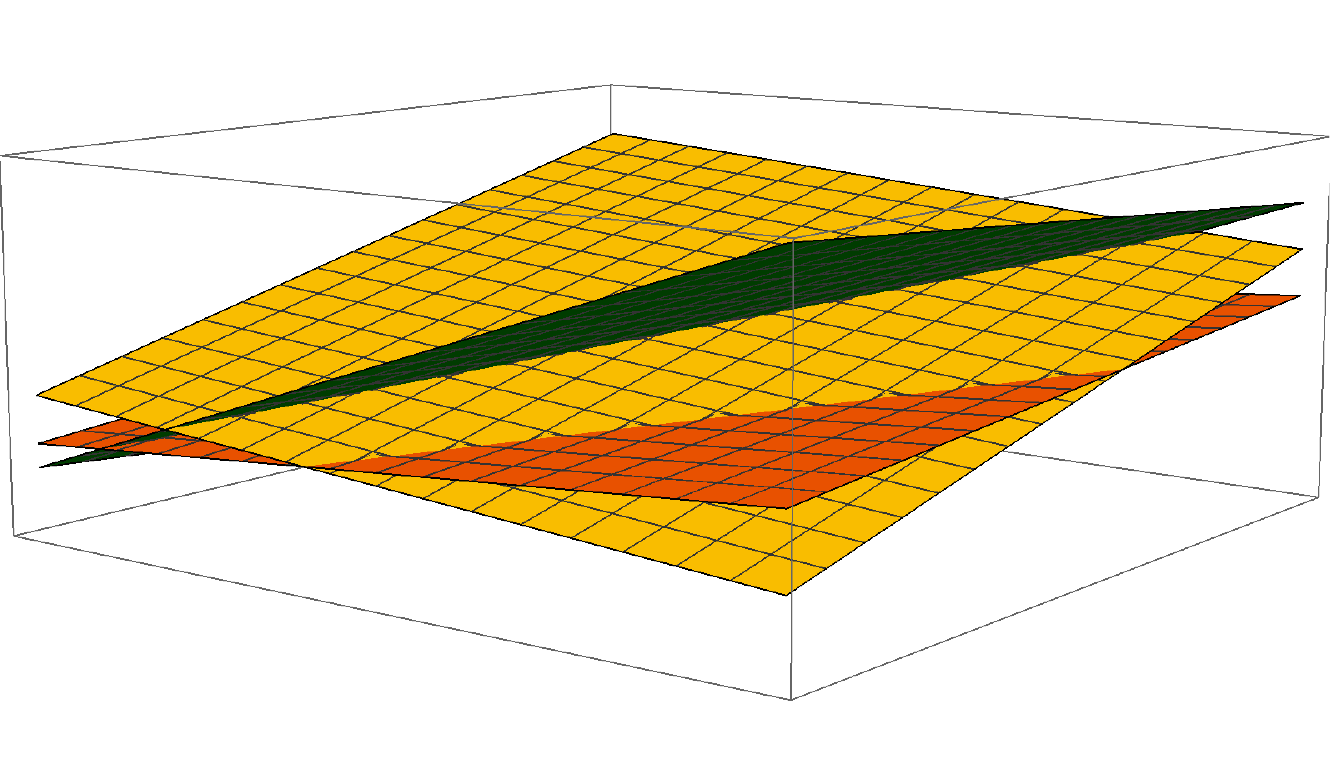}
	\caption{The bounds for $\mathcal{C}/(\mathbb{Z}_2\times\mathbb{Z}_4)$, viewed from ``the side''. This is exactly what the NN gets in Figure \ref{cz2z4PCApositive} (notice that some red and green points are covered by blue and purple points near the middle as the projection is 2d).}
	\label{cz2z4boundpositive}
\end{figure}

For reference, we also plot the Isomap manifold projection for positive real and positive integer coefficients in Figure \ref{cz2z4PCApositive} (a) and (b) respectively.
Again, the machine is viewing the different regions from ``the side''. This agrees with the sketch in Figure \ref{cz2z4boundpositive}, and we see that decision regions are better simulated than before, which therefore gives higher accuracies, as the coefficients are all positive (or equivalently, $|c_i|$).

\subsubsection{Example 4: $K^{4,5,3,2}$}\label{K4532}
For more complicated cases with more interior points, it is harder to analyze the visualizations of different bounds for the genus. Here we briefly discuss this with an example, the so-called $K^{4,5,3,2}$ space whose toric diagram is\footnote{For convenience, the nomenclature for $K^{a,b,c,d}$ follows \cite{Bao:2020kji}.}
\begin{equation}
	\includegraphics[width=2cm]{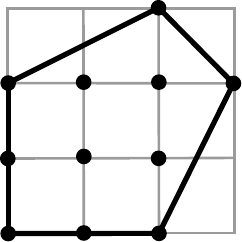}.
\end{equation}
The Newton polynomial is
\begin{equation}
    P(z,w)=c_0+c_1z+c_2z^2+c_3w+c_4zw+c_5z^2w+c_6w^2+c_7zw^2+c_8z^2w^2+c_9z^3w^2+c_{10}z^2w^3.
\end{equation}
Hence, the input is $\{c_0,c_1,\dots,c_{10}\}$, and the classification output is $g=0,1,2,3,4$. We generate a balanced data-set with $\sim25000$ random samples with $c_i\in[-30,30]$. For real coefficients, MLP would only give $0.606(\pm0.003)$ accuracy for 5-fold cross validation while CNN would get $0.671(\pm0.013)$ accuracy. However, if we restrict the coefficients to be positive (or equivalently, taking their absolute values), the accuracy for MLP would be increased to $0.706(\pm0.006)$. Moreover, CNN could reach $0.815(\pm0.006)$ accuracy. It would also be reasonable to expect further improvements of the performance with larger datasets, as well as optimisation over the NN hyperparameters (including the consideration of higher depth networks which can exhibit faster convergence).

As before, we can write the conditions for the genus as
\begin{equation}
    \begin{split}
         |c_4|&>a_1|c_5|+a_2|c_7|+a_3|c_8|+a_4,\\
         |c_5|&>b_1|c_4|+b_2|c_7|+b_3|c_8|+b_4,\\
         |c_7|&>d_1|c_4|+d_2|c_5|+d_3|c_7|+d_4,\\
         |c_8|&>e_1|c_4|+e_2|c_5|+e_3|c_7|+e_4,
    \end{split}\label{K4532genus}
\end{equation}
where $a_i,b_i,d_i,e_i$ are again complicated expressions in $c_i$. For our finite data-set, we may assume that they have small changes compared to the regions of data points and hence think of them as hyperplanes in the 4-dimensional system. Then for fixed $|c_i|$, $g$ is equal to the number of these inequalities being satisfied. It is difficult to visualize the projections as in the previous subsections. Nevertheless, let us still use PCA to project the input vectors to 3d.

\begin{equation}
	\includegraphics[width=8cm]{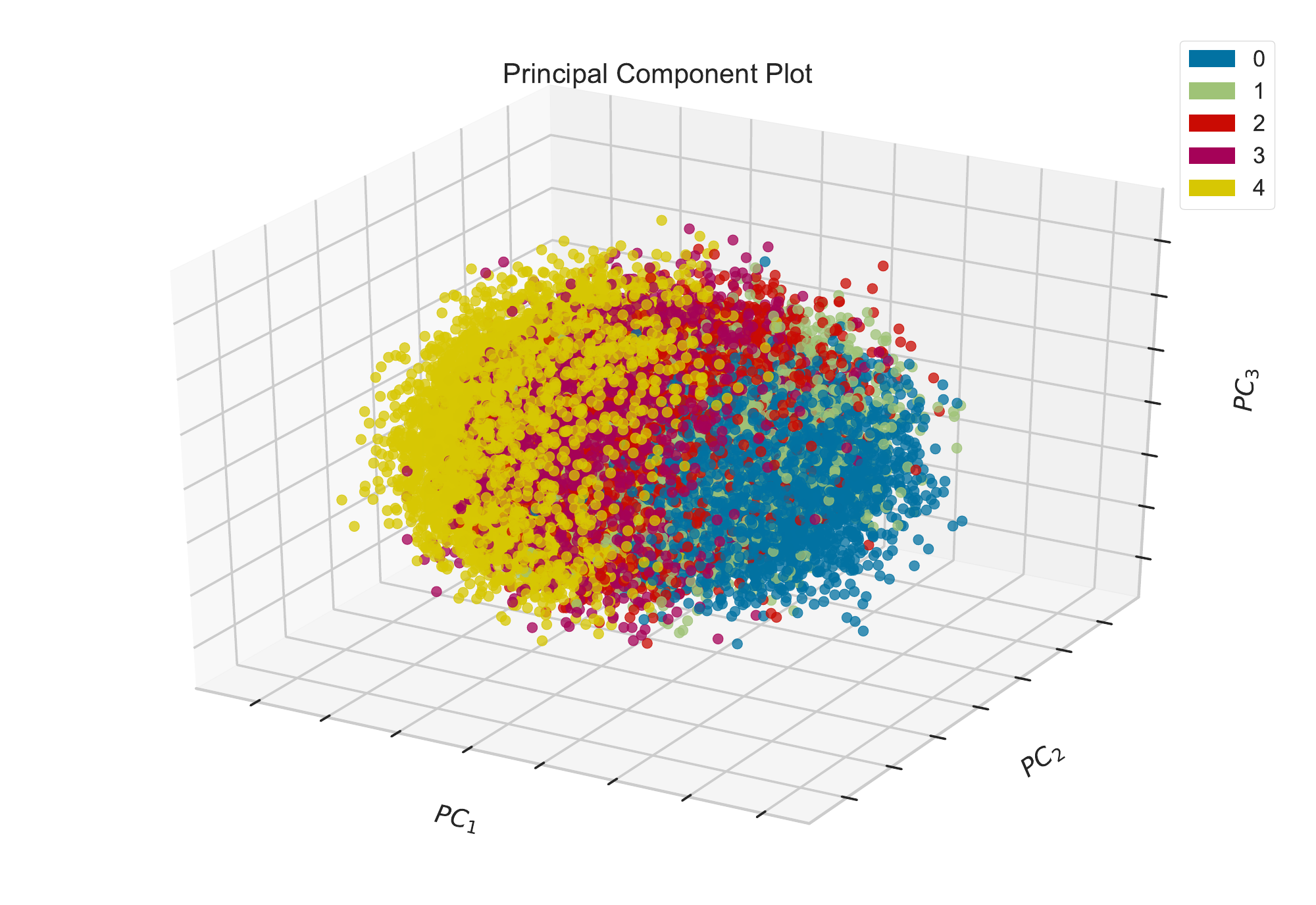}.\label{K4532PCA0}
\end{equation}
Although it is hard to tell what the three principal components are, we can still see that points in different colours are arranged from left to right (with mixings).

\begin{figure}[H]
	\centering
	\begin{subfigure}{5cm}
		\centering
		\includegraphics[width=5cm]{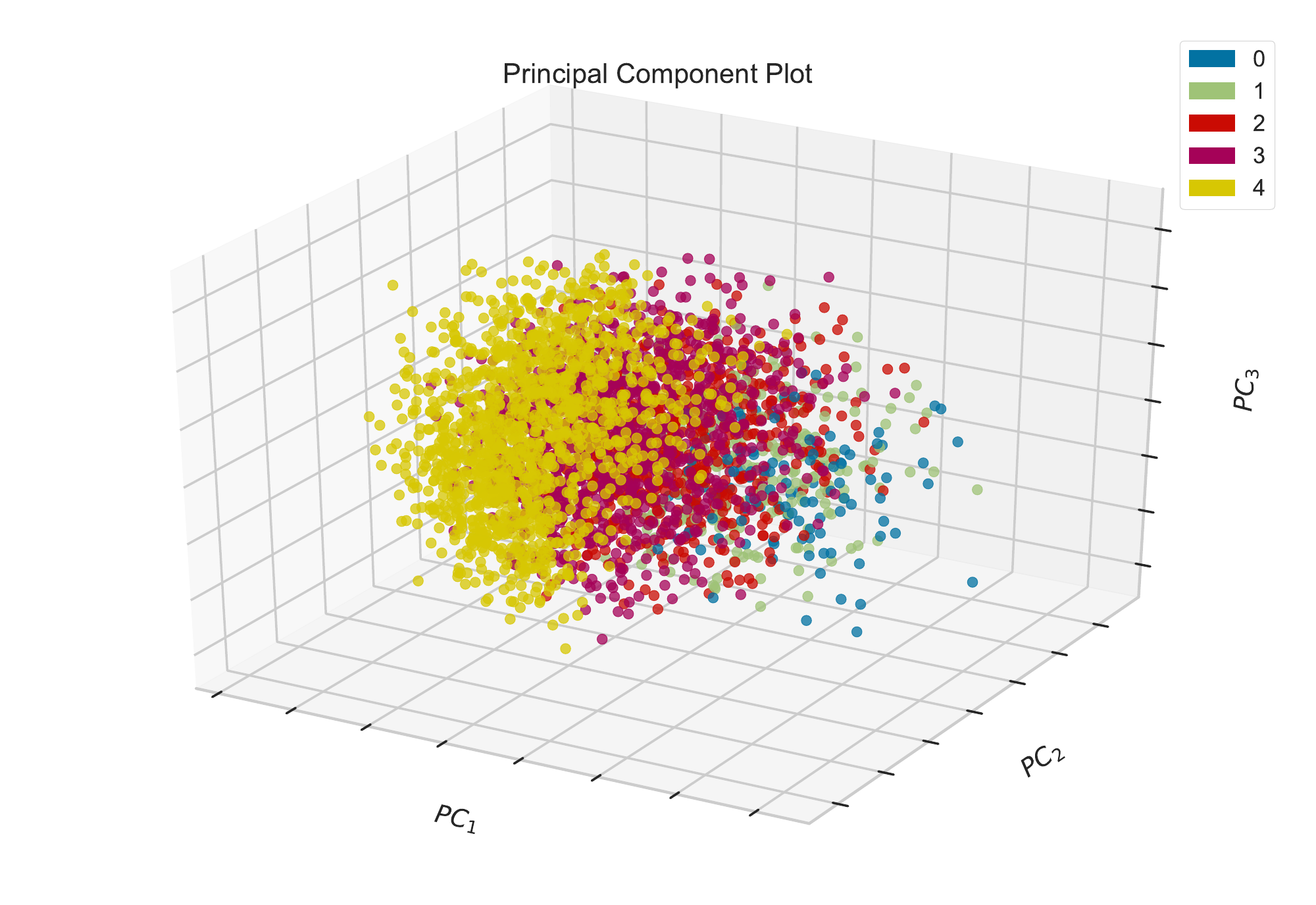}
		\caption{}
	\end{subfigure}
	\begin{subfigure}{5cm}
		\centering
		\includegraphics[width=5cm]{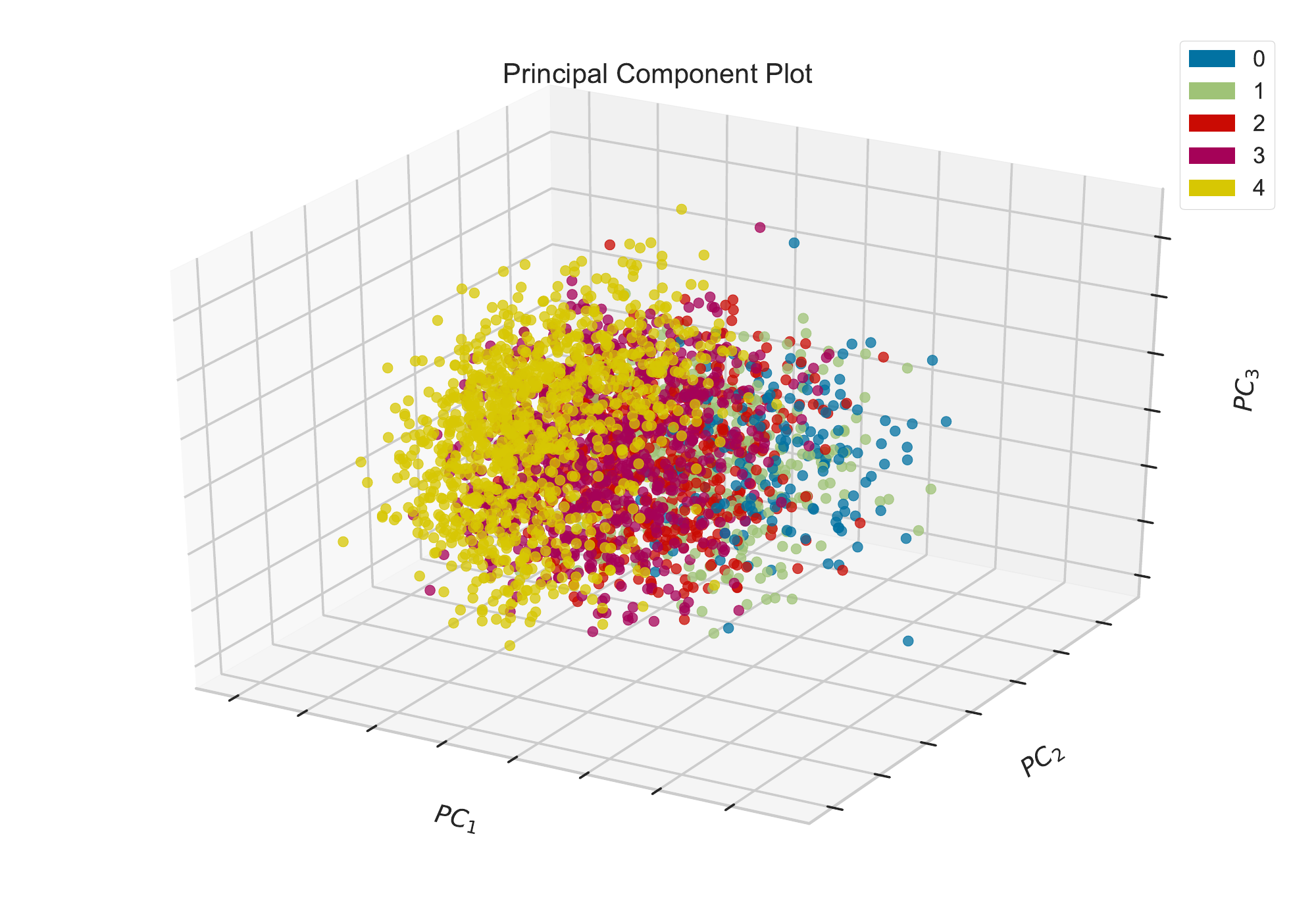}
		\caption{}
	\end{subfigure}
	\begin{subfigure}{5cm}
		\centering
		\includegraphics[width=5cm]{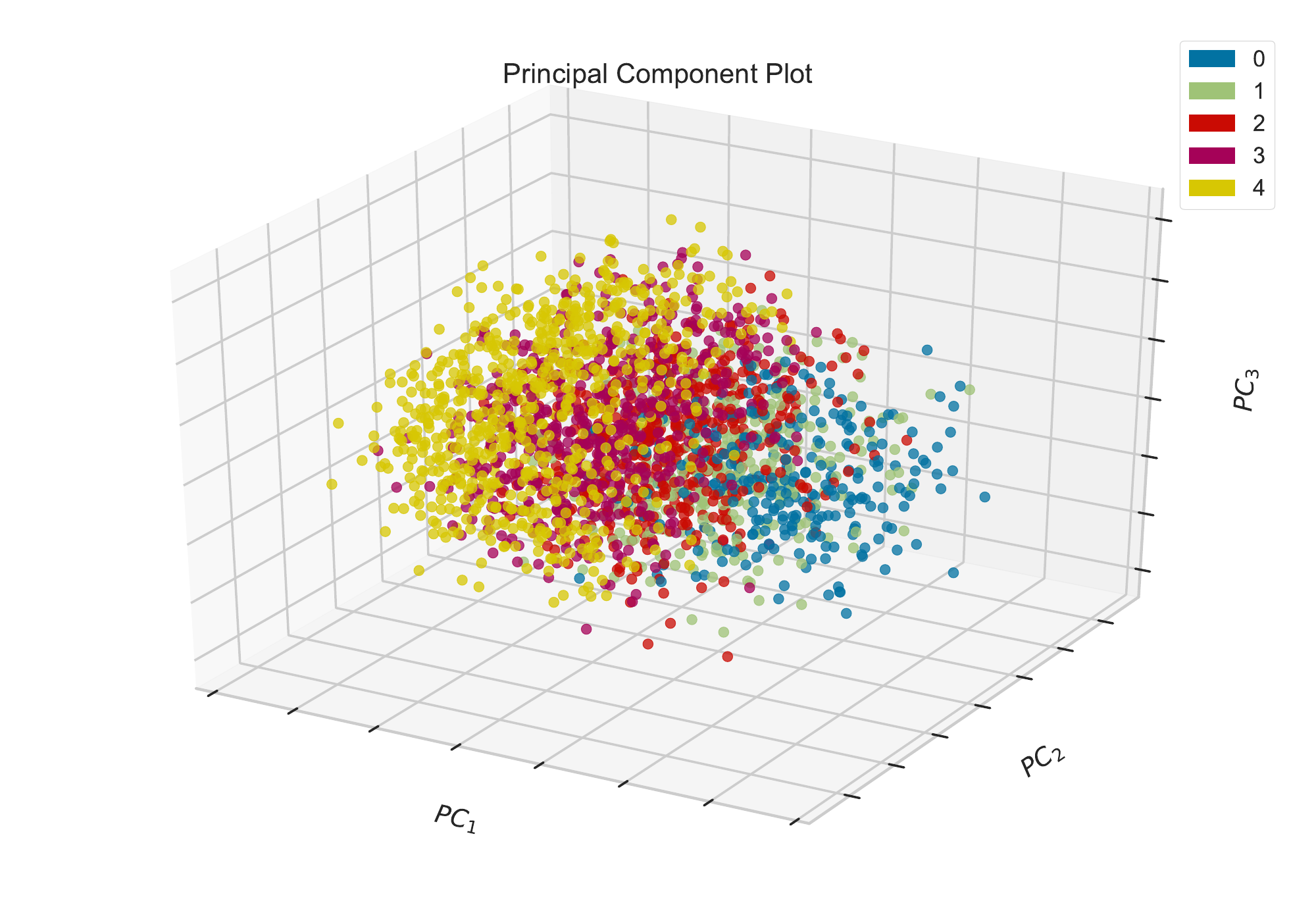}
		\caption{}
	\end{subfigure}
	\begin{subfigure}{5cm}
		\centering
		\includegraphics[width=5cm]{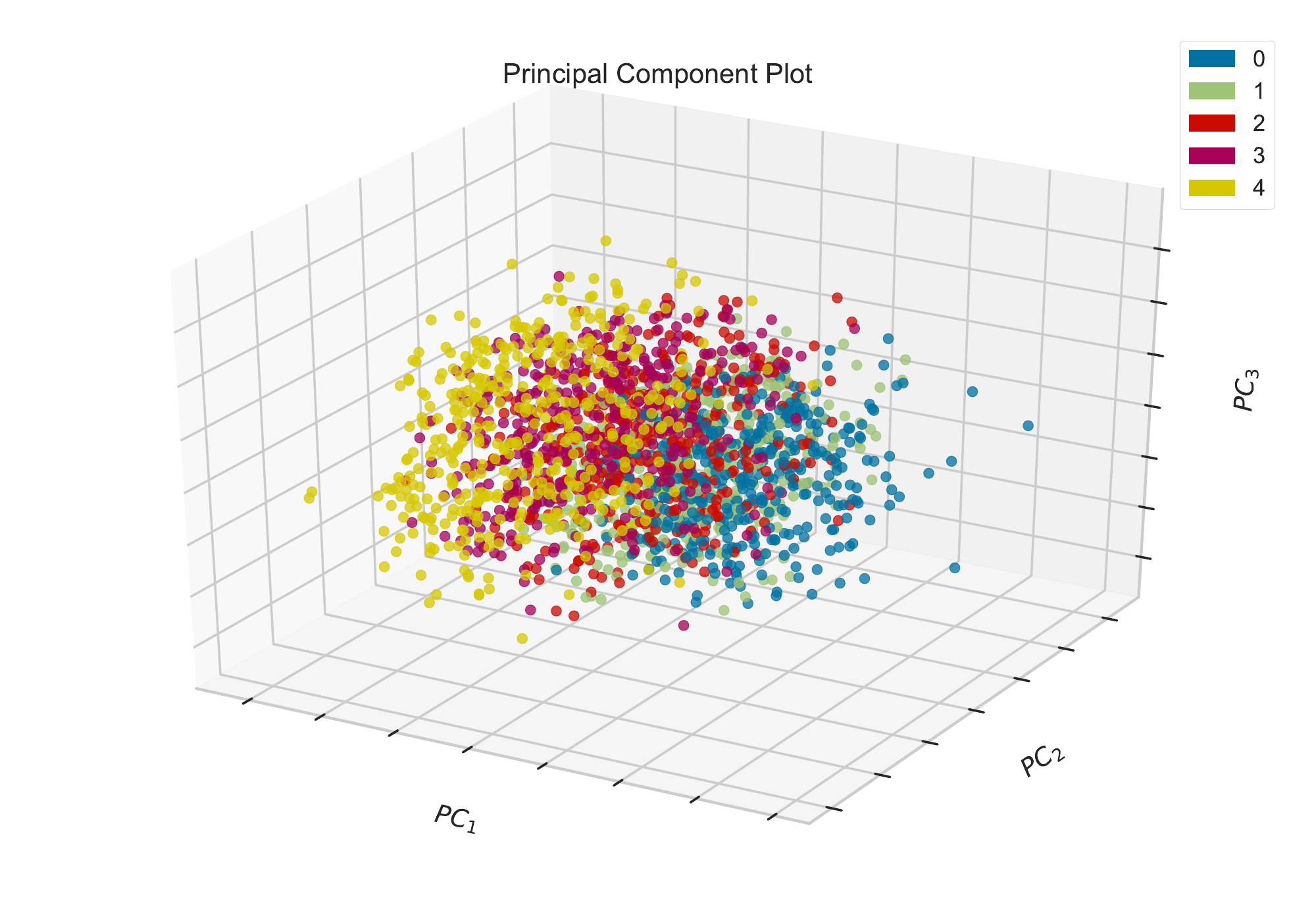}
		\caption{}
	\end{subfigure}
	\begin{subfigure}{5cm}
		\centering
		\includegraphics[width=5cm]{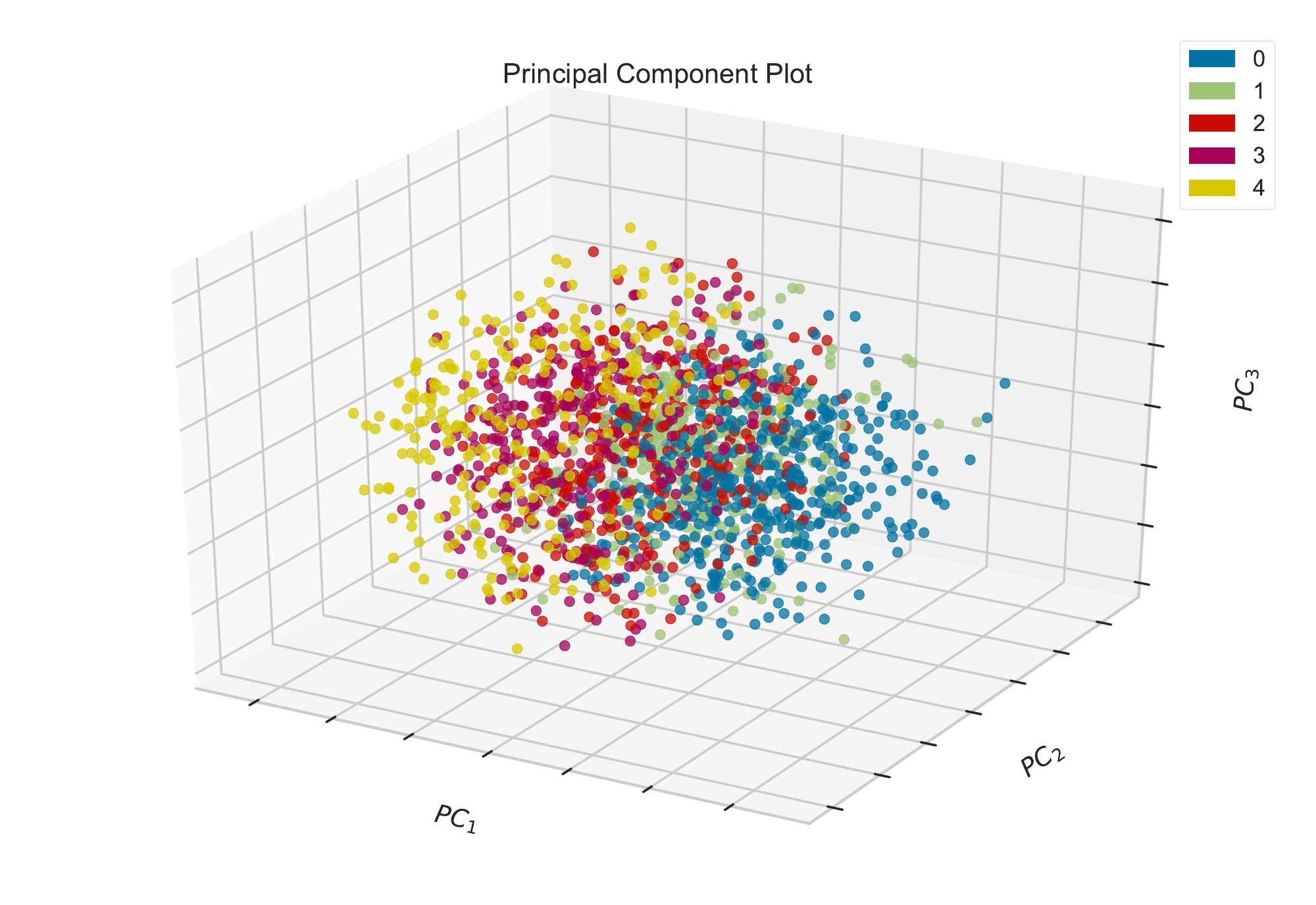}
		\caption{}
	\end{subfigure}
	\begin{subfigure}{5cm}
		\centering
		\includegraphics[width=5cm]{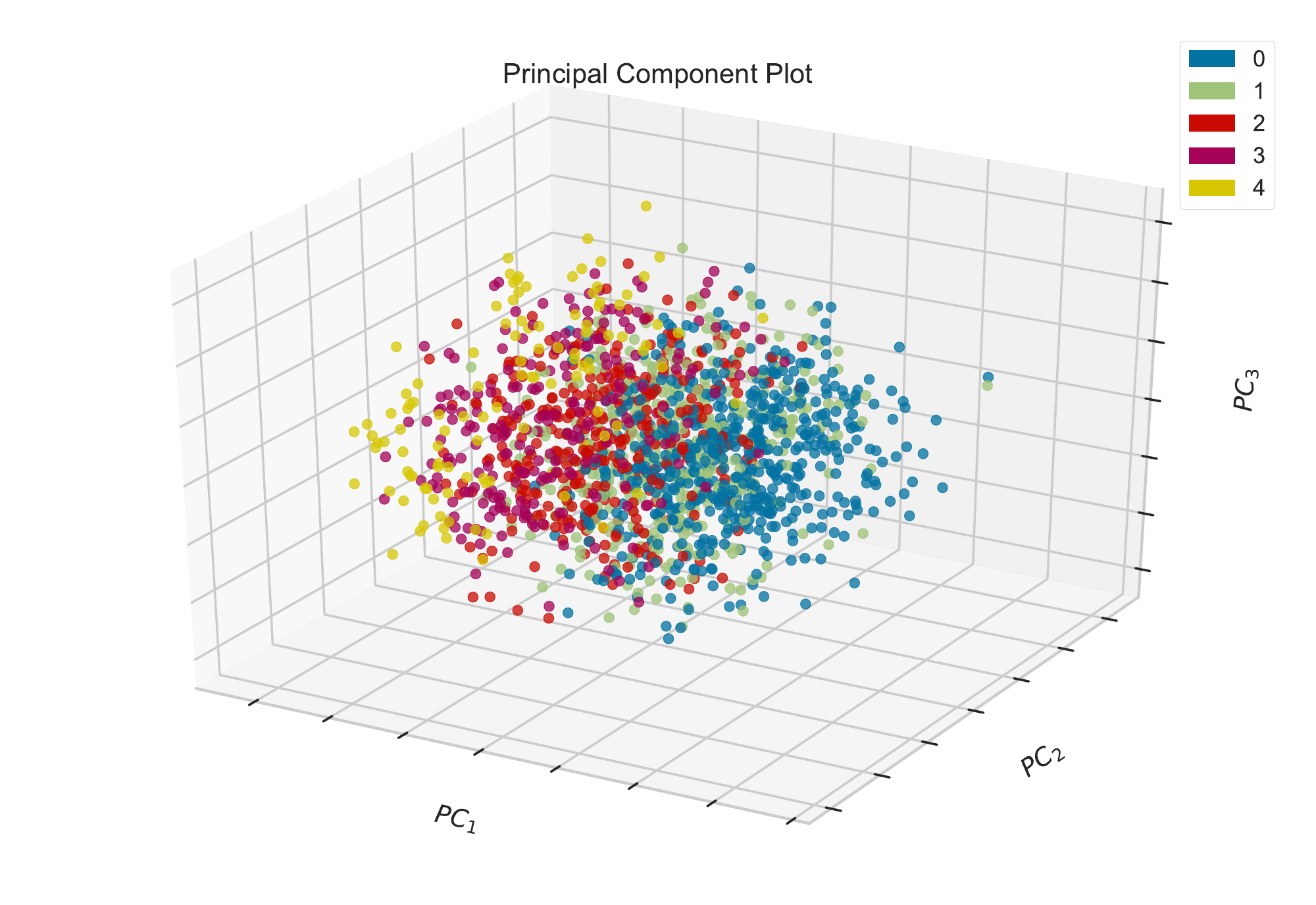}
		\caption{}
	\end{subfigure}
	\caption{The PCAs of different slices for (a) $27<|c_8|\leq30$, (b) $24<|c_8|\leq27$, (c) $21<|c_8|\leq24$, (d) $18<|c_8|\leq21$, (e) $15<|c_8|\leq18$ and (f) $12<|c_8|\leq15$.}\label{K4532PCA}
\end{figure}

Another way to visualize in 3d is to consider 3-dimensional slices in the 4-dimensional system. For instance, we can fix $|c_8|$ and different plots should correspond to different cross sections of the 4d ``plot''. In Figure \ref{K4532PCA}, we plot several slices by restricting $|c_8|$ to small ranges as our data consists of vectors with random real entries. We find that the distribution of coloured points in each slice looks similar to \eqref{K4532PCA0}, where the numbers of points with larger $g$ decrease and those with smaller $g$ increase when we decrease the value of $|c_8|$.
Of course, it is also possible to project an input vector to one with 4 components. We will discuss this in \S\ref{projections}.

\subsection{Reproduction of the Genus: Interpretable ML}\label{reproduce}
Encouraged by the success of the above experiments, it is natural to ask whether we can recover the bounds/conditions for genus using the machine-learning models we built.
The analogous situation in computing cohomology groups for bundles over algebraic surfaces was considered in \cite{Brodie:2019dfx}, where polynomial boundaries in the space of bundle degrees were detected with NNs.
Likewise, hypersurfaces which separate the space of simple finite groups from the non-simple were detected in \cite{He:2019nzx}.
Can our complicated hole-boundaries of amoebae be detected by ML?

\subsubsection{From Projections}\label{projections}
In \S\ref{n1}, we presented many plots where the input data is projected to some lower dimensional representations under certain manifold embeddings. For instance, let us take Figure \ref{F0PCA}(b) for $F_0$ where MDS manifold projection was applied. Each point $(x,y)$ in the plane represents one input vector as $(c_1,c_2,c_3,c_4,c_5)\mapsto(x,y)$. By checking the coordinates of $(x,y)$, it is possible to recover the conditions in \eqref{gF0}, i.e.,\\

$g=\begin{cases}
		0,&|c_5|\leq a\\
		1,&|c_5|>a
	\end{cases}$
with $a := |c_1c_3|^{1/2}+2|c_2c_4|^{1/2}$, as follows.

\begin{figure}[H]
	\centering
	\begin{subfigure}{7cm}
		\centering
		\includegraphics[width=7cm]{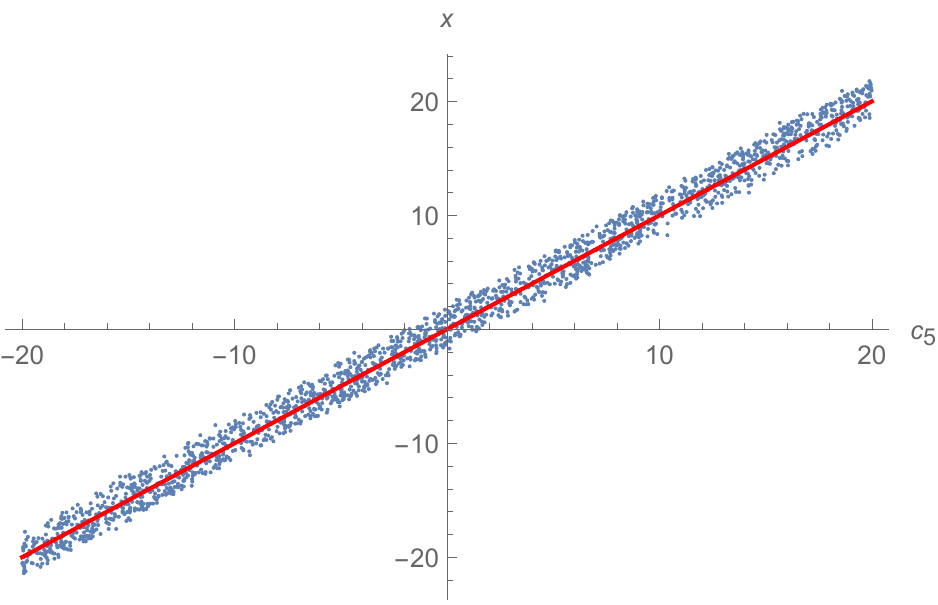}
		\caption{}
	\end{subfigure}
	\begin{subfigure}{7cm}
		\centering
		\includegraphics[width=7cm]{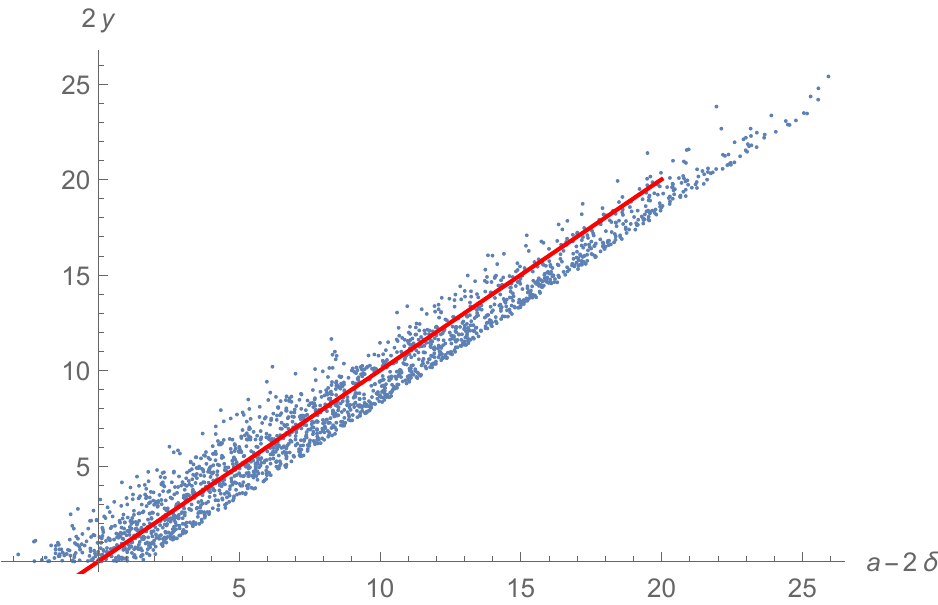}
		\caption{}
	\end{subfigure}
	\caption{For the projection $(c_1,c_2,c_3,c_4,c_5)\mapsto(x,y)$ for $F_0$, we plot
	(a) $x$ versus $c_5$ and
	(b) 2$|y|$ versus $a - 2 \delta$.
	}\label{f:F0bound}
\end{figure}

We plot $x$ against $c_5$ in Fig.~\ref{f:F0bound}(a) and find that the average of $(c_5-x)$ is $-0.172 \pm 1.054$. We also draw the line $x=c_5$ is in red to show a good fit. Indeed, we see that $x$ is essentially recovering the left hand side of the above inequalities.
For $y$, it is more complicated and we find that there is a nice fit with $|c_{1,2,3,4}|$:
\begin{equation}
    |y_\text{fit}|=0.141|c_1c_3|+0.166|c_2c_4|+2.411.\label{yfitF0}
\end{equation}
In fact, this fit actually contains a typical linear approximation of square roots, $\sqrt{m}\approx(0.1k+1.2)\times10^n$ for any real $k\in[1,100)$ and $n\in\mathbb{Z}$ such that $m=k\times10^{2n}$. Here, $c_{1,2,3,4}$ are random reals generated in the range $[-5,5]$ and the right hand side is approximately $(0.1|c_1c_3|+1.2+0.1|c_2c_4|+1.2)$. This agrees with the linear approximation for square roots when $m=k$ and $n=0$. In other words,
\begin{equation}
    |y|\approx\sqrt{|c_1c_3|}+\sqrt{|c_2c_4|}=a/2.
\end{equation}
We show this in Fig.~\ref{f:F0bound}(b).
The blue points are the pairs $(a-2\delta,2|y|)$ and the line $2|y|=a-2\delta$ is in red, where we have defined $\delta:=|y_\text{fit}|-|y|$.
Comparing the fitted results of $(x,y)$ and Figure \ref{F0PCA}, we find that this indeed gives the bound $|c_5|=a$ for $F_0$.

\begin{figure}[h]
	\centering
	\begin{subfigure}{7cm}
		\centering
		\includegraphics[width=7cm]{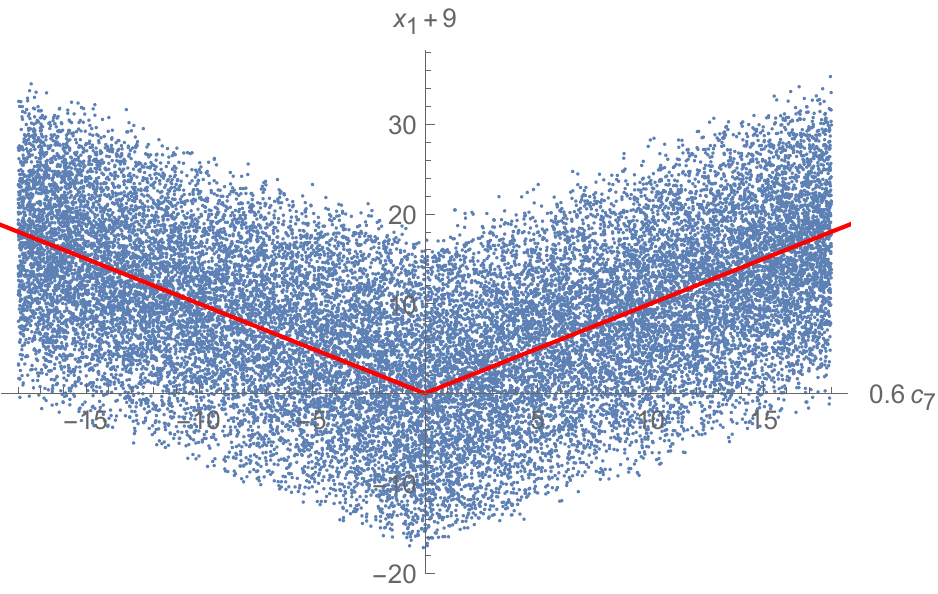}
		\caption{}
	\end{subfigure}
	\begin{subfigure}{7cm}
		\centering
		\includegraphics[width=7cm]{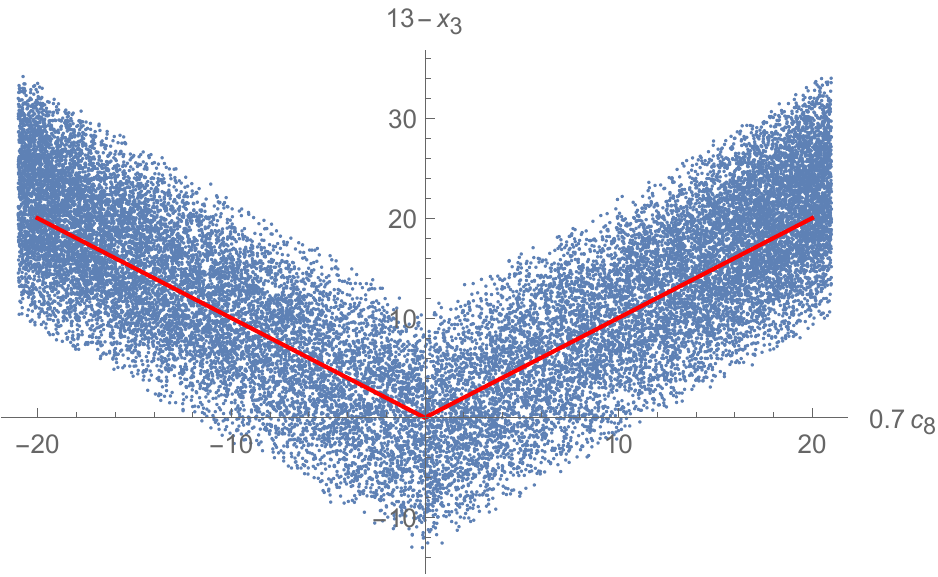}
		\caption{}
	\end{subfigure}
	\begin{subfigure}{7cm}
		\centering
		\includegraphics[width=7cm]{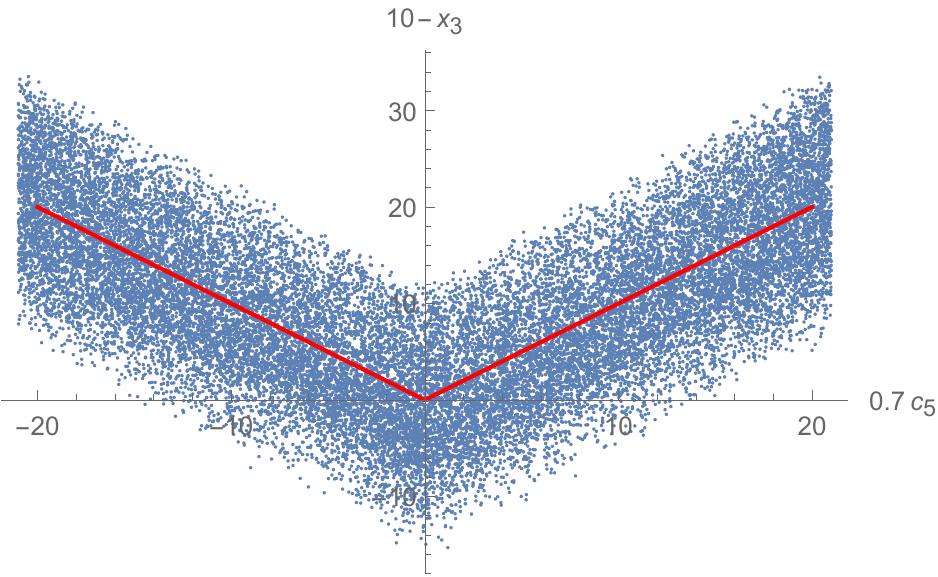}
		\caption{}
	\end{subfigure}
	\begin{subfigure}{7cm}
		\centering
		\includegraphics[width=7cm]{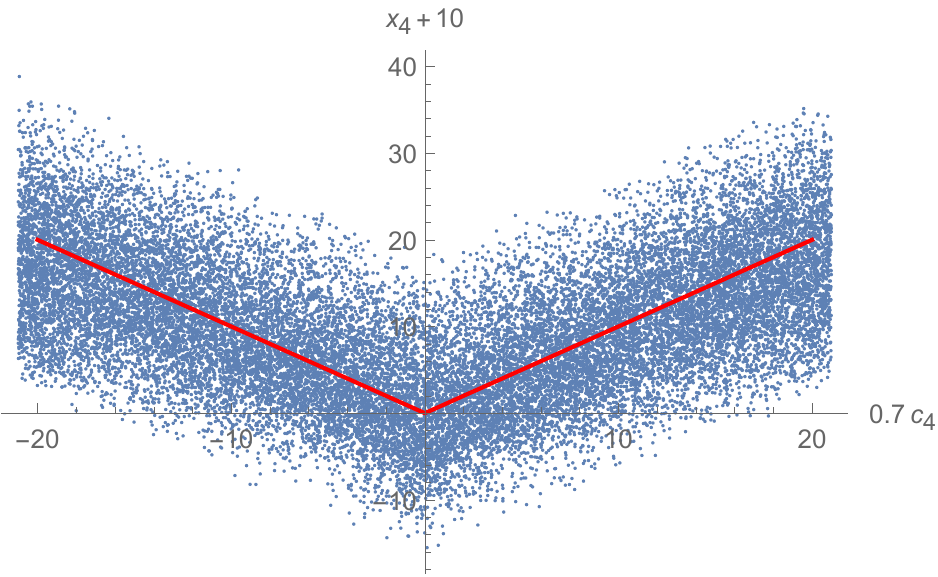}
		\caption{}
	\end{subfigure}
	\caption{The fits of $x_i$'s in terms of $c_j$'s.}\label{K4532x1234}
\end{figure}

Let us now consider a more complicated example, that is, $K^{4,5,3,2}$ discussed in \S\ref{K4532}. Previously, we projected the input vectors to 3d in order to visualize the data. Here, we use PCA to project the 11-component vectors to 4-component ones and see how well this realizes \eqref{K4532genus}. Recall that our input is $(c_0,c_1,\dots,c_{10})$. The projection then maps it to the vector $(x_1,x_2,x_3,x_4)$. The simple fit for each $x_i$ is shown in Figure \ref{K4532x1234}.
In each plot, the horizontal axis is $k_1c_j$ and the vertical axis is $k_2x_i+k_3$ for some constants $k_{1,2,3}$. We also plot the line $Y=|X|$ in red. It is thus not hard to see that we have the fit $k_2x_i+k_3=|k_1c_j|$; which agrees with the four components/variables $|c_{4,5,7,8}|$ in \eqref{K4532genus}. Since the parameters $a_i,b_i,d_i,e_i$ in \eqref{K4532genus} are complicated expressions of $|c_i|$ which are approximately viewed as ``constants'', the V-shape distributions of the points $(k_1c_j,k_2x_i+k_3)$ are thickened strips of the red line.

\subsubsection{Learning from Weights}\label{weights}
From PCA let us return to NNs.
Given the good performance of latter, can we extract information from the NN (hyper-)parametres? 
In this subsection, we show that the bounds/conditions for the genus can be very well approximated by the weight matrices and biases in the NN structure. 
Indeed, an MLP with $n$ hidden layers with activation function $f_{i+1}(W_i\bm{x}_i+\bm{b}_i)$ at the $(i+1)^\text{th}$ layer, where $\bm{x}_i$ is the output from the $i^\text{th}$ layer and $W_i$, $\bm{b}_i$ are the weight matrix and bias respectively. Then the composition of $f_i$'s would give an approximate expression (in the spirit of the universal approximation theorems of NNs) in terms of the coefficients $c_j$ with low calculation cost to compute~$g$.

As an example, let us again illustrate with $F_0$. With input of the form $\{c_1,c_2,c_3,c_4,c_5\}$, it turns out that even an MLP with two hidden layers and three neurons at each layer would give over $95\%$ accuracy. It turns out that $f_\text{output}(W_2\bm{x}_2+\bm{b}_2)$ would always give a non-negative number. When this number is zero, $g$ is classified to be 0, and $g$ is 1 for the number being positive, as we will now see.

Now that we have low dimensional projection of the data, we can use this projected coordinates $(x,y)$ as input and analyze the structure of MLP. This allows us to further simplify the MLP such that only one hidden layer with four neurons can give over $95\%$ accuracy. The precise entries of weight matrices and bias vectors may vary every time we train a new model, but these values are quite stable and have very small changes. For instance, we have
\begin{equation}
    W_1=\begin{pmatrix}
    0.4667&1.4986&-0.3313&-1.3031\\
    -1.3543&-1.0408&-0.5699&-1.5610
    \end{pmatrix}^\text{T}
\end{equation}
and
\begin{equation}
    W_2=\begin{pmatrix}
    -1.1199&0.6904&0.6216&0.2827
    \end{pmatrix}
\end{equation}
for the weight matrices. The bias vectors are
\begin{equation}
    \bm{b}_1=(-0.0349,0.0439,0.3044,-0.5884)^\text{T},~\bm{b}_2=(-4.1003).
\end{equation}
Again, for some input $(x,y)$, this would give a non-negative number
\begin{equation}
    p=W_2\cdot\text{ReLU}(W_1\cdot(x,y)^\text{T}+\bm{b}_1)+\bm{b}_2
\end{equation}
with $g=0$ only when $p\leq0$. In other words, $g=\theta(p)$ where $\theta$ is the Heaviside function. This can actually be visualized as
\begin{equation}
    \includegraphics[width=9cm]{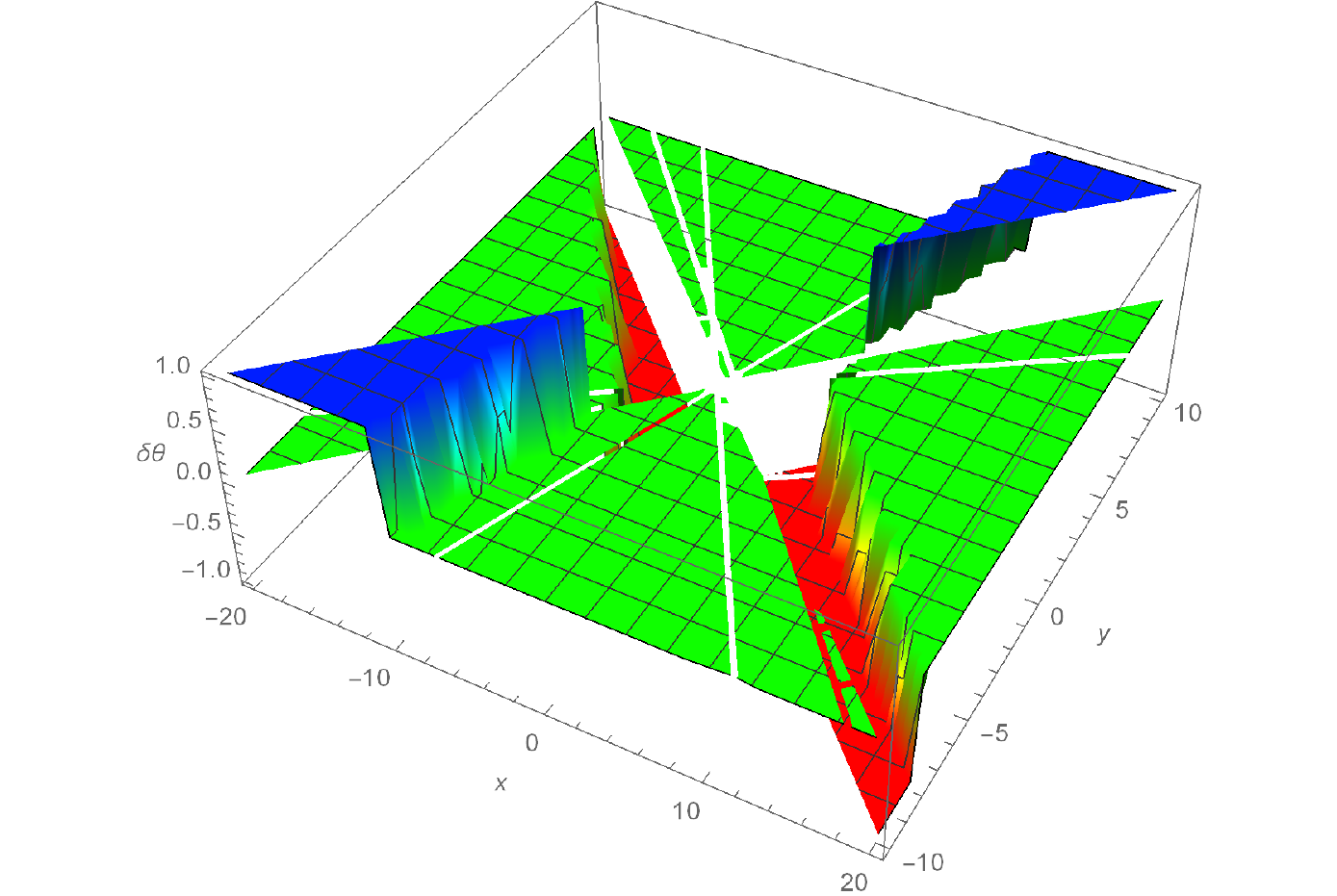},
\end{equation}
where the two horizontal axes are $(x,y)$ and the vertical axis is $\delta\theta:=\theta(p)-\theta(|x|-2|y|)$. Notice that $\theta(|x|-2|y|)=|x|-2|y|=g$ when $c_i$'s are perfectly projected to $|c_5|$ and $\left(\sqrt{|c_1c_3|}+\sqrt{|c_2c_4|}\right)$. As we can see, this gives a precise prediction for $g$ in the green area (assuming that $x$ and $y$ are perfect projections). There are only four slices in blue or red which are not well approximated.

Incidentally, a layer with four neurons is possible to always give $\delta\theta=0$ in principle. This can be realized, for example, when
\begin{equation}
    W_1'=\begin{pmatrix}
    0&1&0&-1\\
    -1&0&1&0
    \end{pmatrix}^\text{T},~
    W_2'=\begin{pmatrix}
    -2&1&-2&1
    \end{pmatrix}
\end{equation}
and $\bm{b}'_{1,2}$ vanish since $|q|=\text{ReLU}(q,0)+\text{ReLU}(-q,0)$ for any real $q$.
Nevertheless, $W_{1,2}$ and $\bm{b}_{1,2}$, together with the projection to $(x,y)$, can already give a nice approximated expression for $g$:
\begin{equation}
    \begin{split}
        &g=\theta(p),~p=W_2\cdot\text{ReLU}(W_1\cdot(x,y)^\text{T}+\bm{b}_1)+\bm{b}_2,\\
        &|x|=|c_5|,~|y|=0.1|c_1c_3|+0.1|c_2c_4|+2.4.
    \end{split}
\end{equation}
We can also use $(c_1,c_2,c_3,c_4,c_5)$ as input directly without any projections. This would give some approximated expression with certain $W_i$ and $\bm{b}_i$ (though it would be hard to visualize as above), and from the results in \S\ref{n1}, we know that this expression approximates $g$ with high accuracy. Likewise, for much more complicated cases, it is also possible to have well-approximated expressions for genus with low calculation cost given sufficient layers and neurons.

\subsection{Lopsided Amoebae: $n>1$}\label{largern}
In the previous section we have shown how ML works successfully with lopsided amoebae for $n=1$ where $\tilde{P}_1 = P$ and we directly addressed the amoeba $\mathcal{A}_P$.
Recall that in practice, it is more efficient to consider $\mathcal{L}\tilde{P}_n$ which approximates $\mathcal{A}_P$ for large enough $n$. What can ML say about the $n>1$ cases?

For any finite $n$, $\mathcal{LA}_{\Tilde{P}_n}$ is always a superset of $\mathcal{A}_P$ which further includes points of distance $\epsilon<(\log(n)+\log(8c))/n$ to the boundary of $\mathcal{A}_P$. In particular, for $n=1$ discussed in the previous section, the genus would be different from the one for $\mathcal{A}_P$ if the ``size'' of the genus is smaller than $\log(8c)$. For example, the amoeba for $P=z+w-z^{-1}-w^{-1}+1$ is of genus one as shown in Fig.~\ref{egAmoeba}.
\comment{
\begin{equation}
    \includegraphics[width=5cm]{f0ex80000.pdf}.
\end{equation}}
However, $\mathcal{LA}_{\Tilde{P}_1}$ would have genus zero as can be seen either from $\epsilon=4\log2\approx1.204$ and the figure, or from checking lopsidedness. Therefore, we now consider larger $n$ for better approximations of the amoebae. Since the basic approach is similar to what we have discussed in \S\ref{n1}, let us again illustrate this with $F_0$ as an example.

\subsubsection{Example: $F_0$}\label{largernF0}
\begin{figure}[h]
    \centering
    \includegraphics[width=5cm]{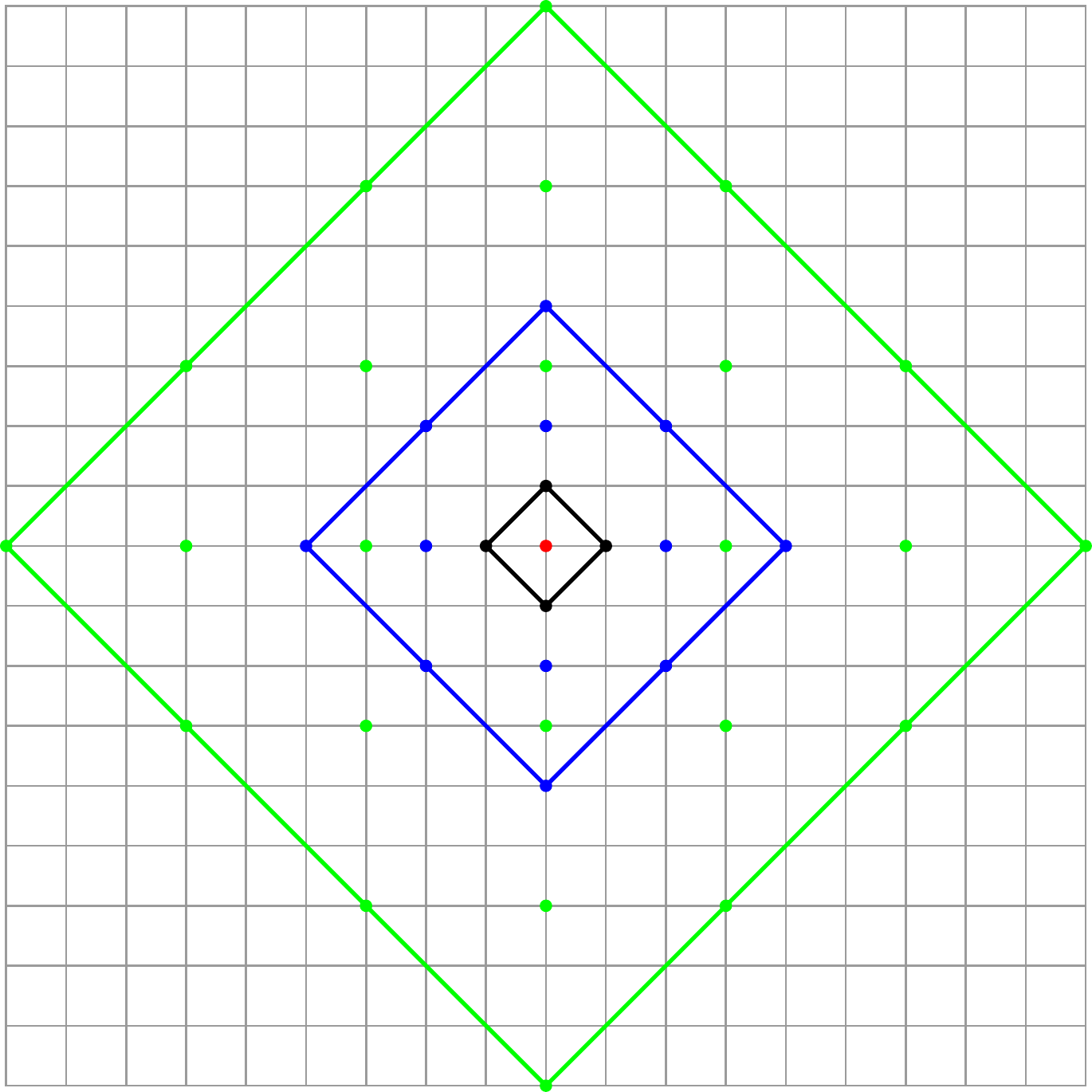}
    \caption{The polygons $(n^2\Delta_{F_0})\cap(n\mathbb{Z})^2$ in the lattice $\mathbb{Z}^2$. The black one is $F_0$, i.e., $n=1$. The blue and green ones correspond to $n=2,3$ respectively. In particular, the centre point in red, which is also the polygon with $n=0$, belongs to $(n^2\Delta_{F_0})\cap(n\mathbb{Z})^2$ for all $n$.}\label{f0dilation}
\end{figure}

In general, for any toric polygon $\Delta$, the cyclic resultant $\Tilde{P}_n$ also corresponds to a Newton polytope which is $(n^2\Delta)\cap(n\mathbb{Z})^2$. For instance, we plot the first several Newton polytopes for $\Tilde{P}_n$ of $F_0$ in Figure \ref{f0dilation}.
In general, $(n^2\Delta)\cap(n\mathbb{Z})^2$ would have $(2n^2+2n+1)$ lattice points (or equivalently, monomials in $\Tilde{P}_n$), among which $4n$ are boundary points.

Using $\Tilde{P}_n$, which are listed in Appendix \ref{cycresF0} by explicit computation of the cyclic resultant, we can get the conditions for the genus (though they are much more complicated than the case of $n=1$). As $n\rightarrow\infty$, we would recover the amoeba $\mathcal{A}_P$. To generate a data-set for $\mathcal{A}_P$, the following trick is used. We check whether the centre point (which is derived in Appendix \ref{lopg}) lives in the amoeba, i.e., whether there exists a solution to $P(z,w)=0$ such that
\begin{equation}
    (\text{Log}|z|,\text{Log}|w|)=\left(\frac{1}{2}\text{Log}\left|\frac{c_3}{c_1}\right|,\frac{1}{2}\text{Log}\left|\frac{c_4}{c_2}\right|\right).
\end{equation}
Again, we use our MLP classifier with input $(c_1,\dots,c_5)$ and output $g$. We also take $(|c_1|,\dots,|c_5|)$ as input to perform the same test. For each test, we generate a balanced dataset with $\sim2000$ samples. The results are listed in Table \ref{F0largerntable}.
\begin{table}[h]
\centering
\begin{tabular}{c||c|c|c|c|c}
Input& $n=1$ & $n=2$ & $n=3$ & $n=4$ & $n\rightarrow\infty$ \\ \hline
$c_i$ & $0.957(\pm0.005)$ & $0.956(\pm0.007)$ & $0.960(\pm0.008)$ & $0.933(\pm0.009)$ & $0.915(\pm0.016)$ \\ \hline
$|c_i|$ & $0.987(\pm0.007)$ & $0.912(\pm0.011)$ & $0.905(\pm0.013)$ & $0.903(\pm0.009)$ & $0.753(\pm0.056)$
\end{tabular}
\caption{The accuracies for MLP using 5-fold cross validation with $95\%$ confidence interval. We use $n\rightarrow\infty$ to denote the case for $\mathcal{A}_P$. The result of $n=1$ from the previous section is also listed here for reference.}\label{F0largerntable}
\end{table}


As we can see, unlike $n=1$, for other cases, the performance using $c_i$ as input is always better than the one using $|c_i|$. In fact, this is very reasonable since only $n=1$ has sets $P\{\bm{z}\}$ with elements $|c_i|$ (multiplied by certain $z^iw^j$) while for larger $n$ the elements can be complicated expressions of $c_i$'s. For instance, as shown in Appendix \ref{cycresF0}, there exists an element $|(4c_1^3c_3-4c_1^2c_2c_4-2c_1^2c_5^2)w^4z^6|$ in $\Tilde{P}_2\{\bm{z}\}$ for $F_0$. Therefore, taking absolute values is not equivalent to restricting to positive coefficients and could mislead the machine.

\begin{figure}[H]
	\centering
	\begin{subfigure}{7cm}
		\centering
		\includegraphics[width=7cm]{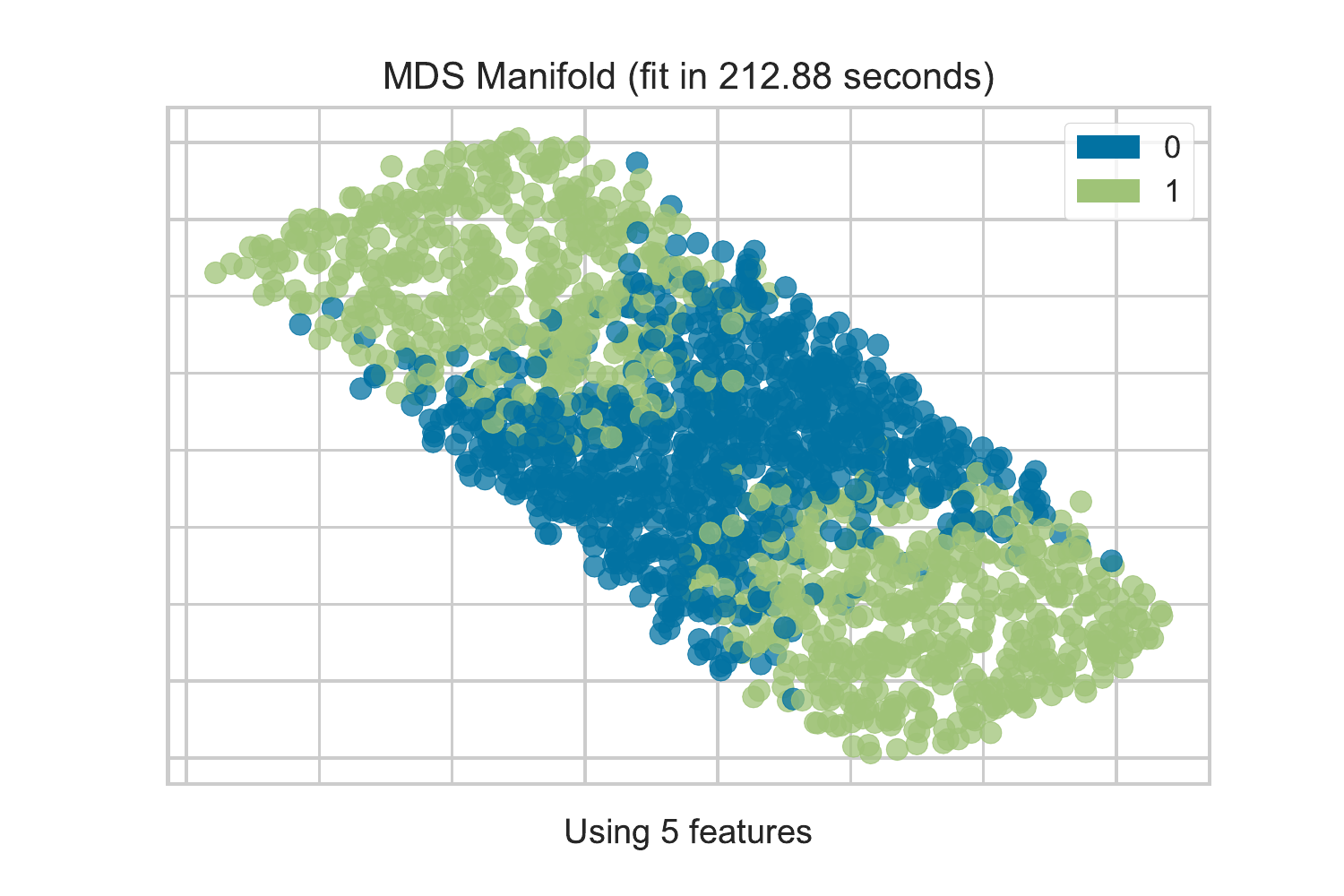}
		\caption{}
	\end{subfigure}
    \begin{subfigure}{7cm}
    	\centering
    	\includegraphics[width=7cm]{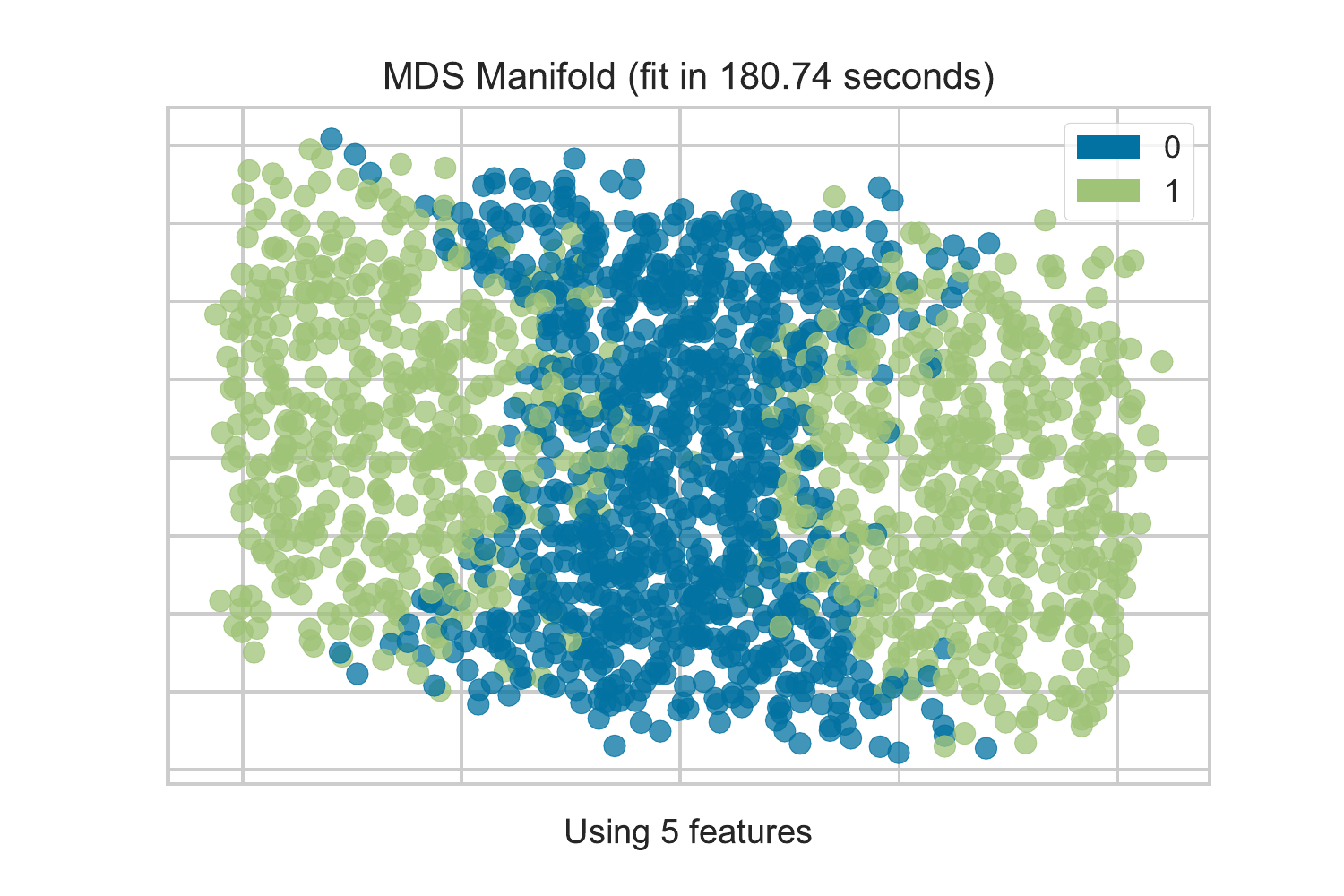}
    	\caption{}
    \end{subfigure}
    \begin{subfigure}{7cm}
    	\centering
    	\includegraphics[width=7cm]{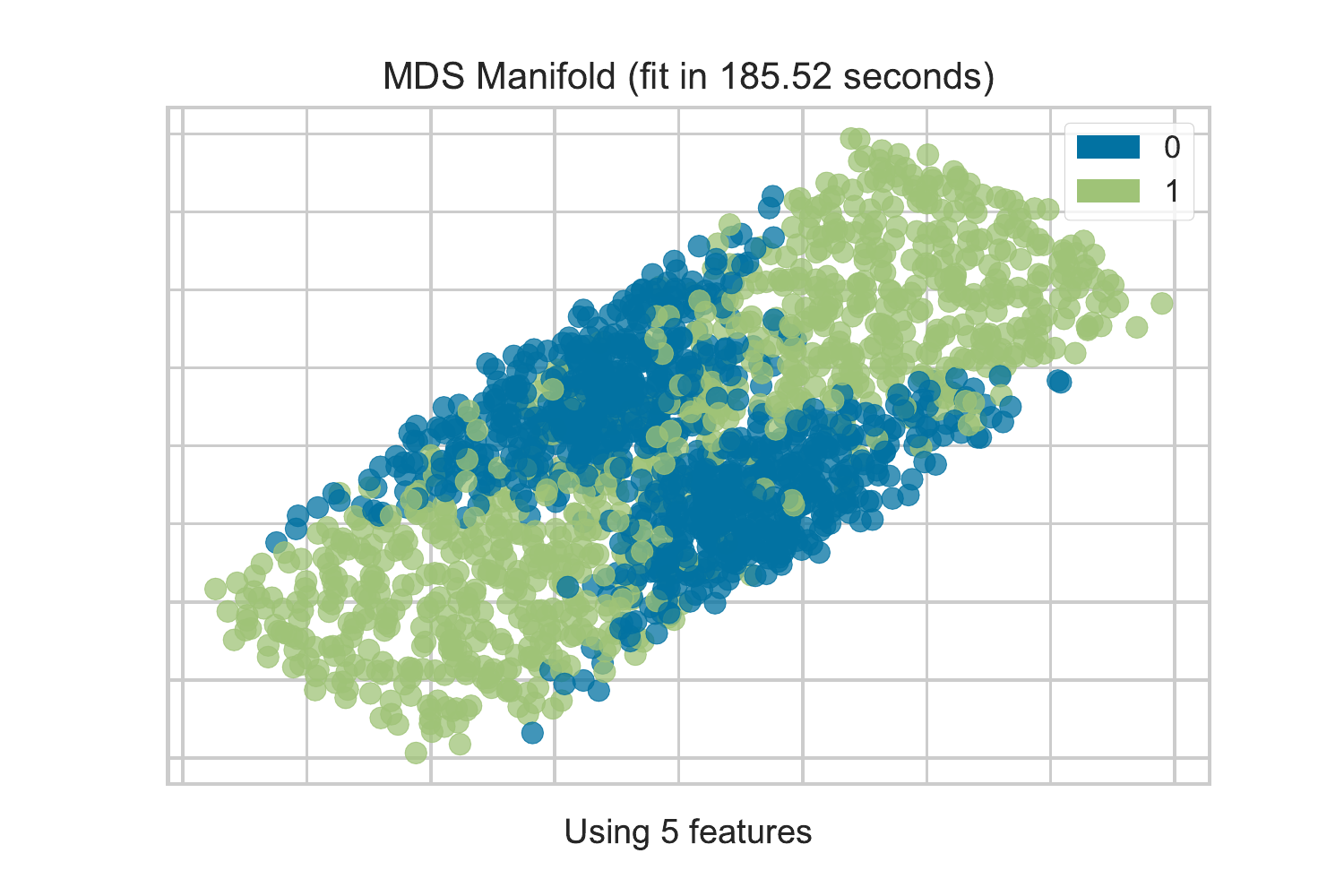}
    	\caption{}
    \end{subfigure}
    \begin{subfigure}{7cm}
    	\centering
    	\includegraphics[width=7cm]{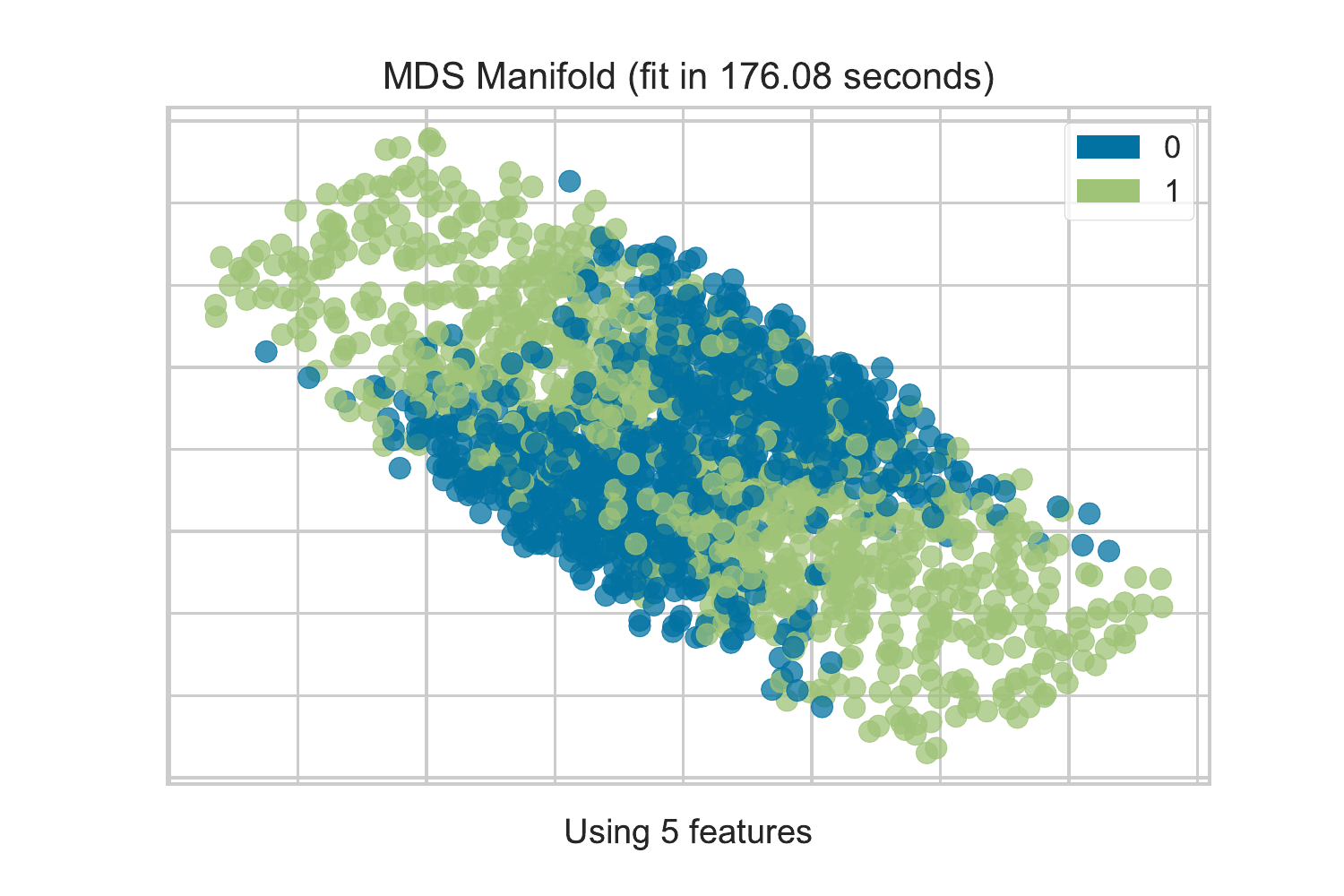}
    	\caption{}
    \end{subfigure}
    \caption{The MDS manifold projections for the datasets with (a) $n=2$, (b) $n=3$, (c) $n=4$ and (d) $n\rightarrow\infty$.}\label{F0mdslargern}
\end{figure}

We may still use the MDS manifold embedding to project the inputs $c_i$ to 2d. This is visualized in Figure \ref{F0mdslargern}.
Although it would be hard to tell what the two components are as the conditions from lopsidedness are much more complicated, it is interesting to see that these distributions of blue and green points have similar shapes/regions compared to $n=1$ \footnote{One may also project the inputs $|c_i|$ to 2d as before. It turns out that the shapes/regions of distributions are also similar to case with $n=1$ but with more mixings at the boundary of two regions. However, we are not really interested in this as taking absolute values is not helpful at all here.}.

\begin{figure}[h]
	\centering
	\begin{subfigure}{5cm}
		\centering
		\includegraphics[width=4cm]{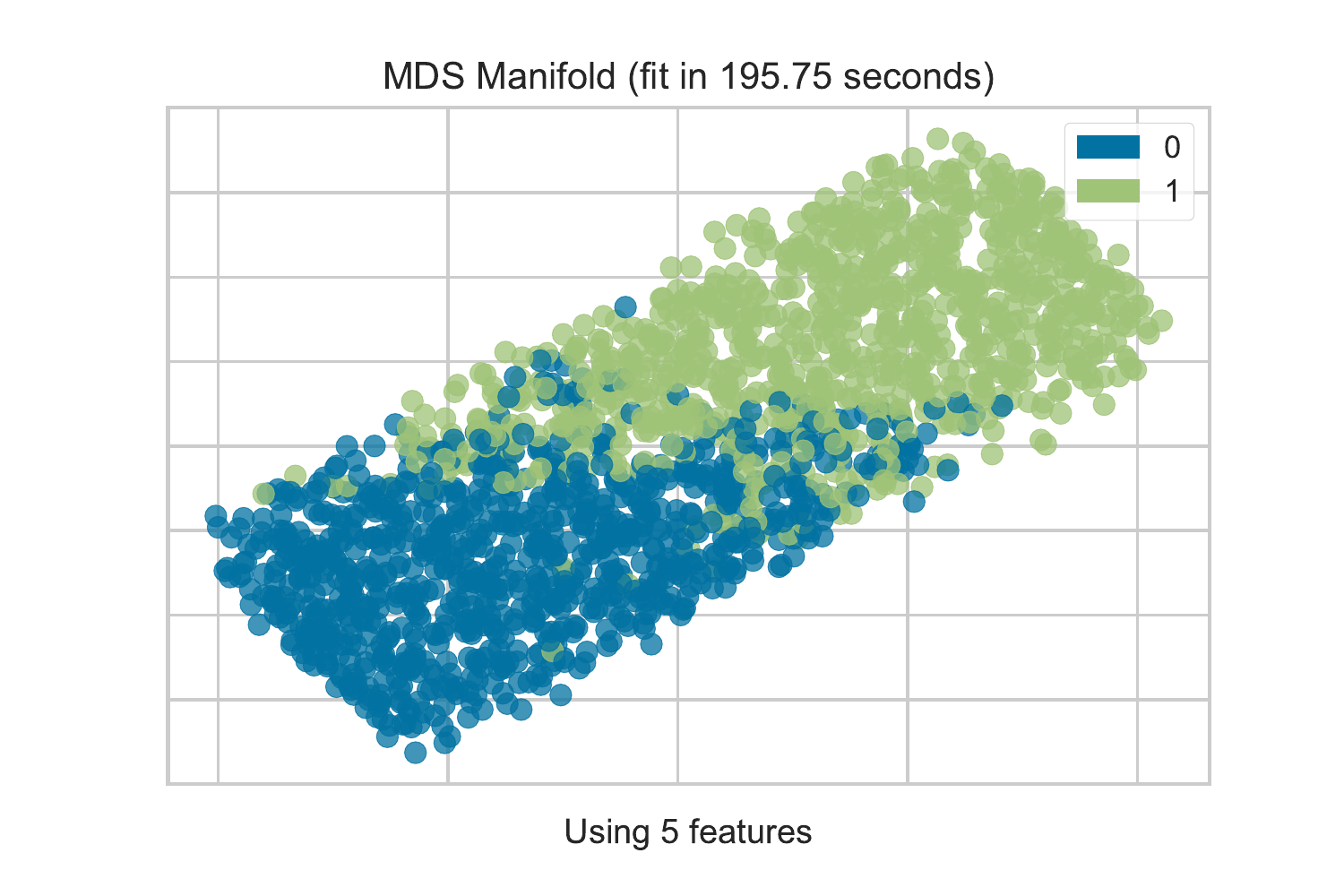}
		\caption{}
	\end{subfigure}
    \begin{subfigure}{5cm}
    	\centering
    	\includegraphics[width=4cm]{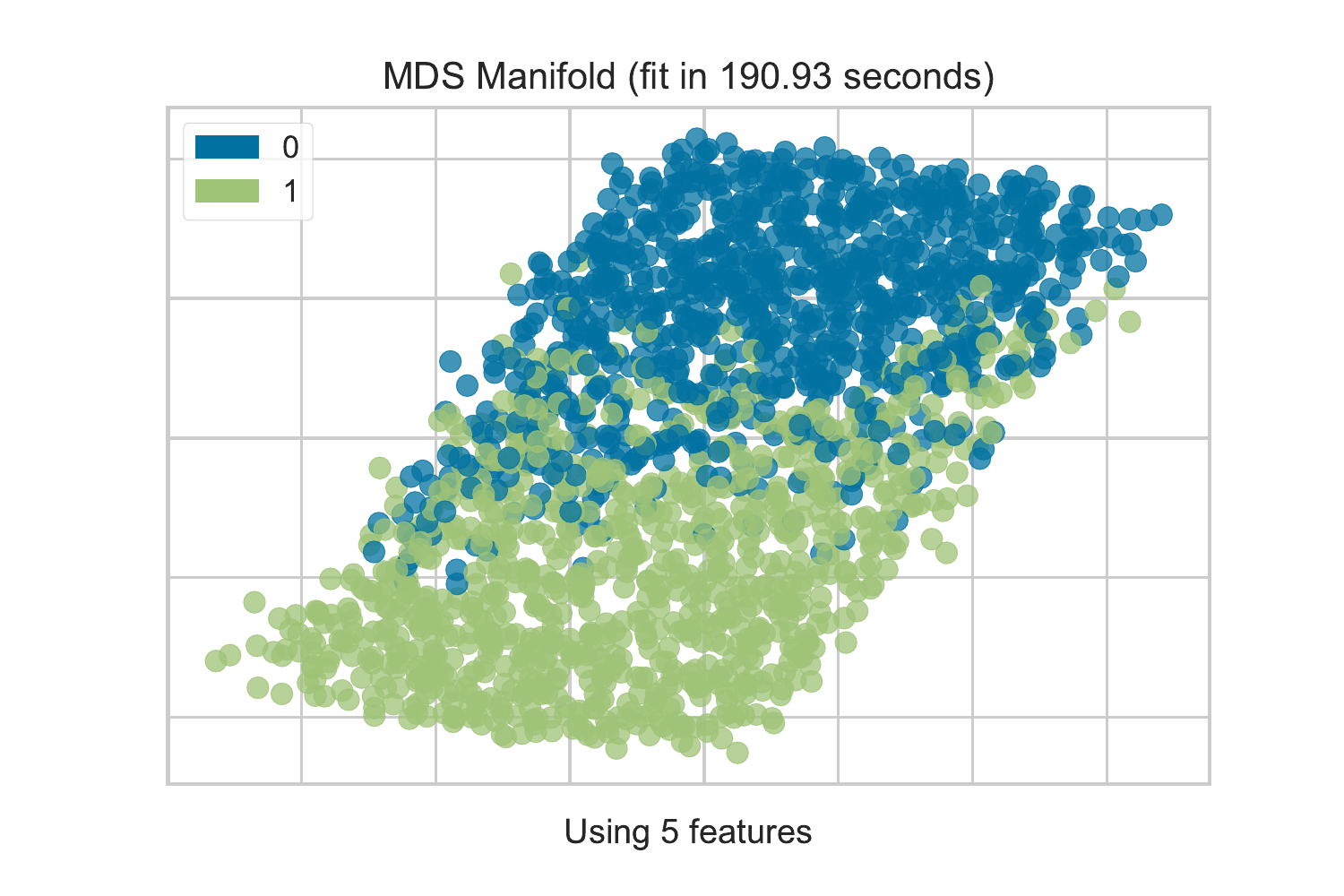}
    	\caption{}
    \end{subfigure}
    \\
    \begin{subfigure}{5cm}
    	\centering
    	\includegraphics[width=4cm]{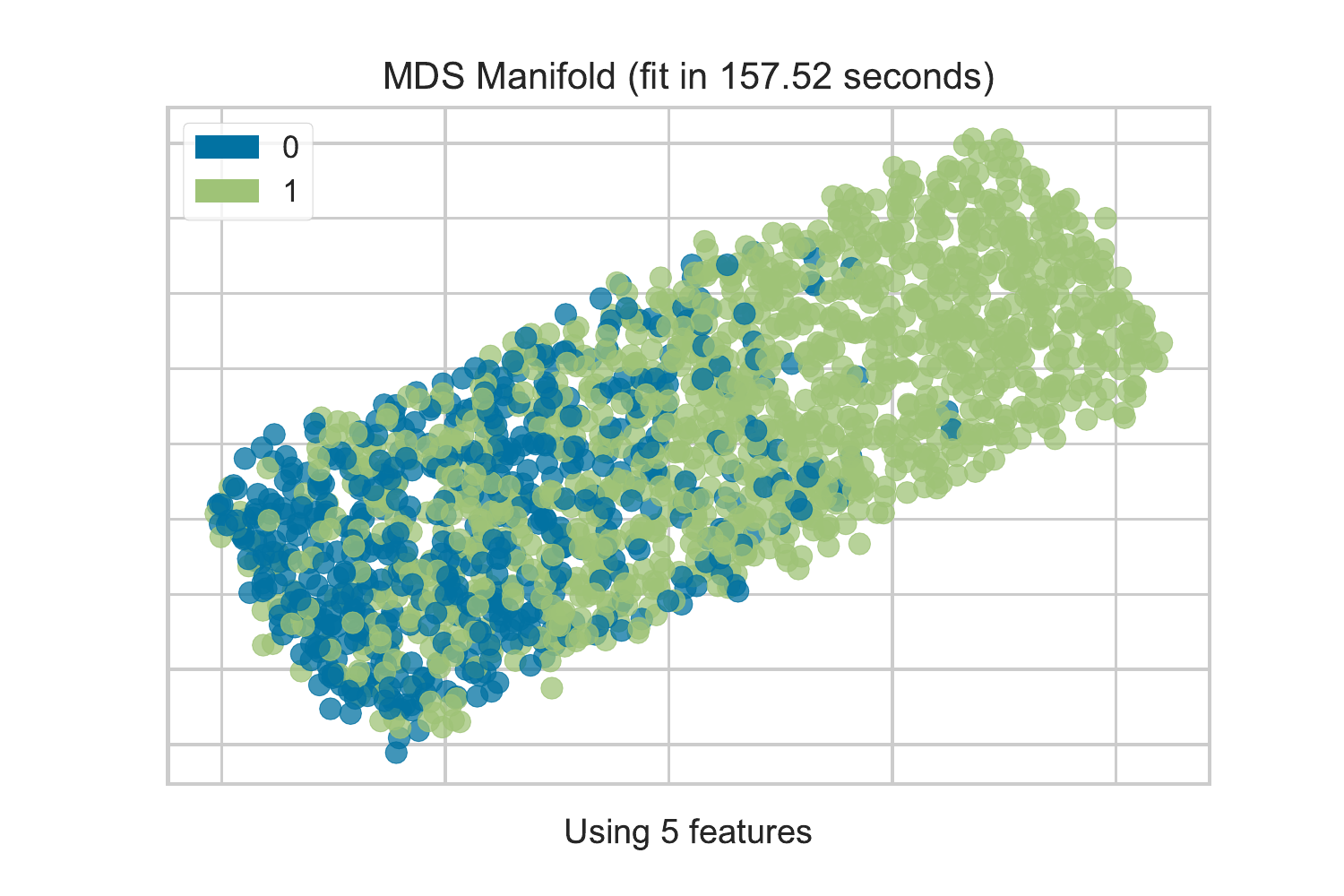}
    	\caption{}
    \end{subfigure}
    \begin{subfigure}{5cm}
    	\centering
    	\includegraphics[width=4cm]{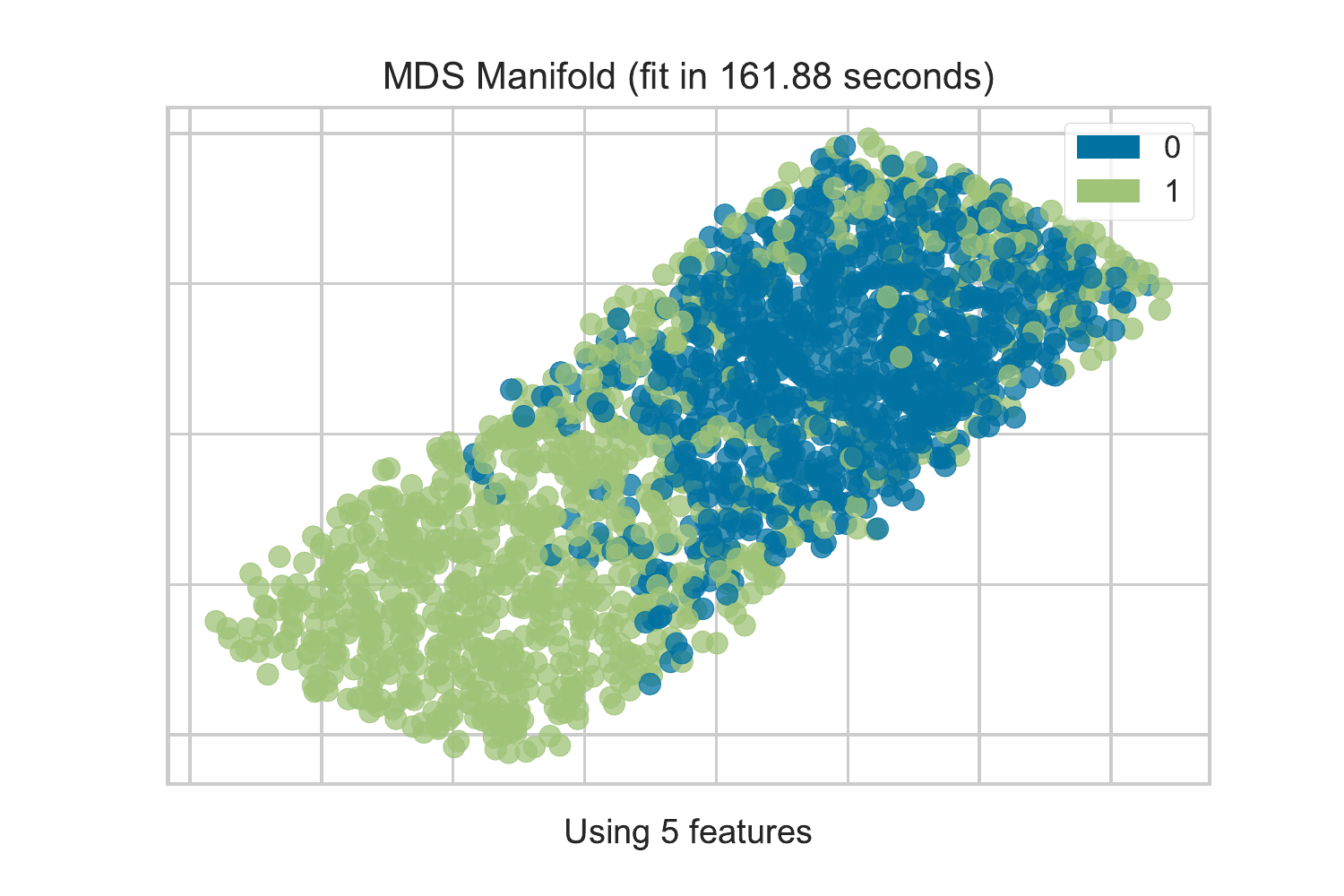}
    	\caption{}
    \end{subfigure}
    \caption{The MDS manifold projections for the datasets with (a) $n=2$, (b) $n=3$, (c) $n=4$ and (d) $n\rightarrow\infty$. Here the coefficients $c_i$ are restricted to be positive.}\label{F0mdslargernpositive}
\end{figure}

Now that taking absolute values is completely different from taking non-negative coefficients, we generate a dataset for non-negative coefficients with $\sim2000$ samples for each $n$ and perform the same tests as above. The result is in Table \ref{F0largerntablepositive}.
\begin{table}[H]
\centering
\begin{tabular}{c||c|c|c|c|c}
Input& $n=1$ & $n=2$ & $n=3$ & $n=4$ & $n\rightarrow\infty$ \\ \hline
$c_i\in\mathbb{R}^+$ & $0.987(\pm0.007)$ & $0.983(\pm0.008)$ & $0.985(\pm0.002)$ & $0.976(\pm0.006)$ & $0.930(\pm0.008)$
\end{tabular}
\caption{The accuracies for MLP using 5-fold cross validation with $95\%$ confidence interval.}\label{F0largerntablepositive}
\end{table}
Compared to the results for real coefficients in Table \ref{F0largerntable}, we see that the performance all get improved. In particular, the cases for finite $n$ listed here all give almost perfect results while the accuracy for $n\rightarrow\infty$ is also slightly improved. As plotted in Figure \ref{F0mdslargernpositive}, the MDS manifold projections for positive coefficients also have similar distributions compared to the one for positive coefficients with $n=1$.

\subsubsection{Reproduction of the Genus}\label{reproduce2}
We may also use the weight matrices and bias vectors of the MLP to get nice approximated expressions of $g$ for different $n$'s as what we previously did for $n=1$. Here, for $F_0$, we present the weight matrices and bias vectors for the amoeba $\mathcal{A}_P$ (viz, $n\rightarrow\infty$) as an example.

\begin{landscape}
It turns out an MLP of one hidden layer with 20 perceptrons could already give high accuracy as in Table \ref{F0largerntable}. We still use ReLU as activation functions. The weight matrices are
\begin{equation}
    W_1=\left({\scriptsize\begin{array}{cccccccccccccccccccc}
    0.597 & 0.911 & -0.981 & 0.280 & -0.458 & -0.359 & -0.143 & -1.183 & -0.559 & -0.864 & -0.378 & -1.076 & -0.046 & 0.624 & -0.804 & -0.656 & -0.067 & 0.274 & -0.516 & -1.021 \\
    -1.031 & 0.910 & 0.444 & 0.887 & 0.898 & 1.117 & -0.530 & 0.182 & -0.355 & 1.347 & -0.435 & 0.769 & 1.379 & 1.558 & -0.384 & -0.779 & 1.054 & 0.747 & 0.146 & -0.346\\
    0.717 & 1.528 & -2.440 & -0.208 & -1.299 & 0.285 & 0.417 & -1.318 & -0.644 & -0.794 & -1.005 & 1.408 & 0.228 & 0.870 & 1.738 & -0.166 & 0.585 & -0.496 & -1.869 & 1.136\\
    -0.596 & 0.757 & 0.284 & 0.814 & -0.940 & -0.802 & -1.363 & -0.550 & 0.375 & 1.116 & 0.730 & 1.308 & -1.127 & 0.834 & -0.373 & -2.297 & 1.203 & 0.780 & 0.332 & 0.630\\
    -0.103 & -0.556 & -0.035 & -0.811 & -0.503 & 0.475 & -0.475 & 1.062 & 0.360 & 0.156 & -0.550 & 0.552 & 0.165 & 0.643 & 0.207 & -0.096 & 1.380 & 0.455 & -1.287 & 0.817
    \end{array}}\right)^\text{T}
\end{equation}
and
\begin{equation}
    W_2=\left({\footnotesize\begin{array}{cccccccccccccccccccc}
    0.144 & -0.958 & -1.435 & 1.160 & 0.629 & -1.311 & 0.183 & 1.193 & 0.701 & -1.597 & 0.331 & 1.064 & 1.394 & -1.296 & 1.108 & -1.078 & 1.262 & 0.891 & 0.511 & -1.776
    \end{array}}\right).
\end{equation}
The bias vectors are
\begin{equation}
    \begin{split}
        \bm{b}_1=&(-0.044, 0.320, 0.234, 0.421, 0.221, 0.057, 0.088, -0.317, 0.387, -0.054,\\
        &-0.154, -0.296, 0.426, 0.051, 0.368, 0.162, 0.006, -0.004, -0.394, -0.072)^\text{T}
    \end{split}
\end{equation}
and
\begin{equation}
    \bm{b}_2=(-0.943).
\end{equation}
This gives a non-negative number $p$
with $g=0$ only when $p=0$:
\begin{equation}
    p=W_2\cdot\text{ReLU}(W_1\cdot(c_1,\dots,c_5)^\text{T}+\bm{b}_1)+\bm{b}_2 \ \Rightarrow \ g=\theta(p) \ .
    \label{F0approx1}
\end{equation}
Recall that computing the genus for $F_0$ is, in fact, finding a solution to
\begin{equation}
    c_1z+c_2w+c_3z^{-1}+c_4w^{-1}+c_5=0
\end{equation}
with conditions
\begin{equation}
    |z|=\left|\frac{c_3}{c_1}\right|^{1/2},~|w|=\left|\frac{c_4}{c_2}\right|^{1/2}.
\end{equation}

If there is no such solution, then $g=1$. Otherwise, $g=0$. Now, we see that this is translated to \eqref{F0approx1} \footnote{More strictly, this approximate expression is valid in the ranges of the coefficients in our data-set, which are $c_{1,2,3,4}\in[-5,5]$ and $c_5\in[-20,20]$. Of course, if we have larger ranges for the coefficients, it is still possible to obtain certain approximations following the same recipe.}.
\end{landscape}

\subsection{The Membership Problem}\label{membership}
So far, we have focused on detecting edges of holes and determining the genus of the amoeba.
We can also use machine learning to study another important problem, that of membership.
We introduced this in \S\ref{rudiments}:
given the Newton polynomial (always bi-variate in this paper), how does one decide whether a point in the plane belongs to the amoeba?
Thus, our input is of form $\{c_1,\dots,c_m,x_1,x_2\}$, where $c_i$ are the coefficients of the Newton polynomial as before and $(x_1,x_2)$ is a point on the Log plane. If $(x_1,x_2)$ belongs to the amoeba, then the output is 1. Otherwise, the output is 0.
This is a binary classification problem.

Again, we illustrate with our running example of $F_0$. 
Here, the input vectors are of the form $\{c_1,\dots,c_5,x_1,x_2\}$. 
For a balanced dataset with $\sim5000$ samples, we find that our MLP with one hidden layer of 8 perceptrons already gives $0.944(\pm0.003)$ accuracy for 5-fold cross validation with $95\%$ confidence interval\footnote{This is probably be the simplest MLP to give the best performance. We also checked MLPs with fewer perceptrons in the hidden layer, but all of them were inferior. On the other hand, if we include more neurons, the performance hardly improved (sometimes even worsened).}. Therefore, the weight matrices and bias vectors could provide good approximated expressions to determine the regions for amoebae.

The weight matrices are
\begin{equation}
    W_1=\begin{pmatrix}
    0.246 & -0.401 & 0.287 & -0.105 & 0.006 & 0.553 & 0.700 & -0.006\\
    -0.123 & -0.039 & 0.088 & -0.139 & 0.002 & 0.147 & -0.031 & 0.008\\
    0.198 & -0.104 & 0.040 & 0.058 & 0.000 & 0.066 & 0.208 & 0.009\\
    0.134 & -0.216 & 0.157 & 0.158 & 0.000 & 0.024 & 0.292 & 0.018\\
    -0.067 & 0.106 & -0.058 & -0.014 & -0.002 & -0.117 & -0.199 & -0.012\\
    -0.011 & -2.796 & -3.378 & -2.026 & -3.72 & -1.365 & 3.522 & -5.367\\
    -2.396 & 0.254 & -2.870 & 3.128 & 3.763 & 3.297 & 1.409 & -5.014
    \end{pmatrix}^\text{T}
\end{equation}
and
\begin{equation}
    W_2=\begin{pmatrix}
   1.273 & 3.525 & 2.168 & 2.275 & -4.418 & 1.993 & -1.313 & -6.769
    \end{pmatrix}.
\end{equation}
The bias vectors are
\begin{equation}
    \bm{b}_1=(0.337, -3.304, 3.064, -0.646, -2.778, 0.453, 2.778, -1.000)^\text{T} \ ;
    b_2 = -1.526 \ .
\end{equation}
This gives a non-negative number
\begin{equation}
    p=W_2\cdot\text{ReLU}(W_1\cdot(c_1,\dots,c_5,x_1,x_2)^\text{T}+\bm{b}_1)+{b}_2.
    \label{membershipapprox1}
\end{equation}
Then $(x_1,x_2)$ is in the complement of amoeba only when $p\leq0$. In other words,
\begin{equation}
    (x_1,x_2)\in\mathcal{A}_P~~\text{if}~~\theta(p)=1.
    \label{membershipapprox2}
\end{equation}

As we can see, the membership problem is now transformed to \eqref{membershipapprox1} and \eqref{membershipapprox2} with lower calculation cost\footnote{To be strict, such approximated expression should be valid in the ranges of the coefficients in our dataset, which are $c_i,x_i\in[-5,5]$.}. Of course, it is natural to expect that we can get better approximations with more training data and a more complicated neural network. Such method can also be applied to $c_i$ and $(x_1,x_2)$ with larger ranges and amoebae for other Newton polytopes.

\subsubsection{Crawling of Amoebae}\label{crawling}
\begin{figure}[H]
    \centering
    \includegraphics[width=8cm]{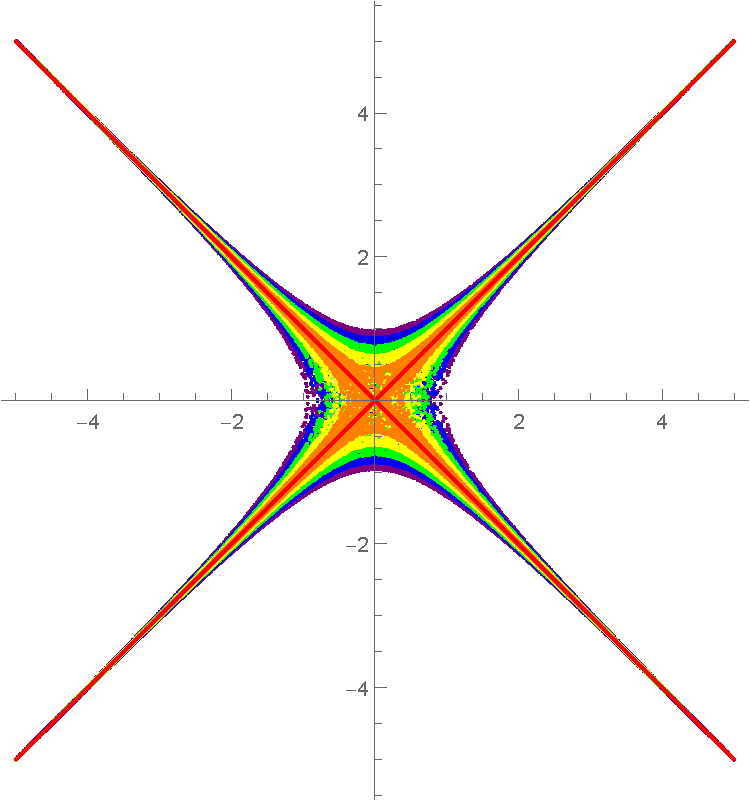}
    \caption{How the amoeba ``crawls'' for $z+w+z^{-1}+w^{-1}+c_5=0$, where $c_5=0,0.2,0.4,0.6,0.8,1$ are in red, orange, yellow, green, blue and purple respectively.}\label{crawlhue}
\end{figure}

Let us end this section by a quick comment on the boundaries of amoebae. For a given Newton polynomial, if we vary the coefficients continuously, then the shape/boundary of the amoeba should also change continuously. We shall refer to this as the ``crawling'' of amoebae. Again, let us use $F_0$ as our example. In Figure \ref{crawlhue}, we show how the amoeba crawls for $z+w+z^{-1}+w^{-1}+c_5=0$ with $c_5\in\{0,0.2,0.4,0.6,0.8,1\}$.

\begin{figure}[H]
	\centering
	\begin{subfigure}{4cm}
		\centering
		\includegraphics[width=4cm]{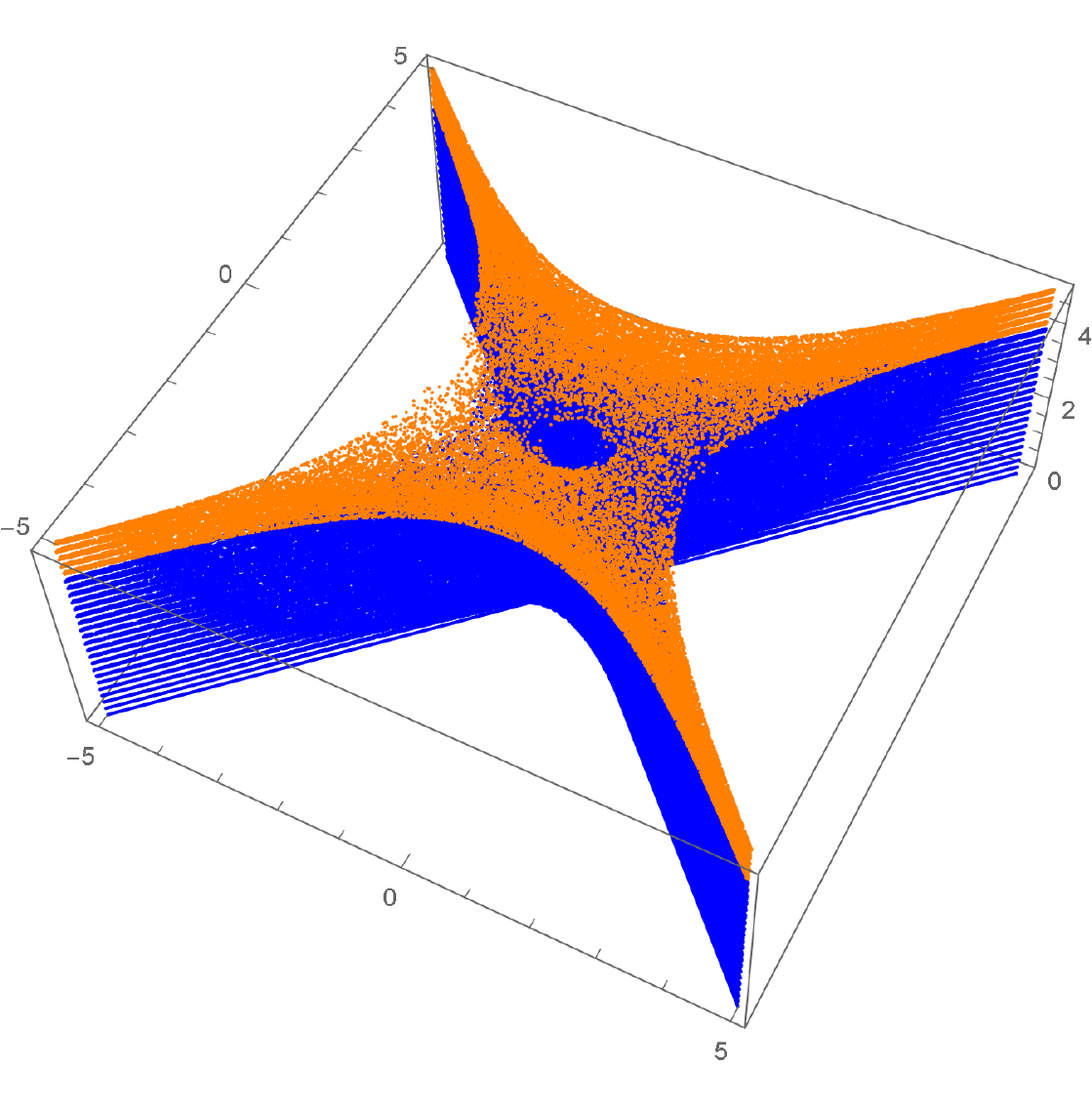}
		\caption{}
	\end{subfigure}
    \begin{subfigure}{4cm}
    	\centering
    	\includegraphics[width=4cm]{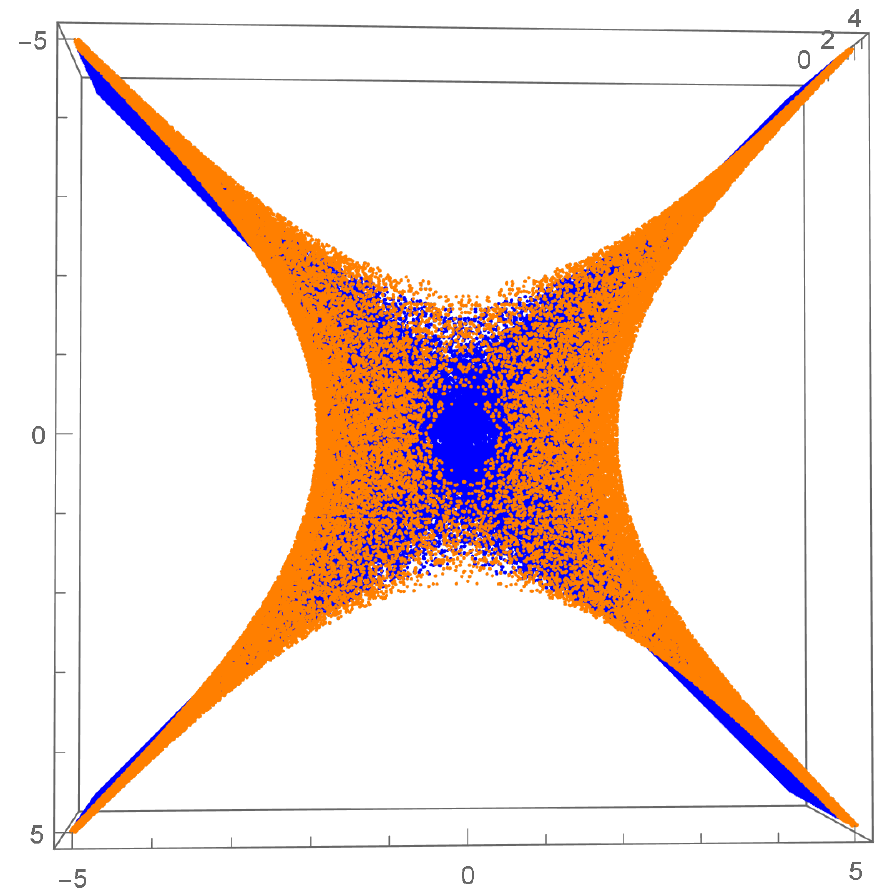}
    	\caption{}
    \end{subfigure}
    \begin{subfigure}{4cm}
    	\centering
    	\includegraphics[width=4cm]{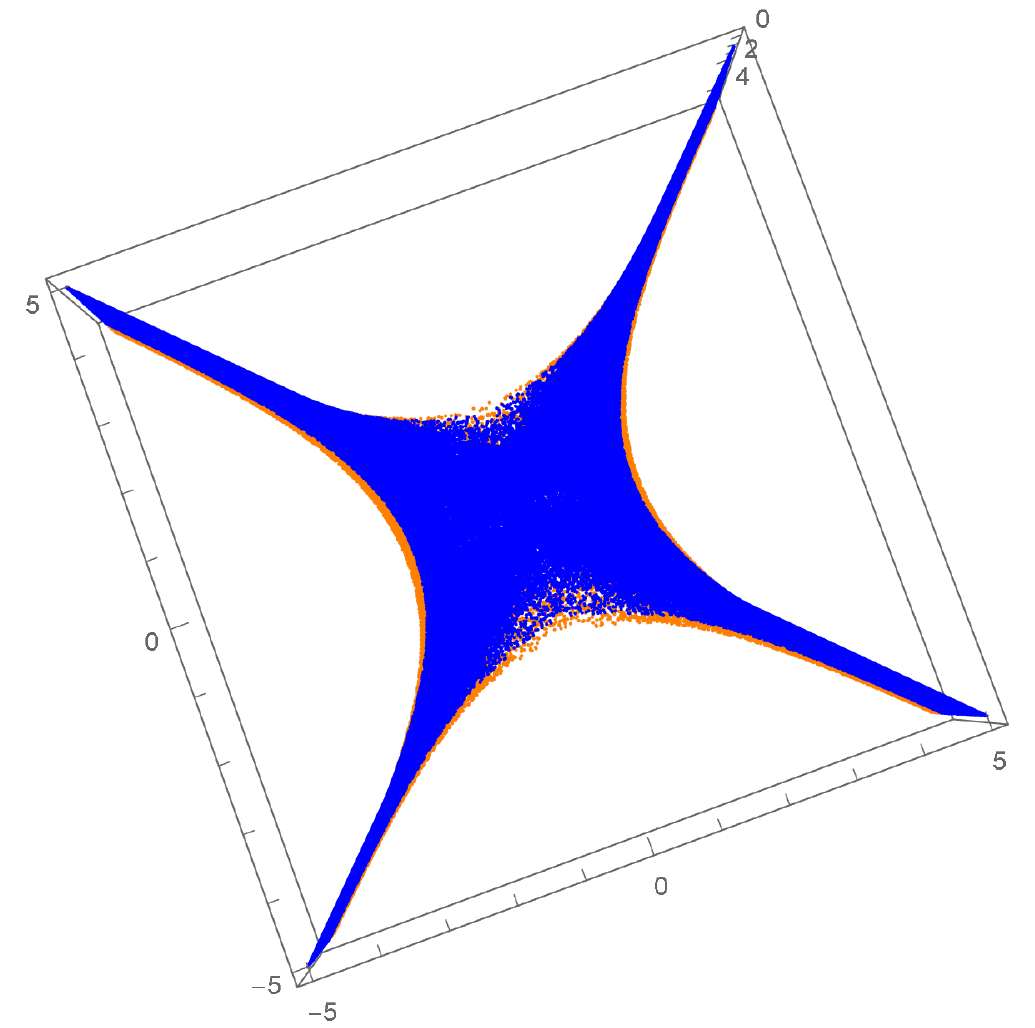}
    	\caption{}
    \end{subfigure}
    \begin{subfigure}{6cm}
    	\centering
    	\includegraphics[width=6cm]{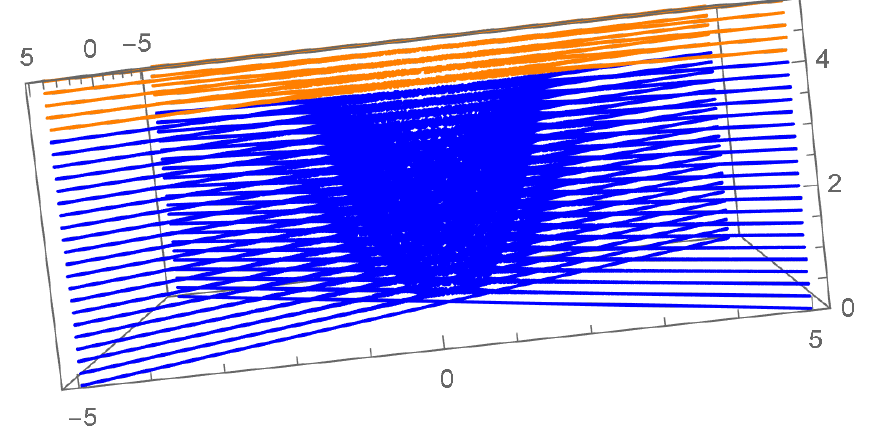}
    	\caption{}
    \end{subfigure}
    \caption{(a) The 3d plot of amoeba varying $c_5$. (b) The same plot viewed from the top. (c) The same plot viewed from the bottom. (d) The same plot viewed from the side. The two colours are chosen based on whether the corresponding lopsided amoeba $\mathcal{LA}_P$ has genus 0 or 1. Hence, some blue slices near the orange part should also have genus 1.}\label{F0crawlc5}
\end{figure}

We may also plot this crawling in 3d with $c_5$ as the third axis. We vary $c_5$ from 0 to 5 with step 0.2 in Figure \ref{F0crawlc5}.

\begin{figure}[h!]
	\centering
	\begin{subfigure}{4cm}
		\centering
		\includegraphics[width=4cm]{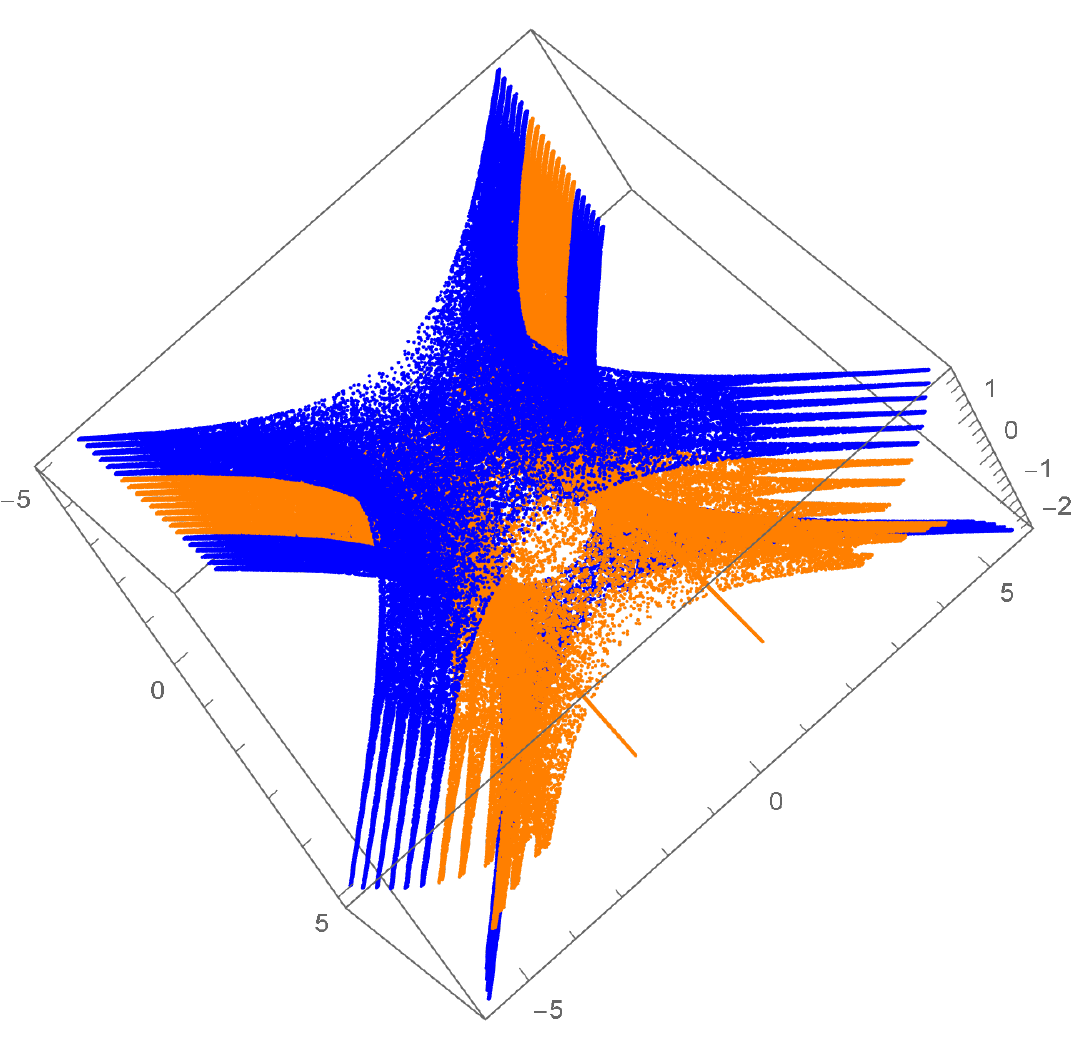}
		\caption{}
	\end{subfigure}
    \begin{subfigure}{4cm}
    	\centering
    	\includegraphics[width=4cm]{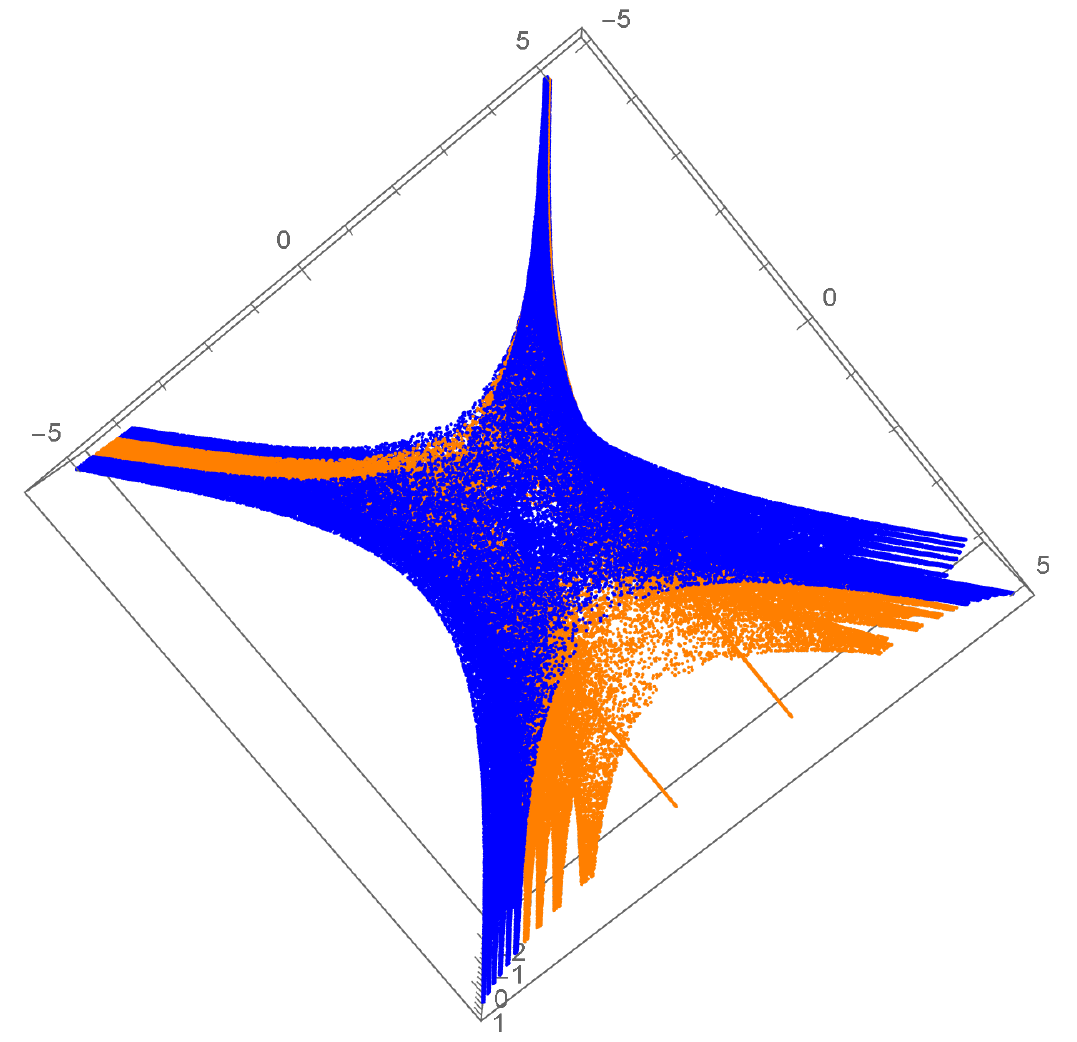}
    	\caption{}
    \end{subfigure}
    \begin{subfigure}{4cm}
    	\centering
    	\includegraphics[width=4cm]{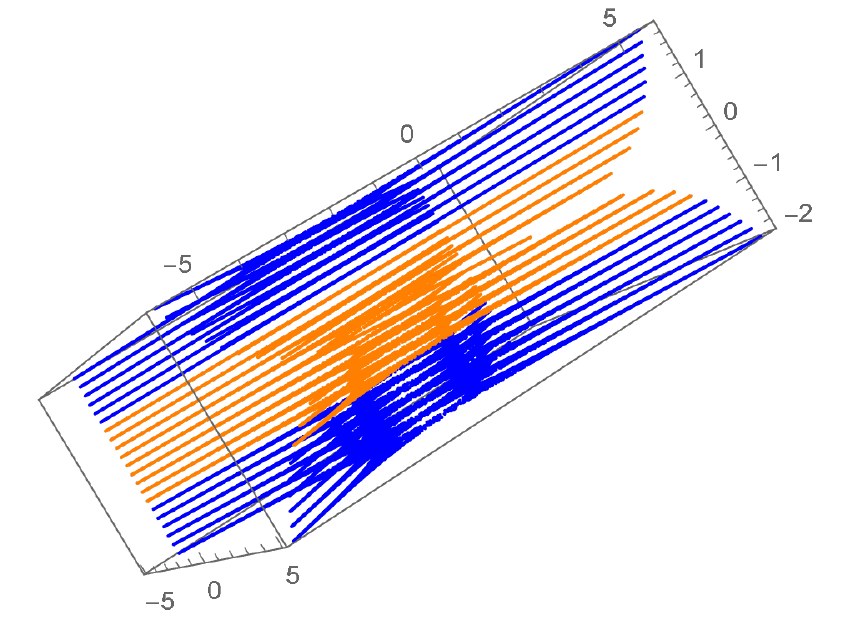}
    	\caption{}
    \end{subfigure}
    \caption{(a) The 3d plot of amoeba varying $c_1$. (b) The same plot viewed from the top. Viewing from the bottom would be the same due to the symmetry of the polynomial. (c) The same plot viewed from the side. The two colours are chosen based on whether the corresponding lopsided amoeba $\mathcal{LA}_P$ has genus 0 or 1. Hence, some blue slices near the orange part should also have genus 1. Notice that when $c_1=0$ there is no genus due to degeneracy.}\label{F0crawlc1}
\end{figure}

Likewise, in Figure \ref{F0crawlc1}, we vary $c_1$ while keeping the other coefficients fixed. The Newton polynomial now reads $c_1z+w+z^{-1}+w^{-1}+4=0$, where $c_1$ is from $-2$ to 2 with step 0.2.

Here, we only vary one coefficient and restrict the coefficients to be real so as to visualize a specific crawling of the amoebae. In general, one can vary multiple coefficients and consider the coefficients to be any complex numbers. If we know how an amoeba crawls when varying any combination of the coefficients (though not easy), we can determine the boundaries of all amoebae for a given Newton polytope/polynomial once we know the boundary information of one amoeba.

In fact, from the above discussion, we may further conjecture that when varying the coefficients, the (change of) boundary of an amoeba (i.e., the crawling of it) forms some (almost) manifold. In other words, it is smooth except when there is degeneracy. By degeneracy, we mean that some of the coefficients vanish in $P(z,w)$. For instance, the amoeba degenerates to $\text{Log}|w|=\pm\text{Log}|z|$ when $c_5=0$ in Figure \ref{F0crawlc5} while the amoeba degenerates to the amoeba of $\mathbb{C}^3/\mathbb{Z}_2~(1,1,1)$ when $c_1=0$ in Figure \ref{F0crawlc1}.

Overall, since the machine learning models here are aimed to show that a minimal number of data points and a minimal structure with very simple weights and bias could already perform well for the membership problem, and since the input is composed of varying coefficients and $(x_1,x_2)$ rather than $(x_1,x_2)$ only, the approximations of the boundaries is still not comparable to the approximations by lopsidedness or plotting directly using Monte Carlo method. However, it is reasonable to expect a more complicated neural network with more training data would refine such approximations, and we hope that the tests and observations here could shed light on the study of amoebae and their complementary regions.

\section{Image Processing Amoebae}\label{ImageCNN}
In addition to machine learning the abstract space of amoebae coefficients, it is interesting to test the success of ML on the amoebae images directly.
For this investigation $F_0$ is used again as our prototypical example, with the database consisting of Monte Carlo generated images of $n=1$ lopsided amoebae.

CNNs find there most common uses in image related tasks, due to their convolutional action over the local structure of the images.
Since this is the behaviour used in collecting Monte Carlo generated points into a full amoeba, this ML architecture is the most natural choice.
To fully investigate their success in identifying genus from these images, varying image resolutions are used for the CNN input.
Determining how the learning measures varied as the amoebae images varied in quality is the main focus of this investigation.

\subsection{Amoeba Image Dataset}
Here, our database consists of Monte Carlo generated {\bf images} of the $F_0$ amoebae for varying coefficients such that the full database consisted of 1000 genus 0 amoebae, and 1000 genus 1 amoebae.
The amoebae genus labelling used the $n=1$ lopsided amoeba approximation, as in equation \ref{gF0}.
Monte Carlo image generation is commonly used for more complicated amoebae where analytical expressions are hard to compute and then plot.
Therefore examining the success of a convolutional neural network (CNN) in identifying genus from these plots is a particularly relevant task.
We now undertake this task of the CNN ML of the labelled data-set
\[
\mbox{Image}(\mathcal{A}_{P(F_0)}) \longrightarrow \text{genus}\; \{0 \mbox{ or } 1 \}.
\]

The amoebae images in the dataset are of varying shape, to create consistency across the dataset (as needed for the CNN tensor inputs), the images are resized such that there is always an equal number of pixels in each dimension.
The images contained the real plane axes since the identification of a plot's origin is useful in determining genus for $F_0$ amoebae by eye, hence giving the CNN realistic information for its learning.
Images are reformatted to greyscale such that inclusion in the amoeba is the only relevant data at each resolution.
The resolutions learnt on varied logarithmically with base 2, choosing this base such that the computational learning algorithms worked most efficiently.

\begin{figure}[t]
	\centering
	\begin{subfigure}{4cm}
		\centering
		\includegraphics[width=4cm]{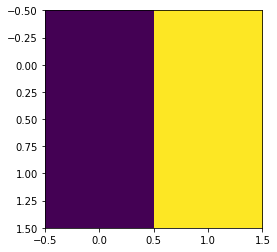}
		\caption{$2 \times 2$}\label{AmbImages2}
	\end{subfigure}
    \begin{subfigure}{4cm}
    	\centering
    	\includegraphics[width=4cm]{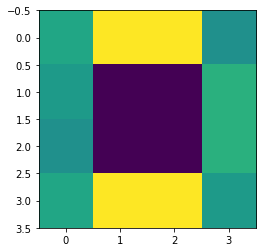}
    	\caption{$4 \times 4$}
    \end{subfigure}
    \begin{subfigure}{4cm}
    	\centering
    	\includegraphics[width=4cm]{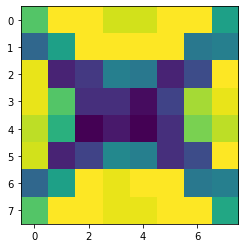}
    	\caption{$8 \times 8$}
    \end{subfigure}
    	\begin{subfigure}{4cm}
		\centering
		\includegraphics[width=4cm]{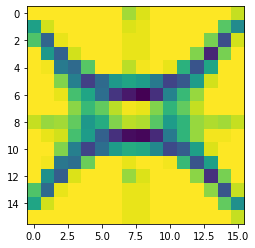}
		\caption{$16 \times 16$}
	\end{subfigure}
    \begin{subfigure}{4cm}
    	\centering
    	\includegraphics[width=4cm]{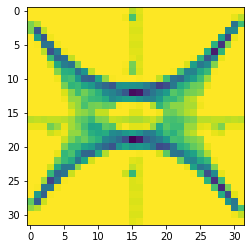}
    	\caption{$32 \times 32$}\label{AmbImages32}
    \end{subfigure}
    \begin{subfigure}{4cm}
    	\centering
    	\includegraphics[width=4cm]{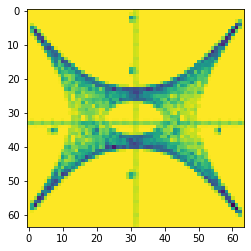}
    	\caption{$64 \times 64$}
    \end{subfigure}
    	\begin{subfigure}{4cm}
		\centering
		\includegraphics[width=4cm]{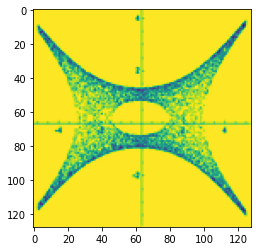}
		\caption{$128 \times 128$}
	\end{subfigure}
    \begin{subfigure}{4cm}
    	\centering
    	\includegraphics[width=4cm]{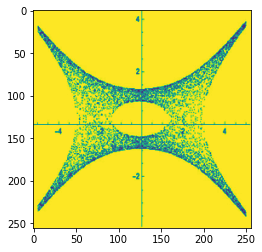}
    	\caption{$256 \times 256$}
    \end{subfigure}
    \begin{subfigure}{4cm}
    	\centering
    	\includegraphics[width=4cm]{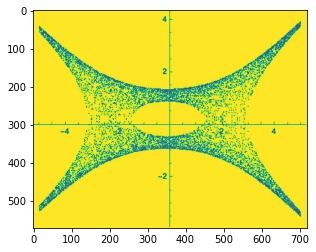}
    	\caption{$750 \times 572$ \\ (original)}
    \end{subfigure}
    \caption{An example $F_0$ amoeba image of genus 1 at varying resolutions. Each subcaption denotes the respective number of pixels in the resizing, the original image is reformatted to a square shape for consistency across all amoebae.}\label{AmbImages}
\end{figure}

An example of a genus 1 amoeba from the dataset is given in Figure \ref{AmbImages}, this amoebae had randomly generated $\mathbb{R}^+$ coefficients: 
\begin{equation}
    (c_1, \ldots, c_5) \sim
\{
0.6104, 
1.8940, 
0.4989, 
2.9777, 
6.9871 
\} \ .
\end{equation}
In the figure the image at the varying resolutions considered is shown, the reshaping to square from the original is clear, and the binary nature of the data is emphasised by the greyscaling (here plotted with a blue-yellow colour scheme).

\paragraph{CNN Architecture: }
The $\texttt{Tensorflow}$ library $\texttt{Keras}$ is used for implementing the CNN architecture \cite{tensorflow2015-whitepaper}. The CNNs trained have a consistent layer structure: 3 $\times$ 2d convolutional layers of size matching the image input size (each followed by a 2d max-pooling layer and a dropout layer (factor 0.01)); following these is a flattening layer, then a dense layer of size equal to the number of pixels in one dimension of the input image, and then a final output dense layer with 1 neuron. 
All layers used a LeakyReLU activation (with factor 0.01), except the final output Dense layer which used sigmoid activation to better map to the required binary classification.
The convolutional and max-pooling layers used the `same' padding regime, and the convolutional layers used a $3 \times 3$ kernel size.

The CNN architecture used for the investigation at each image resolution is trained over 20 epochs of the data on batches of 32, in a 5-fold cross-validation scheme. 
The Adam optimiser \cite{kingma2017adam} is used to minimise the binary cross-entropy loss function for predicting genus 0 or genus 1 respectively.
The learning is measured using the metrics: accuracy and Matthews' Correlation Coefficient (MCC) \footnote{The MCC is equal to the chi-square of the $2 \times 2$ confusion matrix in the classification.}, which are averaged over the 5 runs in the 5-fold cross-validation.

\subsection{Image Processing Results}
\begin{figure}[h!]
    \centering
    \includegraphics[width=14cm]{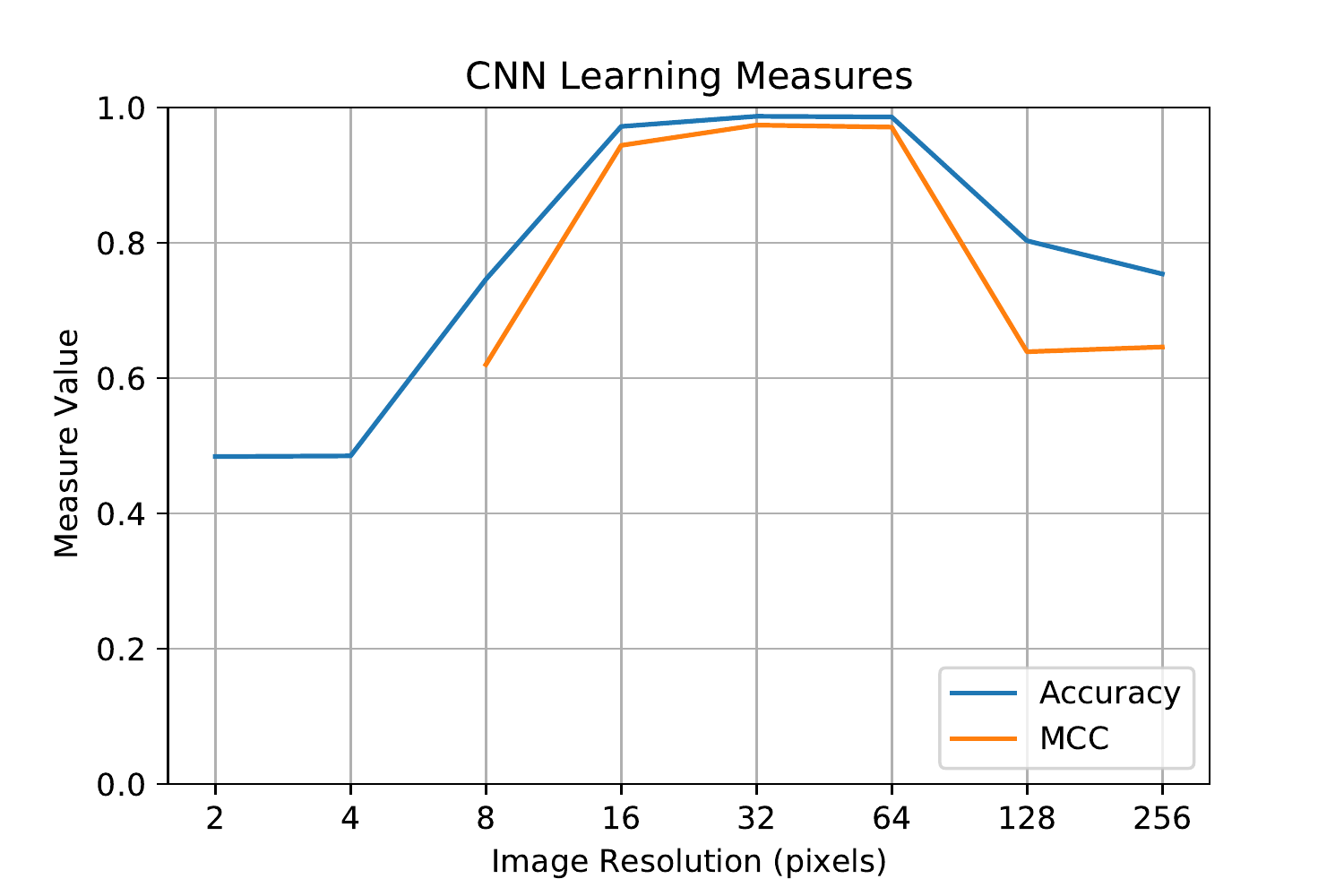}
    \caption{Averaged accuracy and MCC learning measures for the CNNs trained over $F_0$ amoebae images at varying resolutions (denoted by number of pixels in one dimension of the square images). MCC values were incalculable at certain resolutions, and may be artificially high at resolutions $\{8,256\}$.}\label{CNNImage_AccMCC}
\end{figure}

Results for the CNN ML over varying image resolutions are given in Table \ref{ImageMLResults}. The accuracies increase as image resolution improves up to some {\it optimal value}, around $32 \times 32$ pixels where it then falls off; this behaviour is further shown in Figure \ref{CNNImage_AccMCC}.
The high accuracies at optimal resolution ($>0.98$) show the success of CNNs in identifying genus from Monte Carlo generated amoebae images.

The change in accuracy as image resolution decreases from the largest size considered may be explained by initial resolution loss causing averaging over the Monte Carlo generated points to produce a connected amoeba component with a clearer hole structure in genus 1 amoebae (Figure \ref{AmbImages32}).
Then further decrease in resolution loses more of the image information combining any holes into the amoebae themselves until there is no sensible information in the image remaining (Figure \ref{AmbImages2}).
Additionally any poorly sampled parts of the amoeba will be averaged over the resolution decrease to be effectively removed, perhaps causing the amoebae to appear disconnected.

Similar behaviour can be seen over the MCC values where calculable. Although MCC as a measure is less susceptible to bias in the data making it generally preferential, there are issues of incalculability where learning fails and the same class is predicted for all test data.
Where this happened for some of the 5 runs in the 5-fold cross-validation, the MCC average and standard error is calculated over the calculable values, and are denoted with a `*' (note this inflates the measures' value as the cases of no learning are ignored).
Where the MCC is incalculable for all runs an averaged MCC is then incalculable also, these scenarios are denoted with `nan' in the results in Table \ref{ImageMLResults}.

\begin{table}[H]
\begin{tabular}{|c|c|c|c|c|c|c|c|c|}
\hline
\multirow{2}{*}{\begin{tabular}[c]{@{}c@{}}Learning \\ Measures\end{tabular}} & \multicolumn{8}{c|}{Image Resolution} \\ \cline{2-9} & 2$\times$2 & 4$\times$4 & 8$\times$8 & 16$\times$16 & 32$\times$32 & 64$\times$64 & 128$\times$128 & 256$\times$256 \\ \hline
\multirow{2}{*}{Accuracy} & 0.484 & 0.485 & 0.746 & 0.972 & 0.987 & 0.986 & 0.803 & 0.754       \\
& 0.005 & 0.006 & 0.064 & 0.008 & 0.005 & 0.004 & 0.116 & 0.055       \\ \hline
\multirow{2}{*}{MCC} & nan & nan & 0.620* & 0.944 & 0.974 & 0.971 & 0.639 & 0.646*       \\
& nan   & nan   & 0.012* & 0.017 & 0.010 & 0.008 & 0.212 & 0.032*       \\ \hline
\end{tabular}
\caption{Learning results for CNN binary classification of $F_0$ $n=1$ lopsided amoebae images, determining genus 0 or 1 over balanced dataset of 2000 images at varying image resolutions. Resolution is given in terms of the number of pixels. Measures are averaged accuracies and MCCs over the 5-fold cross-validations run, with standard errors. MCCs calculated over less than 5 of the cross-validations are denoted with a `*', for those completely incalculable `nan' is given.}\label{ImageMLResults}
\end{table}

\subsection{Amoeba Image Misclassifications}
\begin{figure}[!tb]
	\centering
	\begin{subfigure}{0.3\textwidth}
		\centering
		\includegraphics[width=\textwidth]{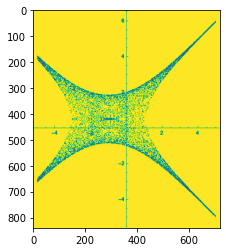}
		\caption{Case 1: $840 \times 720$}\label{g0p1_full}
	\end{subfigure}
	    \begin{subfigure}{0.3\textwidth}
    	\centering
    	\includegraphics[width=0.8\textwidth]{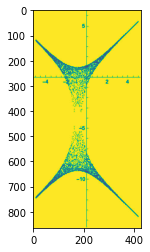}
    	\caption{Case 2: $864 \times 426$}\label{g1p0a_full}
    \end{subfigure}
        \begin{subfigure}{0.3\textwidth}
    	\centering
    	\includegraphics[width=\textwidth]{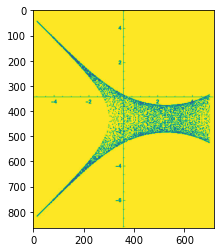}
    	\caption{Case 3: $864 \times 716$}\label{g1p0b_full}
    \end{subfigure}\\
    \begin{subfigure}{0.3\textwidth}
    	\centering
    	\includegraphics[width=\textwidth]{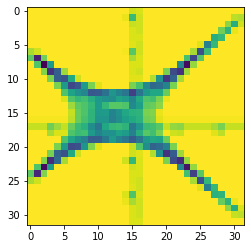}
    	\caption{Case 1: $32 \times 32$}\label{g0p1_32}
    \end{subfigure}
    \begin{subfigure}{0.3\textwidth}
		\centering
		\includegraphics[width=\textwidth]{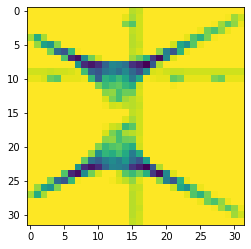}
		\caption{Case 2: $32 \times 32$}\label{g1p0a_32}
	\end{subfigure}
    \begin{subfigure}{0.3\textwidth}
    	\centering
    	\includegraphics[width=\textwidth]{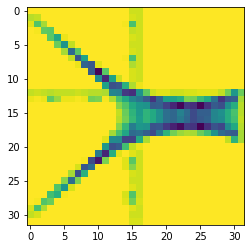}
    	\caption{Case 3: $32 \times 32$}\label{g1p0b_32}
    \end{subfigure}
    \caption{Example CNN misclassifications during model testing, images show the original image and the image at the $32 \times 32$ resolution considered. Case 1, (images a \& d) was labelled as genus 0 with the CNN predicting genus 1; cases 2 $\&$ 3 (images b,e \& c,f) were labelled genus 1 with CNN predictions of genus 0.}\label{AmbImage_misclassifications}
\end{figure}

Examining the image resolution $32 \times 32$ which leads to the most successful CNN classification, the misclassifications where the CNN disagreed with the images labelled genus provide interesting further insight into the learning. 
In one of the 5-fold cross-validation runs 3 of the 400 images in the testing dataset are misclassified such that the predicted amoeba genus differed from the labelled genus.
These images are shown at full resolution as well as the $32 \times 32$ resolution in Figure \ref{AmbImage_misclassifications}.

The first image misclassified, shown in full and lowered resolutions in Figures \ref{g0p1_full} $\&$ \ref{g0p1_32} respectively, is labelled as genus 0, but misclassified to genus 1. 
The original image shows the true genus 0 amoeba, however the Monte Carlo sampling of amoeba points shows a particularly low density in the amoeba's centre, this leads to a poorer collation of points as the resolution drops, shown by fainter parts of the amoebas centre in the lower resolution image perhaps misleading the CNN to predict a non-zero genus.

The second misclassified image, shown in full and lowered resolutions in Figures \ref{g1p0a_full} $\&$ \ref{g1p0a_32} respectively, is labelled as genus 1, but misclassified to genus 0. 
In the original image in particular, the issues with the Monte Carlo image generation become clear, the amoeba looks disconnected into two parts, where the hole in the centre dominates the amoeba.
This leads to part of the hole's boundary being lost in the lower resolution image, making the amoeba appear to be genus 0, misleading the CNN.

The third and final misclassified image in this run, shown in full and lowered resolutions in Figures \ref{g1p0b_full} $\&$ \ref{g1p0b_32} respectively, is also labelled as genus 1, but misclassified to genus 0.
Here the original image shows a genus 0 amoeba, interestingly this scenario is a rare occurrence where the $n=1$ lopsided amoeba (used to produce the genus labels for the dataset as in \ref{gF0}) has a different genus to the true amoeba.
Therefore although this is considered by the CNN as a misclassification, the CNN has actually managed to predict the true amoeba genus from the image, whilst at the same time highlight a case where the match-up to lopsided amoebae fails.

The three misclassification examples from this run demonstrate the subtleties of the Monte Carlo image generation for amoeba, as well as the issues with relying on lopsided amoeba approximation for genus prediction.

\paragraph{Persistent Homology}
Topological data analysis provides an alternative technique for describing the homotopy of data manifolds. 
Within this field, identifying genus of surfaces is often well performed by the technique of persistent homology.

In the persistent homology computation of this image classification problem, the 2-dimensions Monte Carlo generated amoeba points would have 2-dimensional discs of radius $r$ drawn about each of them. 
As the radius is varied through the range: $0 \longmapsto \infty$ the discs begin to intersect, at each new radius where there is a new intersection of $n$ discs the simplicial complex formed from the 0-simplex points has an $n$-simplex drawn between those points (up to the maximum allowed simplex dimensionality: 2).
The discrete updating of the complex as $r$ increases creates a chain complex within which the persistent homology can then be analysed.

To identify the genus, one must analyse the occurrence of 2-dimensional holes; the $H_1$ homology should then have persistent features for each hole contributing to the genus count. 
For these examples where genus is 1 there should hence be a single persistent feature which is far from the persistence diagram diagonal --- corresponding a hole whose boundary is connected quickly and fills in much later.
Conversely for 0 genus there should be no such feature.

Additionally, the $H_0$ homology keeps track of the number of connected components in the complex, for well sampled amoeba points all the features should hence `die' quickly as all the amoeba points become connected to nearby points and form the amoeba as a single component.

\begin{figure}[!tb]
	\centering
	\begin{subfigure}{0.32\textwidth}
		\centering
		\includegraphics[width=1.15\textwidth]{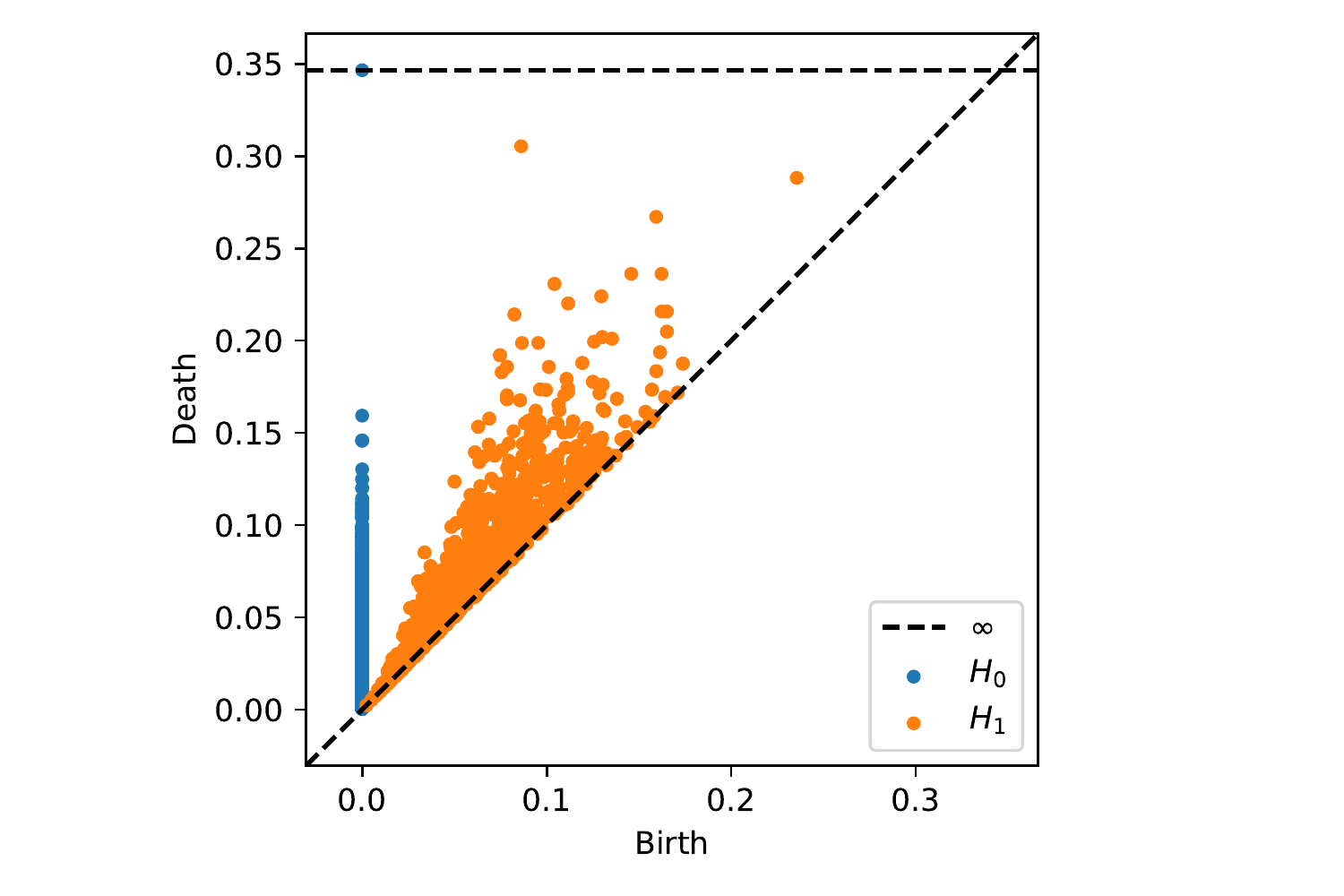}
		\caption{Case 1}\label{g0p1_tda}
	\end{subfigure}
	    \begin{subfigure}{0.32\textwidth}
    	\centering
    	\includegraphics[width=\textwidth]{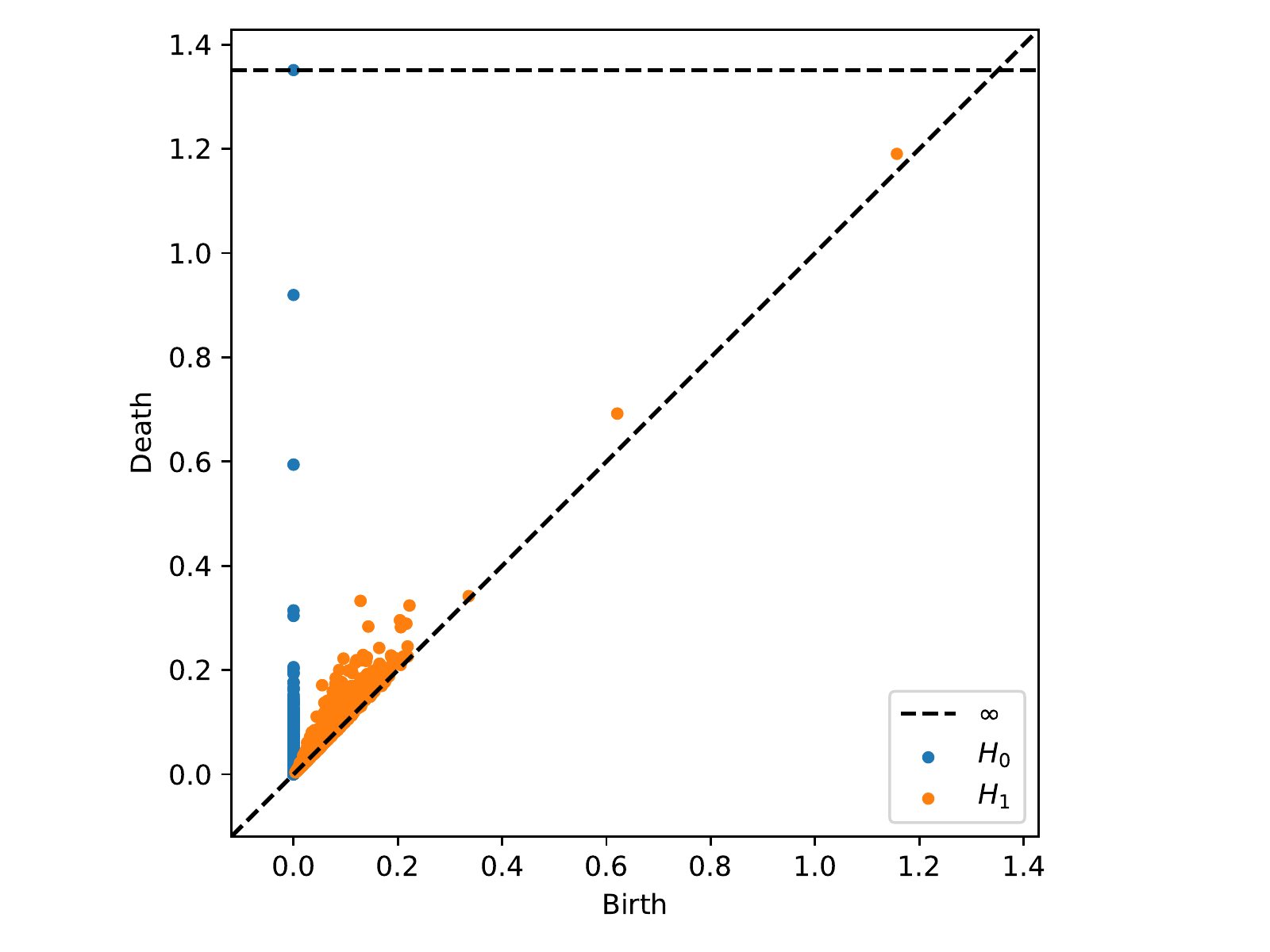}
    	\caption{Case 2}\label{g1p0a_tda}
    \end{subfigure}
        \begin{subfigure}{0.32\textwidth}
    	\centering
    	\includegraphics[width=\textwidth]{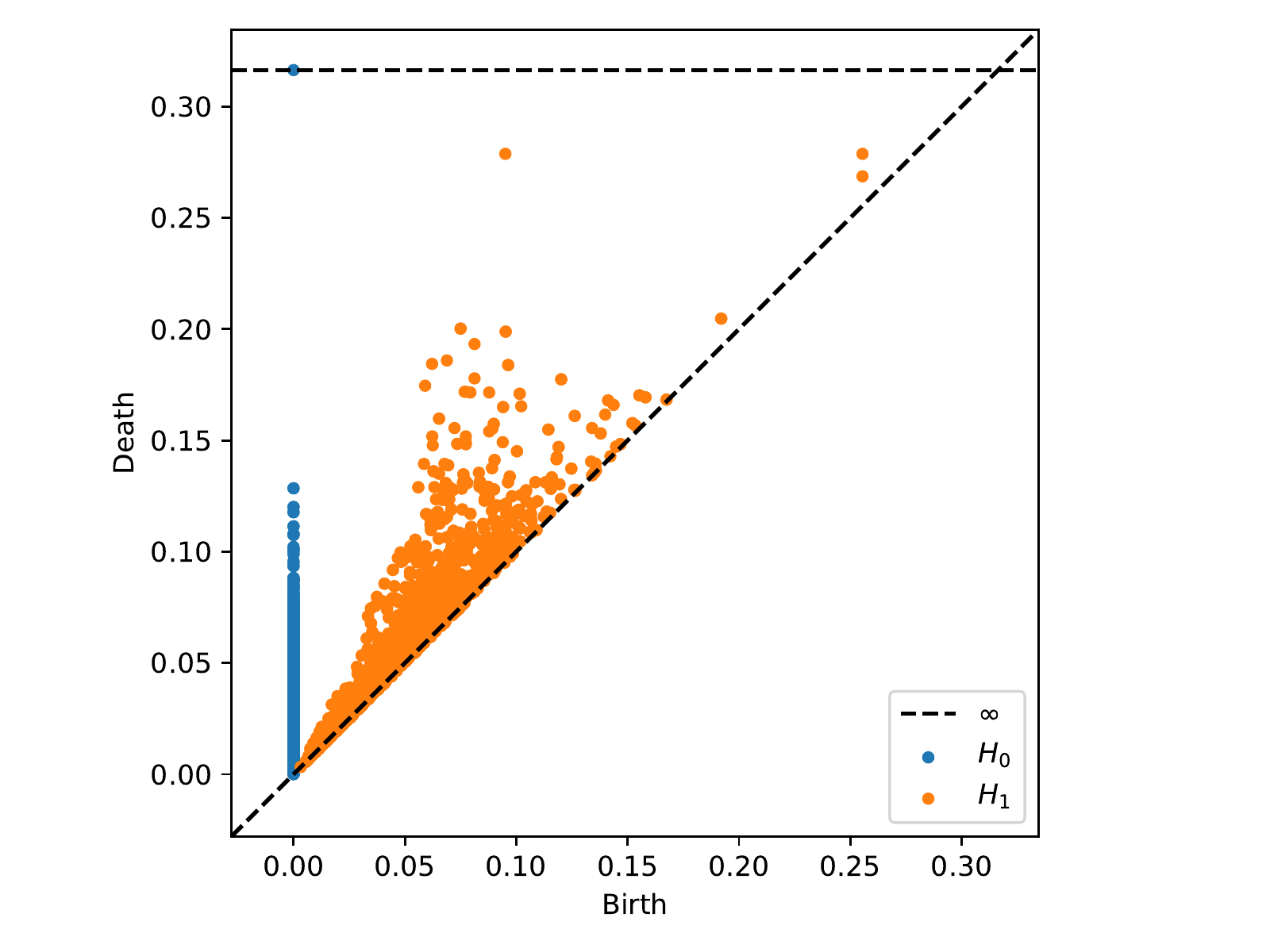}
    	\caption{Case 3}\label{g1p0b_tda}
    \end{subfigure}
    \caption{Persistence diagrams for the 3 misclassification cases considered, showing the $H_0$ and $H_1$ persistent homology of the Monte Carlo sampled amoebae points.}\label{AmbImage_misclassificationsTDA}
\end{figure}

For computation of the persistence diagrams of the 3 misclassification cases considered here the python library \texttt{ripser} was used \cite{ctralie2018ripser}.
The persistence diagrams for each of these cases are plotted in Figure \ref{AmbImage_misclassificationsTDA}.

The $H_0$ homology of cases 1 and 3 show the amoeba is connected quickly, with no gaps in the $H_0$ feature line. 
However for case 2 there are gaps in the line, indicating the sampling is poor and the amoeba appears as disconnected components, agreeing with previous analysis.

For the $H_1$ homology all cases have no significant features far from the diagonal, indicative of a 0 genus amoebae. 
Any significant features should be a distance from the diagonal comparable to the scale of the amoeba.
Since the amoeba centres have size $\sim 2$ and all $H_1$ features die for $r \lesssim 0.3$, despite there being some points slightly further from the diagonal they are not persistent enough to indicate a significant hole and genus 1.
For case 1 the persistent homology now correctly predicts the genus, outperforming the CNN.
However for cases 2 and 3 where the genus is truly 1 the homology prediction is incorrect, likely due to the same reasons as the CNN --- the poor sampling means the case 2 boundary isn't connected properly, and for case 3 the lopsided amoeba boundary is correctly identified even though it does not match the true amoeba.

Both the $H_0$ and $H_1$ data corroborate the CNN results, however the persistent homology performs better with case 1 correctly identifying genus 0.

\section{Conclusions and Outlook}\label{conclusions}
In this paper, we initiated the study of applying neural networks to the analysis of amoebae, especially their genus and membership, which are questions of interest to algebraic geometers and theoretical physicists alike.
We found that a simple MLP or CNN can tell the genus for amoebae with a given Newton polytope with rather high accuracy, from a vector composed of the coefficients of the Newton polynomial $P$. 
In particular, our approach used the amoeba $\mathcal{A}_{P}$ and the lopsided 
amoeba $\mathcal{L}\mathcal{A}_{\tilde{P}_n}$
and showed that ML is able to count the genus with good performance for different $n$'s. 

Importantly, using manifold ML, We showed that one can reproduce the complicated conditions of genus in terms of lopsidedness and the coefficients.
These results are in line with the interpretability of the ML, which is currently one of the biggest challenges in AI. 
In parallel, we saw how one can identify the genus from the amoebae directly as an image processing problem using a CNN.
Here, an optimal image resolution was found which can improve the efficacy on the Monte Carlo generated images; whilst misclassifications highlight subtleties in this process.

As we trained the models with coefficients of the Newton polynomial, the characteristics of the NNs allow one to write down certain approximations to the conditions for deciding the genus. 
It would be interesting to compare this with the approximations from lopsidedness in future work and find more refined approximations using ML.
From our extensive appendices we saw how complicated the decision boundaries of the holes are, and the interpretability of our NN allowed for this approximation.

While we mainly focused on amoebae obtained from the projection onto $\mathbb{R}^2$, one may also take the projection from $(\mathbb{C}^*)^2$ to $\mathbb{T}^2$ and plot the $(\theta,\phi)$ on the real plane. This is known as the \emph{alga} or \emph{coamoeba}, and it is possible to recover brane tilings in string theory therefrom \cite{Feng:2005gw}. A similar study on algae would be an intriguing direction to explore. Besides, we only studied two-dimensional polygons and their associated amoebae in the paper as they are intimately related to D3 branes probing toric CY 3-folds and the mirror geometry involving D6 branes. It is possible to generalize the analysis here to higher dimensions for lattice polytopes and amoebae. More broadly, tropical geometry has been connected to machine learning architecture in \cite{zhang2018tropical,charisopoulos2019tropical}. The piecewise linearity as a common feature for both tropical geometry and neural networks would naturally lead to future communications between the two areas.

\section*{Acknowledgement}
JB is supported by a CSC scholarship. YHH would like to thank STFC for grant ST/J00037X/1. EH would like to thank STFC for the PhD studentship.
We would like to thank Thomas Fink and Forrest Sheldon for stimulating discussion. 

\appendix

\section{Transformations of Amoebae}\label{ap:transformation}

\begin{figure}[H]
	\centering
	(a)\includegraphics[trim=0 0 0 10, clip, width=6cm]{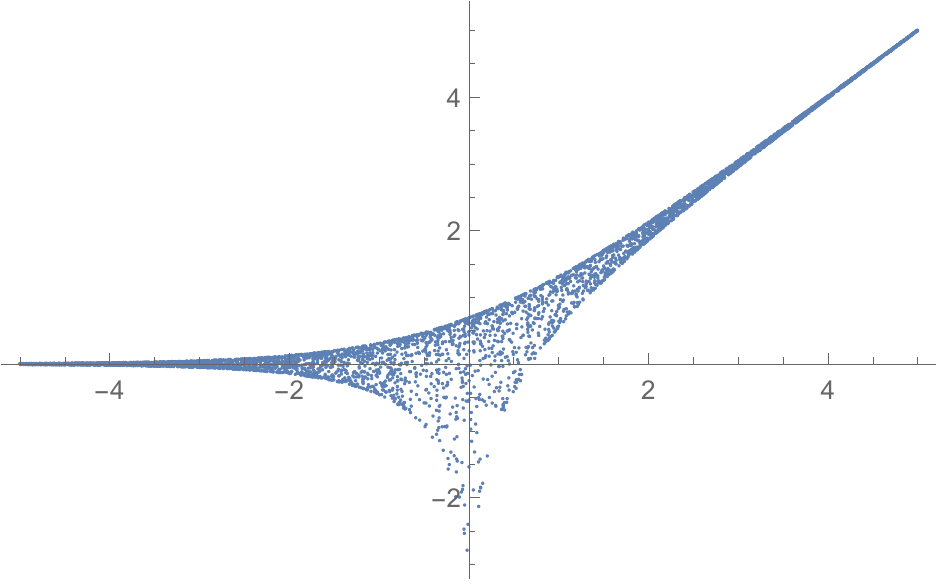}
	(b)\includegraphics[trim=0 0 0 10, clip, width=6cm]{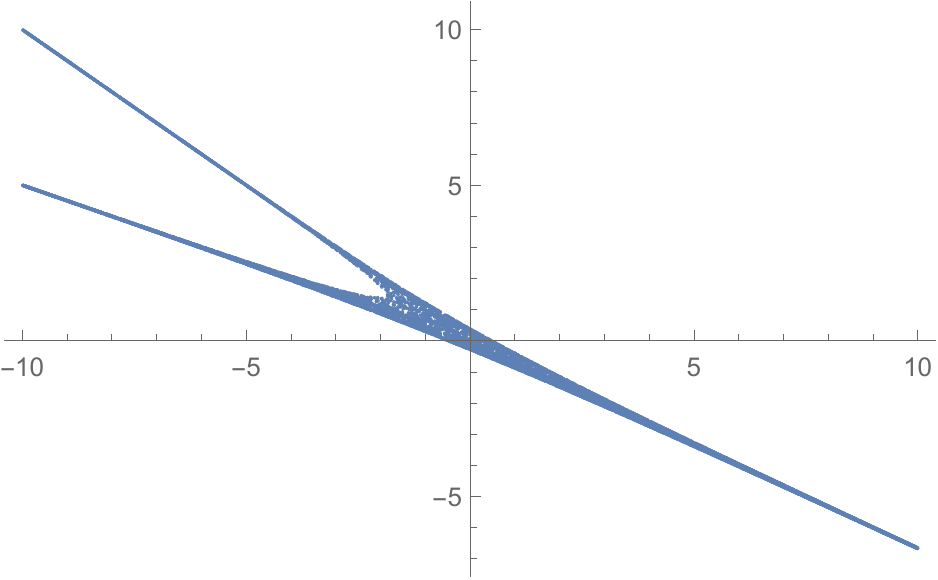}
    \caption{The two amoebae have the same shape. (a) The amoeba given by $z+w+1$. (b) The amoeba given by $z^2w^3+zw^2+1$.}
    \label{isoex}
\end{figure}

When drawing amoebae, it is often much more efficient to use Monte Carlo method than to plot them analytically. However, due to the arguments of $z$ and $w$, the sample points may distribute unevenly on the Log plane. For instance, even for one of the simplest cases, $\mathbb{C}^2$, the amoeba plotted from Monte Carlo method is shown in Figure \ref{isoex}(a).
We can see that one of the tentacles has fewer points compared to the other two. 

As it is often not easy to find suitable distributions for arg$(z)$ and arg$(w)$, one way to resolve this is to consider certain transformations of amoebae \cite{bogaardintroduction}:
\begin{theorem}
	For $(\alpha,M)\in(\mathbb{C}^*)^2\rtimes\textup{GL}(2,\mathbb{Z})$ and $P(\bm{z})\equiv P(z,w)\in\mathbb{C}\left[z,w,z^{-1},w^{-1}\right]$, the map $\Psi:(\mathbb{C}^*)^2\rtimes\textup{GL}(2,\mathbb{Z})\rightarrow\textup{Aut}\left(\mathbb{C}\left[z,w,z^{-1},w^{-1}\right]\right)$ defined by $\Psi(\alpha,M)(P(\bm{z}))=P\left(\alpha\cdot\bm{z}^M\right)$ is an isomorphism. Moreover, their Newton polytopes satisfy $\Delta(\Psi(P))=M\cdot\Delta(P)$.
	Denote the amoeba of $P$ as $\mathcal{A}_P$, then for $\det(M)\neq0$, we have $\mathcal{A}_P=M\mathcal{A}_{\Psi(P)}-\textup{Log}(\alpha)$.
\end{theorem}

Let us take $M=\begin{pmatrix}2&1\\3&2\end{pmatrix}$ and $\alpha=(1,1)$. This sends $z+w+1$ to $z^2w^3+zw^2+1$ whose amoeba is plotted in Figure \ref{isoex}(b). We can see that the points are distributed more evenly. The price is that this would slow the calculations due to the higher degrees.

\section{Manifold Learning}\label{manifoldlearning}
In this appendix, we give a quick introduction to different methods in manifold visualizations, including isomap, t-SNE, MDS and spectral embedding. For more details including other ways of projection, one is referred to the manual of \texttt{Yellowbrick} \cite{bengfort_yellowbrick_2019}.

\paragraph{Isomap} The isometric mapping, also known as isomap, projects the data onto a lower dimensional space while preserving the geodesic distances among all data points. It first applies nearest neighbour search and shortest-path graph search. Then the embedding is encoded in the eigenvectors of the largest eigenvalues of the isomap kernel matrix. See \cite{tenenbaum2000global} for a detailed study on this.

\paragraph{MDS} The multi-dimensional scaling analyzes the similarity of the data, and the data points that are close in the higher dimensional space are also near to each other in the embedding. For the (metric) MDS applied in our main context, it minimizes the cost function known as the stress defined as $\left(\sum\limits_{i,j}(d_{ij}-||x_i-x_j||)^2\right)^{1/2}$ where $d_{ij}$ is the Euclidean distance between the data points $x_i$ and $x_j$. There are also other implementations of MDS. See for example \cite{borg2005modern}.

\paragraph{Spectral embedding} The spectral embedding is a discrete approximation of the low dimensional manifold which uses a graph representation known as the spectral decomposition. The algorithm is known as the Laplacian eigenmaps which transform the data into some graph representation. Then it constructs the graph Laplacian matrix $L=D-A$ where $A$ is the adjacency matrix and $D$ is the degree matrix of the graph. Finally, it applies eigenvalue decomposition on the graph Laplacian. Further discussions can be found in \cite{belkin2003laplacian}.

\paragraph{t-SNE} The t-distributed stochastic neighbour embedding converts the similarity of data points into probabilities. The similarity in the higher dimensional space is represented by Gaussian probabilities while the lower dimensional embedding uses student's t-distributions. Then gradient descent is applied to minimized the relative entropy (aka Kullback-Leibler divergence) of the joint probabilities in the two spaces. In general, t-SNE can group samples and extract local clusters of the data points. Such technique was first proposed and explored in \cite{van2008visualizing}.

\section{Lopsidedness and Cyclic Resultants}\label{lopsided}
In this appendix, we give some examples for lopsidedness and cyclic resultants.

\subsection{Genus for Lopsided Amoeba}\label{lopg}
Using lopsidedness, we can write the conditions for the number of genus for any $\mathcal{LA}_{\Tilde{P}_n}$. We now derive such conditions for some lopsided amoebae where $n=1$.

\paragraph{Example 1: $F_0$} As one of the simplest examples, let us determine the genus of the lopsided amoeba for $F_0$ with $P(z,w)=c_1z+c_2w+c_3z^{-1}+c_4w^{-1}+c_5$ which could have at most one genus corresponding to its sole interior point\footnote{Here, all the coefficients can be any complex numbers. To avoid degenerate cases, we also require $c_{1,2,3,4}\neq0$.}.
Straight away, we can find the centre of the amoeba, which always lie in the hole if $g=1$. We can find the spectral curve (spines) by considering the asymptotic behaviour as
\begin{itemize}
	\item $z,w\rightarrow\infty,z/w\sim\mathcal{O}(1)$: this yields $\text{Log}|w|=\text{Log}|z|-\text{Log}\left|\frac{c_2}{c_1}\right|$.
	
	\item $1/z,1/w\rightarrow\infty,z/w\sim\mathcal{O}(1)$: this yields $\text{Log}|w|=\text{Log}|z|+\text{Log}\left|\frac{c_4}{c_3}\right|$.
	
	\item $z,1/w\rightarrow\infty,zw\sim\mathcal{O}(1)$: this yields $\text{Log}|w|=-\text{Log}|z|-\text{Log}\left|\frac{c_1}{c_4}\right|$.
	
	\item $1/z,w\rightarrow\infty,zw\sim\mathcal{O}(1)$: this yields $\text{Log}|w|=-\text{Log}|z|+\text{Log}\left|\frac{c_3}{c_2}\right|$.
\end{itemize}
In particular, the first two lines are parallel to each other, and so are the other two. The remaining 4 pairs give rise to 4 intersection points (which may or may not coincide). In other words, we have obtained the equations for the four spines of the amoeba and how they surround a rectangle with these 4 intersection points as vertices. 
One can check that this rectangle is centred at $\left(\frac{1}{2}\text{Log}\left|\frac{c_3}{c_1}\right|,\frac{1}{2}\text{Log}\left|\frac{c_4}{c_2}\right|\right)$.

In general, to determine the genus, we should find all possibilities for the lopsided lists. 
Here, using \eqref{lopPoly}, we have that $P\{x_1,x_2\}=\{|c_1z|,|c_2w|,|c_3/z|,|c_4/w|,|c_5|\}$. 
Suppose $|c_5|$ is the largest number, then we have the lopsided condition:
\begin{equation}
	|c_5|>|c_1z|+|c_2w|+|c_3/z|+|c_4/w|.\label{ineqF0}
\end{equation}
However, the right hand side reaches a minimum when $|c_1z|=|c_3/z|,|c_2w|=|c_4/w|$, i.e., $|z|=|c_3/c_1|^{1/2},|w|=|c_4/c_2|^{1/2}$, which is exactly the aforementioned centre of the amoeba. Therefore, $|c_5|$ should at least be greater than this minimum for genus 1, and this bound is precisely $|c_5|>a :=
	2|c_1c_3|^{1/2}+2|c_2c_4|^{1/2}$.
In other words,
\begin{equation}
	g=\begin{cases}
		0,&|c_5|\leq a\\
		1,&\text{otherwise}
	\end{cases}.
\end{equation}
In particular, the centre point is precisely the point where the right hand side of \eqref{ineqF0} reaches its minimum.

For completeness, there are four more possibilities for $P\{x_1,x_2\}$ to be lopsided, but we can see that they would not lead to a non-zero genus by the same argument. For example, suppose the largest is
\begin{equation}
	|c_1z|>|c_2w|+|c_3/z|+|c_4/w|+|c_5|.
\end{equation}
Let us now fix $|w|$, viz, contemplating a horizontal line on the Log plane. If we keep increasing $|z|$ (or equivalently Log$|z|$), this inequality would always hold. Therefore, this region, as a complementary component of the amoeba on the Log plane, would go to infinity. Hence, it is not bounded and cannot be a hole of the amoeba. Likewise, the other three inequalities would not give a hole either by considering the asymptotic behaviour of $|w|$ (or $|z|$) going to infinity or zero while keeping $|z|$ (or $|w|$) fixed. This also verifies that the lopsided amoeba for $F_0$ can have at most genus 1.

\paragraph{Example 2: $L^{3,3,2}$} Let us now consider $L^{3,3,2}$ as it is a non-reflexive polytopes (see \cite{Bao:2020kji}) and hence has more interior points. Its Newton polynomial is $P=c_1z+c_2w+c_3z^{-1}+c_4w^{-1}+c_5z^2+c_6=0$. Therefore, $P\{\bm{x}\}=\{|c_1z|,|c_2w|,|c_3/z|,|c_4/w|,|c_5z^2|,|c_6|\}$. One possibility for these numbers to be lopsided is
\begin{equation}
	|c_1z|>|c_2w|+|c_3/z|+|c_4/w|+|c_5z^2|+|c_6|.
\end{equation}
As $|z|$ cannot be zero, we can divide both sides by $|z|$ and then find the minimum for the right hand side. For $|w|$, it is easy to see that this requires $|w|=w_0\equiv(|c_4/c_2|)^{1/2}$. Then for $|z|$, we have the cubic equation
\begin{equation}
	|c_5||z|^3-\left(2|c_2c_4|^{1/2}+|c_6|\right)|z|-2|c_3|=0.
\end{equation}
Write
\begin{equation}
	\begin{split}
		p&=-\frac{\left(2|c_2c_4|^{1/2}+|c_6|\right)}{|c_5|},~q=-2\left|\frac{c_3}{c_5}\right|,\\
		\Delta&=\left(\frac{q}{2}\right)^2+\left(\frac{p}{3}\right)^3=\left|\frac{c_3}{c_5}\right|^2-\frac{\left(2|c_2c_3|^{1/2}+|c_6|\right)^3}{27|c_5|^3}.
	\end{split}
\end{equation}
Based on the sign of the discriminant, we have three different cases. If $\Delta>0$, or equivalently, $27\left|c_3^2c_5\right|>\left(2|c_2c_4|^{1/2}+|c_6|\right)^3$, then there is only one real root to this equation:
\begin{equation}
	z_0=\sqrt[3]{-\frac{q}{2}+\sqrt{\Delta}}+\sqrt[3]{-\frac{q}{2}-\sqrt{\Delta}}.
\end{equation}
Since $q<0$, $z_0$ is always positive. If $\Delta=0$, there are three real roots. Again, due to negative $q$, we always have a positive root
\begin{equation}
	z_0=-2\sqrt[3]{\frac{q}{2}}.
\end{equation}
If $\Delta<0$, then we would have three distinct roots, $z_{1,2,3}$. Since $z_1+z_2+z_3=0$, there must be at least one positive root, which we shall still call $z_0$. Hence, there would be (at least) one hole if
\begin{equation}
	|c_1|>a_1 := |c_2|w_0/z_0+|c_3|/z_0^2+|c_4|/(z_0w_0)+|c_5|z_0+|c_6|/z_0.
\end{equation}

It is also possible that these numbers are lopsided as
\begin{equation}
	|c_6|>|c_1z|+|c_2w|+|c_3/z|+|c_4/w|+|c_5z^2|.
\end{equation}
Likewise, the right hand side reaches its minimum when $|w|=w_0=(|c_4/c_2|)^{1/2}$ and $|z|=z_0'$ where $z_0'$ is a positive number satisfying\footnote{One can show that there is always a positive root for this cubic equation. A quick way to see this is to consider the function $y=x^2(2|c_5|x+|c_1|)$, which is a cubic curve tangent to the $x$-axis at the origin (and always increasing for positive $x$). It also crosses the negative $x$-axis once while increasing. Then we can simply move this curve down along the $y$-axis to get $y=x^2(2|c_5|x+|c_1|)-|c_3|$. Hence, this would always give a positive root. Using this method, one can also check that both $z_0$ and $z_0'$ give local minima in the two cases.}
\begin{equation}
	2|c_5|z_0'^3+|c_1|z_0'^2-|c_3|=0.
\end{equation}
Then there would be (at least) one hole if
\begin{equation}
	|c_6|>a_2 := |c_1|z_0'+|c_2|w_0+|c_3|/z_0'+|c_4|/w_0+|c_5|z_0'^2.
\end{equation}

One may check that other ways for $P\{\bm{x}\}$ to be lopsided would lead to unbounded complementary regions. To summarize\footnote{If different holes combined with each other, then there would be a point (i.e., fixed $|z|,|w|$) in the hole satisfying more than one inequality. However, this is not possible for fixed $|z|,|w|$.},
\begin{equation}
	g=\begin{cases}
		0,&|c_1|\leq a_1~\text{and}~|c_6|\leq a_2\\
		1,&(|c_1|>a_1~\text{and}~|c_6|\leq a_2)~\text{or}~(|c_1|\leq a_1~\text{and}~|c_6|>a_2)\\
		2,&|c_1|>a_1~\text{and}~|c_6|>a_2
	\end{cases}.
\end{equation}
The punchline is that in this example, we are dealing with cubic (and quadratic) equations. Hence, we can always write down a full analytic condition for the genus. In general, for most of the polygons (even including reflexive ones) as well as $\Tilde{P}_n$, we can always write certain equations to determine the genus for any coefficients, but there may not be general formulae to solve them analytically.

\subsection{Cyclic Resultants for $F_0$}\label{cycresF0}
In this subsection, for reference, we list the cyclic resultants $\Tilde{P}_n$ for $F_0$ for $n=1,2,3$: 
\begin{equation}
    \Tilde{P}_1=P=c_1wz^2+c_2w^2z+c_3w+c_4z+c_5wz,
\end{equation}

\begin{equation}
\begin{split}
    \Tilde{P}_2=&c_2^4 w^8 z^4-2 c_1^2 c_2^2 w^6 z^6-2 c_2^2 c_3^2 w^6 z^2+\left(4 c_4 c_2^3-2 c_5^2 c_2^2-4 c_1 c_3 c_2^2\right) w^6 z^4+c_1^4 w^4 z^8\\
    &+\left(4 c_3 c_1^3-2
   c_5^2 c_1^2-4 c_2 c_4 c_1^2\right) w^4 z^6+\left(4 c_1 c_3^3-2 c_5^2 c_3^2-4 c_2 c_4 c_3^2\right) w^4 z^2\\
   &+\left(c_5^4-4 c_1 c_3 c_5^2-4 c_2 c_4 c_5^2+6
   c_1^2 c_3^2+6 c_2^2 c_4^2-8 c_1 c_2 c_3 c_4\right) w^4 z^4-2 c_1^2 c_4^2 w^2 z^6\\
   &-2 c_3^2 c_4^2 w^2 z^2+\left(4 c_2 c_4^3-2 c_5^2 c_4^2-4 c_1 c_3
   c_4^2\right) w^2 z^4+c_3^4 w^4+c_4^4 z^4,
\end{split}
\end{equation}

\begin{equation}
    \begin{split}
        \Tilde{P}_3=&z^9 c_2^9 w^{18}+3 z^{12} c_1^3 c_2^6 w^{15}+3 z^6 c_2^6 c_3^3 w^{15}+z^9 \left(-9 c_4 c_5 c_2^7+3 c_5^3 c_2^6+18 c_1 c_3 c_5 c_2^6\right) w^{15}\\
        &+3 z^3 c_2^3c_3^6 w^{12}+3 z^{15} c_1^6 c_2^3 w^{12}+z^{12} \left(9 c_2^3 c_3 c_5 c_1^4-21 c_2^3 c_5^3 c_1^3+9 c_2^4 c_4 c_5 c_1^3\right) w^{12}\\
        &+z^6 \left(9 c_3^3 c_4c_5 c_2^4-21 c_3^3 c_5^3 c_2^3+9 c_1 c_3^4 c_5 c_2^3\right) w^{12}+z^9 \left(3 c_4^3 c_2^6-27 c_1 c_3 c_4^2 c_2^5+27 c_4^2 c_5^2 c_2^5\right.\\
        &-18 c_4 c_5^4c_2^4-54 c_1 c_3 c_4 c_5^2 c_2^4+54 c_1^2 c_3^2 c_4 c_2^4+3 c_5^6 c_2^3+9 c_1 c_3 c_5^4 c_2^3-21 c_1^3 c_3^3 c_2^3\\
        &\left.+27 c_1^2 c_3^2 c_5^2 c_2^3\right)w^{12}+z^{18} c_1^9 w^9+c_3^9 w^9+z^{15} \left(-9 c_3 c_5 c_1^7+3 c_5^3 c_1^6+18 c_2 c_4 c_5 c_1^6\right) w^9\\
        &+z^3 \left(-9 c_1 c_5 c_3^7+3 c_5^3 c_3^6+18c_2 c_4 c_5 c_3^6\right) w^9+z^{12} \left(3 c_3^3 c_1^6+27 c_3^2 c_5^2 c_1^5-27 c_2 c_3^2 c_4 c_1^5\right.\\
        &-18 c_3 c_5^4 c_1^4+54 c_2^2 c_3 c_4^2 c_1^4-54 c_2 c_3c_4 c_5^2 c_1^4+3 c_5^6 c_1^3+9 c_2 c_4 c_5^4 c_1^3-21 c_2^3 c_4^3 c_1^3\\
        &\left.+27 c_2^2 c_4^2 c_5^2 c_1^3\right) w^9+z^6 \left(3 c_1^3 c_3^6+27 c_1^2 c_5^2c_3^5-27 c_1^2 c_2 c_4 c_3^5-18 c_1 c_5^4 c_3^4+54 c_1 c_2^2 c_4^2 c_3^4\right.\\
        &\left.-54 c_1 c_2 c_4 c_5^2 c_3^4+3 c_5^6 c_3^3+9 c_2 c_4 c_5^4 c_3^3-21 c_2^3 c_4^3c_3^3+27 c_2^2 c_4^2 c_5^2 c_3^3\right) w^9+z^9 \left(c_5^9-9 c_1 c_3 c_5^7\right.\\
        &-9 c_2 c_4 c_5^7+27 c_1^2 c_3^2 c_5^5+27 c_2^2 c_4^2 c_5^5+27 c_1 c_2 c_3 c_4c_5^5-21 c_1^3 c_3^3 c_5^3-21 c_2^3 c_4^3 c_5^3\\
        &-27 c_1 c_2^2 c_3 c_4^2 c_5^3-27 c_1^2 c_2 c_3^2 c_4 c_5^3-18 c_1^4 c_3^4 c_5-18 c_2^4 c_4^4 c_5+117 c_1c_2^3 c_3 c_4^3 c_5\\
        &\left.-162 c_1^2 c_2^2 c_3^2 c_4^2 c_5+117 c_1^3 c_2 c_3^3 c_4 c_5\right) w^9+3 z^{15} c_1^6 c_4^3 w^6+3 z^3 c_3^6 c_4^3 w^6+z^{12} \left(9c_3 c_4^3 c_5 c_1^4\right.\\
        &\left.-21 c_4^3 c_5^3 c_1^3+9 c_2 c_4^4 c_5 c_1^3\right) w^6+z^6 \left(9 c_1 c_4^3 c_5 c_3^4-21 c_4^3 c_5^3 c_3^3+9 c_2 c_4^4 c_5 c_3^3\right)w^6+z^9 \left(3 c_2^3 c_4^6\right.\\
        &+27 c_2^2 c_5^2 c_4^5-27 c_1 c_2^2 c_3 c_4^5-18 c_2 c_5^4 c_4^4+54 c_1^2 c_2 c_3^2 c_4^4-54 c_1 c_2 c_3 c_5^2 c_4^4+3 c_5^6c_4^3+9 c_1 c_3 c_5^4 c_4^3\\
        &\left.-21 c_1^3 c_3^3 c_4^3+27 c_1^2 c_3^2 c_5^2 c_4^3\right) w^6+3 z^{12} c_1^3 c_4^6 w^3+3 z^6 c_3^3 c_4^6 w^3+z^9 \left(-9 c_2 c_5c_4^7+3 c_5^3 c_4^6\right.\\
        &\left.+18 c_1 c_3 c_5 c_4^6\right) w^3+z^9 c_4^9,
    \end{split}
\end{equation}

\addcontentsline{toc}{section}{References}
\bibliographystyle{utphys}
\bibliography{references}

\end{document}